\documentclass[11pt]{article}
\usepackage{amssymb}
\usepackage{mathrsfs}
\usepackage{latexsym,amsmath,amssymb,amsfonts,epsfig,graphicx,cite,psfrag}
\usepackage{eepic,color,colordvi,amscd}
\usepackage{ebezier}
\usepackage{verbatim}
\usepackage{subcaption}
\usepackage{enumitem}
\usepackage{algorithm}
\usepackage{algorithmic}
\usepackage{hyperref}
\usepackage{amsthm}
\usepackage{longtable}
\usepackage{comment}
\usepackage{color}
\usepackage{xcolor}
\usepackage{graphicx} 
\usepackage{epstopdf}
\usepackage{longtable}
\usepackage{booktabs}
\usepackage{mathtools}
\usepackage{multicol, multirow}
\usepackage{capt-of}
\definecolor{blue}{rgb}{0,0,0.9}
\definecolor{red}{rgb}{0.9,0,0}
\definecolor{green}{rgb}{0,0.9,0}

\theoremstyle{plain}
\newtheorem{remark}{Remark}
\newtheorem{example}{Example}

\newtheorem{assumption}{Assumption}

\newtheorem{proposition}{Proposition}
\newtheorem{theorem}{Theorem}

\newtheorem{lemma}{Lemma}
\usepackage{longtable}

\def\<{\big\langle}
\def\>{\big\rangle}

\def\M{\mathcal{M}}
\def\T{{\rm T}}

\def\DD{{\rm Diag}}

\def\B{\mathcal{B}}

\def\K{\mathcal{K}}
\def\L{\mathcal{L}}
\def\F{\mathcal{F}}
\def\I{\mathcal{I}}

\def\R{\mathbb{R}}

\def\U{\tilde{U}}

\def\Re{{\operatorname{Rtr}}}

\def\g{{\rm grad}}

\def\N{\mathcal{N}}

\def\S{\mathbb{S}}

\def\P{{\rm Proj}}

\def\O{\mathcal{O}}

\def\({\left(}
\def\){\right)}
\def\hR{\widehat{R}}

\def\wt{\widetilde}

\let\svthefootnote\thefootnote
\newcommand\blankfootnote[1]{
	\let\thefootnote\svthefootnote%
}
\usepackage[many]{tcolorbox}    	
\newtcolorbox{boxA}{
    fontupper = \bf,
    boxrule = 1.5pt,
    colframe = black 
}

\def\<{\big\langle}
\def\>{\big\rangle}

\def\M{\mathcal{M}}
\def\T{{\rm T}}

\def\DD{{\rm Diag}}

\def\K{\mathcal{K}}
\def\L{\mathcal{L}}
\def\F{\mathcal{F}}
\def\Z{\mathcal{Z}}
\def\I{\mathcal{I}}

\def\R{\mathbb{R}}

\def\Re{{\operatorname{Rtr}}}

\def\g{{\rm grad}}

\def\N{\mathcal{N}}

\def\S{\mathbb{S}}

\def\P{{\rm Proj}}

\def\O{\mathcal{O}}

\def\({\left(}
\def\){\right)}
\def\hR{\widehat{R}}

\def\wt{\widetilde}

\def\blambda{{{\boldsymbol{\lambda}}}}
\def\hblambda{{{\boldsymbol{\hat\lambda}}}}

\def\bmu{{{\boldsymbol{\mu}}}}

\def\tmu{{{\boldsymbol{\tilde\mu}}}}
\def\hmu{{{\boldsymbol{\hat\mu}}}}
\def\bv{\bar{v}}
\def\v{{v}}
\def\U{{U}}
\def\rv{{\Theta}}
\def\brv{{\bar{\Theta}}}
\def\B{{B}}
\def\Br{\mathcal{B}_r}
\def\p{q}
\def\Mab{\mathcal{M}_r}
\def\RR{{\widehat R\widehat R^\top}}
\DeclareMathOperator{\diag}{diag}
\DeclareMathOperator{\Diag}{Diag}
\def\tol{\tt{tol}}
\def\timelimit{\tt{TimeLimit}}

\begin{document}
        \title{A low-rank augmented Lagrangian method for doubly nonnegative relaxations of mixed-binary quadratic programs}
        
        \author{Di Hou\thanks{Department of Mathematics, National
         University of Singapore, Singapore
         119076 ({\tt dihou@u.nus.edu}).
         }, \quad 
	 Tianyun Tang\thanks{Institute of 
Operations Research and Analytics, National
         University of Singapore, Singapore
         119076 ({\tt ttang@u.nus.edu}).
         }, \quad 
	 Kim-Chuan Toh\thanks{Department of Mathematics, and Institute of 
Operations Research and Analytics, National
         University of Singapore, 
       Singapore
         119076 ({\tt mattohkc@nus.edu.sg}). The research of this author is supported by the Ministry of Education, Singapore, under its Academic Research Fund Tier 3 grant call (MOE-2019-T3-1-010).}
 	}
	\date{\today}
	
\maketitle

\begin{abstract}
Doubly nonnegative (DNN) programming problems are known to be challenging to solve because of their huge number of $\Omega(n^2)$  constraints and $\Omega(n^2)$ variables. In this work, we introduce RNNAL, a method for solving DNN relaxations of large-scale mixed-binary quadratic programs by leveraging their solutions' possible low-rank property. RNNAL is a globally convergent Riemannian augmented Lagrangian method (ALM) that penalizes the nonnegativity and complementarity constraints while preserving all other constraints as an algebraic variety. After applying the low-rank decomposition to the ALM subproblem, its feasible region becomes an algebraic variety with favorable geometric properties. Our low-rank decomposition model is different from the standard Burer-Monteiro (BM) decomposition model in that we make the key improvement to equivalently reformulate most of the quadratic constraints after the BM decomposition into fewer and more manageable affine constraints. This modification is also important in helping us to alleviate the violation of Slater's condition for the primal DNN problem. Moreover, we make the crucial step to show that the metric projection onto the algebraic variety, although non-convex, can be transformed into a solvable convex optimization problem under certain regularity conditions, which can be ensured by a constraint-relaxation strategy. RNNAL is able to handle general semidefinite programming (SDP) with additional polyhedral cone constraints, thus serving as a prototype algorithm for solving general DNN problems. Numerous numerical experiments are conducted to validate the efficiency of the proposed RNNAL method.
\end{abstract}

\bigskip
\noindent{\bf keywords:} semidefinite programming, augmented Lagrangian, 
doubly nonnegative programming, algebraic variety, Riemannian optimization
\\[5pt]
{\bf Mathematics subject classification: 90C20, 90C22, 90C30}

\section{Introduction}

\subsection{Mixed-binary nonconvex quadratic program}

In this paper, we consider the following mixed-binary nonconvex quadratic programs:
\begin{equation}\label{prob-MBQP}
\min\left\{x^{\top} Q x+2 c^{\top}x :
\begin{array}{l}
     Ax=b,\ x_i \in\{0,1\},\ \forall i \in \B,  \\
     x_ix_j=0,\ \forall (i,j)\in E,\ x\in \R^{n}_+ 
\end{array}
\right\},\tag{MBQP}
\end{equation}
where $Q\in \mathbb{S}^{n},\,c\in\R^n,\, A := \left(a_1, \dots, a_m\right)^{\top}\in \R^{m\times n},\,  b := \left(b_1, \dots, b_m\right)^{\top}\in \R^m$, $\B\subseteq[n]$ is the index set of binary variables, and $E\subseteq \{(i,j)\mid 1\leq i<j\leq n \}$ is the index set of incompatible pairs. We assume that $A$ has full row rank and $b\geq 0$, without loss of generality. Problem {\eqref{prob-MBQP} is general} because other problems with partially nonnegative constraints and inequality constraints can be converted into \eqref{prob-MBQP} by splitting free variables and adding slacks to inequalities. \eqref{prob-MBQP} covers various interesting problems such as 0-1 mixed integer programming (MIP), nonconvex quadratic programming (QP), binary integer nonconvex quadratic programming (BIQ), and more.

Since \eqref{prob-MBQP} is in general nonconvex and NP-hard, various convex relaxations have been proposed for finding its global minimizer \cite{bomze2002solving,bundfuss2009adaptive}. In the next subsection, we will describe a convex relaxation that is tractable and frequently used in the literature.

\subsection{Doubly nonnegative relaxation}
\label{subsec-DNN}

In \cite{burer2009copositive,burer2010optimizing}, Burer showed that under a mild assumption, \eqref{prob-MBQP} is equivalent to the following convex linear optimization problem {subject to constraints involving the completely positive cone}:
\begin{equation}\label{prob-CP}
\min\left\{
\<C,Y\>  :\  
Y\in \F_0\cap \Z\cap \mathcal{CP}
\right\},
\end{equation}
where the cost matrix $C  = [0,c^\top;c,Q] 
\in \mathbb{S}^{n+1},$ $\mathcal{CP}$ is the convex cone of completely positive matrices defined by 
$\mathcal{CP} := \operatorname{conv}\{x x^T  :\  x \in \R^{n+1}_+\}$, $\F_0$ and $\Z $ are defined as 
\begin{align*}
    \Z&:=\left\{\begin{pmatrix}z& x^\top \\ x & X\end{pmatrix}\in  \mathbb{S}^{n+1} :\ X_{ij}=0,\ \forall (i,j)\in E \right\},\\
    \F_0 &:= \left\{\begin{pmatrix}1& x^\top \\ x & X\end{pmatrix}\in  \mathbb{S}^{n+1} :\  Ax=b,\ \diag(AXA^\top)=b^2,\ x_i =X_{ii}, \  \forall i \in \B\right\}
\end{align*}
with $b^2 = (b_1^2,\dots,b_m^2)^\top$. Although \eqref{prob-CP} is convex, it is still NP-hard because it has been proven in \cite{dickinson2014computational} that even checking membership in the $\mathcal{CP}$ cone is NP-hard. A practical approach to tackle \eqref{prob-CP} is to relax $\mathcal{CP}$ by the doubly 
nonnegative (DNN) cone $\mathbb{S}_{+}^{n+1}\cap \mathbb{N}^{n+1}$, thus resulting in the following DNN programming problem:
\begin{equation}\label{prob-dnn}
\min\left\{
\<C,Y\>  :\  \
Y\in \F_0\cap \Z\cap \mathbb{S}_{+}^{n+1} \cap\mathbb{N}^{n+1}
\right\},
\end{equation}
where $\mathbb{N}^{n+1}$ denotes the cone of nonnegative matrices in $\R^{(n+1) \times (n+1)}$ and $\mathbb{S}_{+}^{n+1}$ denotes the cone of positive semidefinite matrices
in $\mathbb{S}^{n+1}$. As $\mathcal{CP} \subseteq (\mathbb{S}_{+}^{n+1} \cap \mathbb{N}^{n+1}) $, \eqref{prob-dnn} serves as a lower bound for \eqref{prob-MBQP}. The lower bound provided by the DNN relaxation is usually tight in practice \cite{yoshise2010optimization, kim2016lagrangian, ito2019algorithm} and can be computed using solvers like SDPNAL+ \cite{SDPNALp2}. While these solvers have been successful in solving a variety of DNN problems, they encounter several challenges in solving \eqref{prob-dnn} because of the following two difficulties: 
\begin{enumerate}
\item Convex solvers that directly handle the $n\times n$ matrix variable become inefficient as $n$ increases significantly because the dimensionality of the matrix variable and the number of constraints are of the order $\Omega(n^2)$;
\item The Slater condition for \eqref{prob-dnn} fails, so the strong duality may not hold, and many solvers are unable to produce a solution \cite{tanaka2012application, drusvyatskiy2017many}.
\end{enumerate}

While the first issue is quite obvious, the second issue, which has been shown in \cite{burer2009copositive}, is more subtle. To see it, define
\begin{equation*}
    S = \sum_{i=1}^{m}\begin{pmatrix}
        -b_i\\
        a_i
    \end{pmatrix}\begin{pmatrix}
        -b_i\\
        a_i
    \end{pmatrix} ^\top=\begin{pmatrix}
        -b&
        A
    \end{pmatrix}^\top\begin{pmatrix}
        -b&
        A
    \end{pmatrix}\succeq 0.
\end{equation*}
Then for any feasible solution $Y$ of \eqref{prob-dnn}, we have
\begin{equation*}
    \langle S,Y \rangle=\sum_{i=1}^{m} \begin{pmatrix}
        -b_i\\a_i
    \end{pmatrix}^\top \begin{pmatrix}1& x^\top \\ x & X\end{pmatrix}\begin{pmatrix}
        -b_i\\a_i
    \end{pmatrix}=\sum_{i=1}^{m} \left( b_i^2-2b_ia_i^\top x+a_i^\top X a_i \right)=0,
\end{equation*}
which suggests that $Y$ is not positive definite and hence Slater's condition does not hold. In the next several subsections, we will discuss ideas to alleviate these two issues as well as our contributions.

\subsection{An equivalent formulation of (\ref{prob-dnn})}

The issue of lacking an interior point in \eqref{prob-dnn} has been extensively studied, with two theoretical approaches developed: (i) the facial reduction algorithm (FRA) aims to identify the minimal cone such that the problem restricted to the minimal cone has a strictly feasible point with the same optimal solution as the original problem, see \cite{borwein1981regularizing,borwein1981facial,drusvyatskiy2017many} for general FRA and \cite{bomze2017fresh} for the application of FRA to \eqref{prob-dnn}; (ii) other methods like \cite{ramana1997exact} consider the dual problem and develop an extended dual formulation to ensure strong duality without the need for assuming any constraint qualification. The connection between the two approaches is well-explained in \cite{pataki2013strong}. 

However, implementing FRA can be as challenging as solving the optimization problem itself. Alternatively, various computational approaches have been proposed to reduce FRA's cost, preserve sparsity, and accurately reformulate the original problem \cite{drusvyatskiy2017many,permenter2018partial,zhu2019sieve}. Among them, one effective approach to alleviate the issue is the partial application of facial reduction algorithm (PFRA)  \cite{burer2010optimizing, tanaka2012application}. PFRA restricts the feasible region of problem \eqref{prob-dnn} to a smaller dimensional face exposed by $S$, thus leading to a reformulation with a smaller duality gap and potentially may satisfy the Slater condition. This could help to improve the numerical stability of the solver employed to solve the reformulated problem.

In this paper,  instead of using PFRA, we consider an equivalent reformulation of \eqref{prob-dnn}, which was proposed in \cite{bomze2017fresh} as follows:
\begin{equation}\label{prob-dnn-new}
    \min\left\{
\<C,Y\>  :\  
Y\in \F\cap \Z\cap \mathbb{S}_{+}^{n+1} \cap\mathbb{N}^{n+1}
\right\},
\end{equation}
where $\F$ is defined as
\begin{equation*}
    \F := \left\{\begin{pmatrix}1& x^\top \\ x & X\end{pmatrix}\in  \mathbb{S}^{n+1} :\  Ax=b,\ AX=bx^\top,\ x_i =X_{ii},\  \forall i \in \B\right\}.
\end{equation*}
The new formulation \eqref{prob-dnn-new} is equivalent to \eqref{prob-dnn} in terms of the optimal value and solution. Additionally, \eqref{prob-dnn-new} has several advantages when $E=\emptyset$ \cite{bomze2017fresh,bomze2019notoriously,chen2023semidefinite}: 
\begin{enumerate}
    \item The primal and dual problems of \eqref{prob-dnn-new} are equivalent to those obtained through PFRA, whose primal Slater condition holds when $B$ is an empty set, i.e., no binary variables;
    \item Under some assumptions like the boundedness of the set $\{x\in \R^n_+  :\  Ax=b \}$, the Slater condition for the dual problem of \eqref{prob-dnn-new} holds, thus ensuring attainability of the primal optimal solution and zero duality gap;
    \item When neither the primal nor dual Slater condition holds, \eqref{prob-dnn-new} has the smallest duality gap compared to other reformulations, including \eqref{prob-dnn};
    \item The formulation \eqref{prob-dnn-new} keeps the sparsity structure of the constraint matrices, unlike PFRA, which introduces dense transformation matrices that destroy sparsity.
\end{enumerate}

Even though \eqref{prob-dnn-new} offers 
numerous benefits, to the best of our knowledge, there are currently no methods specifically designed to address it. In contrast, most existing algorithms, such as those described in \cite{burer2010optimizing, tanaka2012application, bomze2019notoriously}, focus on solving \eqref{prob-dnn} through PFRA, rather than exploring its equivalent form \eqref{prob-dnn-new}. In this paper, we propose an efficient algorithm specially designed to tackle \eqref{prob-dnn-new}. While a potential drawback of (\ref{prob-dnn-new}) is the introduction of $\Omega(mn)$ affine constraints $AX=bx^\top$, which may also explain why it has not been explored computationally, we discover that the computational burden can be significantly reduced by an equivalent reformulation of these constraints (under the BM decomposition)
that will be detailed in the following subsection.

\subsection{Low-rank augmented Lagrangian method (ALM)}
Our main question is how to efficiently solve \eqref{prob-dnn-new}. Renowned SDP solvers like SDPT3 \cite{TTT}, SeDuMi \cite{sturm1999using}, and DSDP \cite{benson2008algorithm}, which utilize interior point methods, are rarely used for DNN problems due to their high computational costs per iteration, scaling as $\mathcal{O}(n^6)$. Instead, first-order methods based on the alternating direction method of multipliers (ADMM) \cite{SDPNALp1,chen2017efficient} are preferred for DNN problems. Although solvers like SDPNAL+ have been quite effective in solving medium-sized DNN problems (with $n \leq 2000$), solving large-scale instances (say with $n \geq 3000$) remains highly challenging. This difficulty primarily arises from the costly eigenvalue decompositions required by ADMM-type or augmented Lagrangian
methods to perform projection onto $\S^n_+$, as well as slow convergence issues caused
by degeneracy of  solutions.

In Section~\ref{sec-alg}, in order to reduce the dimension of (\ref{prob-dnn-new}) and avoid expensive spectral decomposition, we will design an algorithm based on low-rank ALM \cite{BM1,BM2}. In detail, in every outer iteration, our algorithm solves the following subproblem:
\begin{equation}\label{ALM-sub}
\min\left\{ \<C,Y\>+\frac{\sigma}{2}\| \Pi_{{(\mathbb{N}^{n+1}\cap \Z)^*} }( \sigma^{-1}W-Y ) \|^2 :\ Y\in \F\cap \S^{n+1}_+\right\},
\end{equation}
where $\sigma>0$ is the penalty parameter, {$(\mathbb{N}^{n+1}\cap \Z)^*$ is the dual cone of $\mathbb{N}^{n+1}\cap \Z$,} and $W\in \S^{n+1}$ is the Lagrangian multiplier of the {nonnegativity} and complementarity constraints $Y\in \mathbb{N}^{n+1}\cap \Z $. Suppose that the subproblem (\ref{ALM-sub}) has an optimal solution of rank $r\in \mathbb{N}^+.$ Then we can apply the BM factorization to get the following equivalent model:
\begin{equation}\label{ALM-sub-LR}
\min\left\{ \<C,\hR \hR^\top\>+\frac{\sigma}{2}\| \Pi_{{(\mathbb{N}^{n+1}\cap \Z)^*} }( \sigma^{-1}W-\hR\hR^\top) \|^2 :\ R\in \N_r\right\},
\end{equation}
where $\hR:=[e_1^\top;R]$, with $e_1$ being the first standard unit vector in $\R^r$ and $\N_r$ is defined as follows:
\begin{equation}\label{algebraic-variety-N}
    \mathcal{N}_r :=  \left\{R\in\R^{n\times r} :\  ARe_1=b,\ ARR^\top=b(Re_1)^\top,\
     \operatorname{diag}_\B(RR^{\top})=R_\B e_1 \right\}.
\end{equation}
We refer the reader to Section~\ref{sec-notations} for the meaning of the notation $\operatorname{diag}_\B(\cdot).$
{Here and in other parts of this paper, given two matrices $P$ and 
$Q$ with the same number of 
columns, the notation $[P;Q]$ denotes the matrix that is obtained by appending Q to the 
last row of $P$.}

In (\ref{ALM-sub-LR}), we deviate from the traditional low-rank ALM approach described in \cite{BM1,BM2,wang2023decomposition}, which penalizes all constraints except possibly the simple diagonal constraints. Instead, we only penalize the {nonnegativity} and complementarity constraints, ensuring that all other constraints are strictly satisfied within the subproblem. Using (\ref{ALM-sub-LR}) has mainly two advantages. First, our numerical experiments indicate that this formulation in (\ref{ALM-sub-LR}) significantly reduces the penalty parameter's magnitude and the number of both outer and inner iterations compared to the traditional low-rank ALM. This efficiency gain stems from preserving more affine constraints in the subproblem, which potentially decreases the number of outer iterations needed for ALM to converge. Second, although the subproblem (\ref{ALM-sub-LR}) is a constrained optimization problem, we find that $\mathcal{N}_r$ has many good geometric properties so that we can preserve the constraints in $\mathcal{N}_r$ with almost negligible computational cost. One important observation is that, although $\N_r$ contains $\Omega(mn)$ quadratic constraints in $ARR^\top=b(Re_1)^\top$, it is equivalent to the following simpler set:
\begin{equation}\label{algebraic-variety-M}
\mathcal{M}_r := \left\{R \in \mathbb{R}^{n \times r} :\  AR=be_1^\top,\ \operatorname{diag}_\B(RR^{\top})=R_\B e_1 \right\},
\end{equation}
which contains only $mr+|\B|$ constraints, among which $mr$ constraints in $AR=be_1^\top$ are linear. Therefore, problem (\ref{ALM-sub-LR}) can be further simplified as follows:
\begin{equation}\label{ALM-sub-LR-M}
\min\left\{ \<C,\hR \hR^\top\>+\frac{\sigma}{2}\| \Pi_{{(\mathbb{N}^{n+1}\cap \Z)^*} } (\sigma^{-1}W- \hR\hR^\top ) \|^2 :\ R\in \M_r\right\}.
\end{equation}
We should emphasize that although (\ref{ALM-sub-LR-M}) is equivalent to (\ref{ALM-sub-LR}), it is necessary to recover the Lagrangian multipliers for the constraints $ARe_1=b$ and $ARR^\top=b(Re_1)^\top$ within $\N_r$ from the KKT solution of (\ref{ALM-sub-LR-M}). This recovery is essential for verifying the optimality conditions of the original convex subproblem (\ref{ALM-sub}) to ensure global optimality. 

In Section~\ref{sec-theorey}, we will present an explicit formula for computing the Lagrangian multipliers of (\ref{ALM-sub}) based on those derived from (\ref{ALM-sub-LR-M}). Additionally, our analysis introduces a rank-adaptive strategy that enables us to escape from non-optimal saddle points and thereby ensuring the convergence of our algorithm to a global optimal solution.

In the next subsection, we will discuss how to solve the ALM subproblem (\ref{ALM-sub-LR-M}) based on Riemannian optimization. 

\subsection{Riemannian optimization on $\M_r$}

One prominent approach for low-rank SDP is the feasible method based on Riemannian optimization \cite{Staircase2,WYmani,tang2024solving,xiao2024dissolving}. This method, however, is traditionally limited to special constraints due to the assumption that the feasible set of the factorized SDP forms a smooth manifold. This methodology was expanded in \cite{tang2023feasible} to encompass SDPs with general constraints, where the feasible sets of the factorized models may not be smooth manifolds.
Nonetheless, feasible methods remain less effective for solving doubly nonnegative (DNN) problems, primarily due to the extensive number of constraints of the order $\Omega(n^2)$.

In addressing low-rank SDP problems with numerous constraints, the Riemannian ALM is a commonly employed approach  \cite{liu2020simple,wang2023solving,wang2023decomposition,zhou2023semismooth}. Recent variants utilizing ALMs with BM factorization include \cite{monteiro2024low,han2024low}. This approach separates the constraints into two categories: one forms a Riemannian manifold and the other is managed through an augmented Lagrangian penalty function. This separation aligns with our approach in (\ref{ALM-sub-LR-M}). However, it is important to note that existing Riemannian ALMs only utilize simple manifolds such as the Cartesian product of unit spheres and Stiefel manifold, which correspond to simple block-diagonal constraints in linear SDP problems. In contrast, the structure of $\M_r$ we consider here is more complex due to the affine constraints $AR=be_1^\top.$ It is currently uncertain whether $\M_r$ qualifies as a manifold, which complicates the execution of operations like projection and retraction around it. 

In Section~\ref{sec-alg-var}, we delve into the geometric properties of $\M_r$. While $\M_r$ may not qualify as a manifold in general, we introduce a constraint-relaxation strategy in Subsection~\ref{subsec-avoid-non-reuglar} to ensure its smoothness. This approach involves the introduction of {slack} variables. Through this transformation, the feasible set of (\ref{ALM-sub-LR-M}) is assured to conform to a manifold structure, thus enabling the application of Riemannian optimization methods for its solution. 

When using the feasible method to solve (\ref{ALM-sub-LR-M}), two important operations are the projection and retraction \cite{manibook,intromani}. The projection onto the tangent space of $\M_r$ involves solving an $(mr+|\B|)$ by $(mr+|\B|)$ positive definite linear system, whose computational cost is in general $\O((mr+|\B|)^3).$ However, in Subsection~\ref{subsec-proj}, we will show that by utilizing the special structure of $\M_r,$ the computational cost of {the} projection can be reduced to $$\O\left(\min\left\{ |\B|^3+m^2r+mr|\B|, (mr)^2|\B|+(mr)^3\right\}\right),$$ which is much smaller than $\O((mr+|\B|)^3)$ when either $|\B|$ or $mr$ is small.

As for retraction, it is typically more complicated than projection because of the nonlinearity and nonconvexity of $\M_r$. In Subsection~\ref{subsec-retrac}, we demonstrate that the non-convex metric projection problem onto $\M_r$ can be equivalently transformed into a convex generalized geometric medium problem. This allows us to adapt the generalized Weiszfeld algorithm to tackle the convex problem with a convergence guarantee. In addition, our analysis is applicable to a broader class of algebraic varieties $\mathcal{M}^g_r$ defined in \eqref{manifold-general}, whereby encompassing the feasible set described in \cite{tang2024feasible} as a special case.

\subsection{Summary of our contributions}

Our paper's contributions are summarized as follows: 
\begin{enumerate}
    \item Unlike existing algorithms such as \cite{burer2010optimizing, tanaka2012application, bomze2019notoriously} that aim to solve DNN relaxations of \eqref{prob-MBQP} using PFRA on \eqref{prob-dnn}, our algorithm focuses on the equivalent form \eqref{prob-dnn-new}. This form retains the sparsity structure of the constraints and has the same smallest duality gap property as PFRA. To the best of our knowledge, there is currently no method specifically designed to solve the DNN relaxations in the form of \eqref{prob-dnn-new}, possibly due to its large number of $\Omega(mn)$ affine constraints.
    
    \item We introduce a Riemannian based augmented Lagrangian method, RNNAL, to solve the DNN problem \eqref{prob-dnn-new}. We design rank-adaptive strategies for escaping from saddle points and develop a technique to recover the dual variables of the DNN problem \eqref{prob-dnn-new}. We also prove the global convergence of RNNAL. Moreover, RNNAL can handle general SDPs with additional polyhedral cone constraints as shown in \eqref{prob-P-dnn-YZ}, thus serving as a prototype algorithm for solving general DNN problems like those in SDPNAL+ \cite{SDPNALp1}.
    
    \item Without the requirement that the feasible set of the DNN problem \eqref{prob-dnn-new} after BM factorization must be a smooth manifold, we avoid the non-smoothness of the algebraic variety $\N_r$ of the ALM subproblem by deriving an equivalent reformulation $\M_r$. We propose a strategy in Subsection \ref{subsec-avoid-non-reuglar} to ensure the smoothness of $\M_r$ via
    reformulating \eqref{prob-MBQP} by
    adding {slack} variables so that its corresponding $\M_r$ is smooth. As far as we know, such a technique has not been employed in the literature.
    
    \item We analyze the smoothness and geometric properties of the algebraic variety $\M_r$. Importantly, we demonstrate that the non-convex retraction problem onto $\M_r$ can be solved via solving a convex generalized geometric medium problem. We adapt the generalized Weiszfeld algorithm to tackle the convex problem and offer theoretical guarantees for its convergence. 
    
    \item We conduct numerous numerical experiments to evaluate the performance of our RNNAL method for solving the DNN relaxations of various classes of MBQP problems. 
    \end{enumerate}

\subsection{Organization}

This paper is organized as follows. In Section~\ref{sec-nota}, we provide some notations and preliminaries. In Section~\ref{sec-alg}, we introduce the augmented Lagrangian framework with the low-rank factorization for solving \eqref{prob-dnn-new}. In Section~\ref{sec-theorey}, we conduct the theoretical analysis of our algorithm RNNAL. Section~\ref{sec-alg-var} analyses the geometric properties of the algebraic variety $\M_r.$ Section \ref{sec-numerical} presents several experiments to demonstrate the efficiency and extensibility of the proposed method. Finally, we conclude the paper in Section \ref{sec-conclusion}.

\section{Notations and preliminaries}\label{sec-nota}

\subsection{Notations}
\label{sec-notations}
Let $\langle A, B\rangle := \operatorname{Tr}\left(A B^{\top}\right)$ denote the matrix inner product and $\|\cdot\|$ be its induced Frobenius norm in $\mathbb{S}^n$. Define $e$ as a column vector of all ones, and $e_1$ as a column vector with 1 as its first entry and zero otherwise. For a matrix $X \in \mathbb{R}^{m \times n}$, $\operatorname{vec}(X)$ denotes the vector in $\mathbb{R}^{mn}$ formed by stacking the columns of $X$. We use $\circ$ to denote the element-wise multiplication operator of two vectors/matrices of the same size. We use $\delta_{\mathcal{C}}(\cdot)$ to denote the indicator function of a set $\mathcal{C}$. Let $[n] := \{1,2, \ldots, n\}$ for any positive integer $n$. For a matrix $X\in\S^{n+1}$, we denote its block decomposition as follows:
\begin{equation}\label{notation-block}
    X=\begin{pmatrix}
        X_{11}&X_{12}\\
        X_{21}&X_{22}
    \end{pmatrix}\in \begin{pmatrix}
        \mathbb{R}&\mathbb{R}^{1\times n}\\
        \mathbb{R}^{n\times 1}&\S^{n}
    \end{pmatrix}.
\end{equation}
Next we define some operators. Given an index set $\B \subseteq[n]$ with its cardinality denoted by $|\B|$, define $\operatorname{diag}_\B: \mathbb{R}^{n \times n} \rightarrow \mathbb{R}^{|\B|}$ such that $\operatorname{diag}_\B(X)=\left(X_{ii}\right)_{i \in \B}$. Its adjoint mapping is denoted as $\operatorname{diag}_\B^*: \mathbb{R}^{|\B|} \rightarrow \mathbb{R}^{n \times n}$, i.e., letting $\B=\{B(1),\cdots,B(t)\}$ for $t=|\B|$, then $\diag_\B^*(\bmu)=\Diag (\tmu)$, where 
\begin{equation*}
    \tmu_i=\begin{cases}
        \bmu_k & \text{if } i=B(k) \text{ for some } k\\
        0     & \text{if } i\not\in \B.
    \end{cases}
\end{equation*}
The index $\B$ is omitted if $\B=[n]$. For a matrix $R \in \mathbb{R}^{n \times r}$, $R_\B \in \mathbb{R}^{|\B| \times r}$ denotes the submatrix of $R$ corresponding to rows in index set $\B$, and $R_i\in\R^{1\times r}$ denotes the $i$-th row of $R$.
Define $\widehat R:= (e^\top_1;R)$ and the linear map $\mathcal{L}_R:\mathbb{R}^{n\times r}\rightarrow \S^{n+1}$ such that \begin{equation}\label{def-Lmap}
    \mathcal{L}_R(H) := 
    \begin{pmatrix}
        0&e_1^\top H^\top\\
        He_1&HR^\top+RH^\top
    \end{pmatrix}\quad \forall H\in \mathbb{R}^{n\times r}.
\end{equation}
The above linear map is the differential of the factorization $R\rightarrow \widehat R\widehat R^\top$. The adjoint map $\mathcal{L}^*_R:\S^{n+1}\rightarrow \mathbb{R}^{n\times r}$ is given by
\begin{equation}\label{def-Lmap*}
    \mathcal{L}_R^*(S) := 
    2[0_{n\times 1},I_n]S\widehat R \;=\; 2(S_{21} e_1^\top + S_{22} R).
\end{equation}

\subsection{Preliminiaries on Riemannian optimization}
In this subsection, we provide basic material on Riemannian optimization. Although this area concerns optimization problems over general Riemannian manifolds, we use the problem (\ref{ALM-sub-LR-M}) as a special example to illustrate the main ideas. We rewrite the problem (\ref{ALM-sub-LR-M}) as follows
\begin{equation}\label{Roptprob}
\min\left\{ f(R):\ R\in \M_r \right\},
\end{equation}
where $r\in \mathbb{N}^+,$ $\M_r\neq \emptyset$ is defined as in (\ref{algebraic-variety-M}) and $f:\R^{n\times r}\rightarrow \R$ is a continuously differentiable function. A point $R$ in the feasible set $\M_r$ is considered regular (or smooth) if it satisfies the linear independence constraint qualification (LICQ) condition. We assume that the LICQ condition holds everywhere inside $\M_r$ so that $\M_r$ is a Riemannian manifold embedded in $\mathbb{R}^{n\times r}$, with the metric being the Euclidean metric.  We will discuss how to ensure the LICQ property in Subsection~\ref{subsec-avoid-non-reuglar} later. For every point $R\in \M_r,$ its {\bf tangent space} is defined as follows:
\begin{equation}\label{pre-tansp}
\T_R\M_r:=\left\{H \in \mathbb{R}^{n\times r} :\  AH=0,\ 2\operatorname{diag}_\B(H R^{\top})-H_\B e_1=0 \right\},
\end{equation}
which is the set of directions of smooth curves in $\M_r$ passing through $R.$ The {\bf projection} $\P_{\T_R\M_r}:\R^{n\times r}\rightarrow \T_R\M_r$ is the orthogonal projection onto the tangent space $\T_R\M_r.$ The {\bf Riemannian gradient} $\g f(R)$ is the projection of the Euclidean gradient $\nabla f(R)$ onto $\T_R\M_r,$ i.e., $\P_{\T_R\M_r} (\nabla f(R)).$ A {\bf retraction} is a smooth map $\Re_R:\T_R \M_r\rightarrow \M_r$  that satisfies
\begin{itemize}
\item[(i)] $\Re_R(0_{n\times r})=R .$
\item[(ii)] ${\rm D}\Re_R(0_{n\times r})[H]=H,$ for any $H\in \T_R \M_r$,
\end{itemize}
where ${\rm D}\Re_R(0_{n\times r})[\cdot]$ denotes the Frechét differential mapping of $\Re_R$ at $0_{n\times r}$. Note that the above two conditions essentially say that $\Re_R(\cdot)$ is a first order approximation of the exponential mapping of $\M_r$ at $R.$ We call such a retraction a {\bf first order retraction}. Moreover, if for any $H\in \T_R\M_r,$ 
\begin{equation}\label{pre-2nd-retrac}
\Re_R\(tH\)=R+tH+\frac{t^2}{2}V+o(t^2)
\end{equation}
for some $V\in \(\T_R\M_r\)^\perp,$ then we call $\Re_R$ a {\bf second order retraction}. A commonly used second order retraction is the metric projection $\P_{\M_r}(\cdot)$ (see \cite[Section 5.12]{intromani}). With the availability of Riemannian gradient and retraction, the Riemannian gradient descent method, which will be used in our algorithm, updates the iterate $R$ as follows:
\begin{equation}\label{pre-riegd}
R^+:=\Re_R\(-t\cdot\g f(R)\),
\end{equation}
where $t>0$ is some stepsize. By Taylor expansion, we have that
\begin{equation}\label{pre-Tyex}
f(R^+)=f(R)-t\| \g f(R) \|^2+o(t),
\end{equation}
which implies that a decrease in the function value is guaranteed as long as $\| \g f(R) \|\neq 0$ and $t>0$ {is} small enough.

The above is a simple introduction to Riemannian optimization. For more information on this topic, we refer the reader to books such as \cite{manibook, intromani}.

\section{Algorithm}\label{sec-alg}
In this section, we design an augmented Lagrangian method with the low-rank factorization to solve \eqref{prob-dnn-new} based on Riemannian optimization. Consider the more general problem:
\begin{equation}\label{prob-P-dnn-YZ}
    \min\left\{
\<C,Y\> :\ 
Y-Z=0,\ Y\in \F\cap\mathcal{K},\ Z\in  \mathcal{P}
\right\},\tag{P}
\end{equation}
where $\mathcal{K}=\mathcal{K}^*=\mathbb{S}_{+}^{n+1}$ is the positive semidefinite matrix cone, $\mathcal{P}$ is a polyhedral convex cone in $\mathbb{S}^{n+1}$. This polyhedral convex cone $\mathcal{P}$ includes the cone of symmetric matrices with non-negative and zero entries, denoted as $\mathbb{N}^{n+1}\cap \Z$ in equation \eqref{prob-dnn-new}, as a special instance. We remark that our algorithm can also be directly extended to solve \eqref{prob-P-dnn-YZ} with a more general closed convex set $\mathcal{P}$. We assume that $\F$ has the compact form:
\begin{equation*}
    \F=\left\{Y\in\S^{n+1} :\  \mathcal{A}(Y)= {d}  \right\}.
\end{equation*} 
For \eqref{prob-dnn-new}, the operator $\mathcal{A} :\S^{n+1}\rightarrow \R^{m+mn+|\B|+1}$ and $ {d} \in\R^{m+mn+|\B|+1}$ are defined as
\begin{equation*}
    \mathcal{A}(Y) := \begin{pmatrix}
        AY_{21}\\
        \operatorname{vec}(AY_{22}-bY_{12})\\
        \diag_\B (Y_{22})-(Y_{21})_\B \\
        Y_{11}
    \end{pmatrix},\quad  {d}  := \begin{pmatrix}
        b\\
        0_{mn}\\
        0_{|\B|}\\
        1
    \end{pmatrix},
\end{equation*}
where $Y_{11},Y_{12},Y_{21},Y_{22}$ are the submatrices of $Y$ defined in \eqref{notation-block}. The optimality conditions (KKT conditions) for \eqref{prob-P-dnn-YZ} can be written as follows:
\begin{equation}\label{KKT-DNN}
    \begin{aligned}
        & Y-Z=0,\ \mathcal{A}(Y)- {d} =0,\ C-\mathcal{A}^*(y)-S-W=0, \\
        &\langle Y,{S}\rangle=0,\ Y\in \mathcal{K},\ S\in \K^*,\ \langle Z, W\rangle=0,\ Z\in\mathcal{P},\ W\in \mathcal{P}^*,
    \end{aligned}
\end{equation}
where $y\in \R^{m+mn+|\B|+1} $, $W\in\mathbb{S}^{n+1}$ and $S\in\K^*$ are dual variables. We make the following assumption throughout the paper.
\begin{assumption}\label{assumption-mild}
    The problem \eqref{prob-P-dnn-YZ} admits an optimal solution, and its objective function is bounded from below.
\end{assumption}

\subsection{Augmented Lagrangian method}
We first present the augmented Lagrangian method for solving \eqref{prob-P-dnn-YZ}. Define $\mathcal{M} := \F\cap \mathcal{K}$, then \eqref{prob-P-dnn-YZ} can be equivalently written as
\begin{equation}\label{prob-ALM-reform}
    \min\left\{\<C,Y\> + \delta_\mathcal{M}(Y)+\delta_{\mathcal{P}}(Z) :\  Y-Z=0\right\}.
\end{equation}
Let $\sigma > 0$ be a given penalty parameter. The augmented Lagrange function is defined by
\begin{align*}
    L_\sigma(Y,Z;W) :&= \<C,Y\> + \delta_\mathcal{M}(Y)+\delta_{\mathcal{P}}(Z) - \<W,Y-Z\> + \frac{\sigma}{2}\left\|Y-Z\right\|^2\\
    &=\<C,Y\> + \delta_\mathcal{M}(Y)+\delta_{\mathcal{P}}(Z) + \frac{\sigma}{2}\left\|Y-Z-\sigma^{-1}W\right\|^2-\frac{1}{2\sigma}\|W\|^2.
\end{align*}
We can apply the following augmented Lagrangian method to solve \eqref{prob-ALM-reform}. Specifically, given $\sigma_0>0$, $W^0\in \mathcal{P}^*$, perform the following steps at the $(k+1)$-th iteration:
\begin{align}
    (Y^{k+1},Z^{k+1})&=\arg\min\ \{L_{\sigma_k}(Y,Z;W^k) :\ Y\in \M,\ Z\in \mathcal{P}\},\label{alg-alm-prim}\\
    W^{k+1}&=W^k-\sigma_k(Y^{k+1}-Z^{k+1}),\notag
\end{align}
where $\sigma_k \uparrow \sigma_{\infty} \leq+\infty$ are positive penalty parameters. For a general discussion on the augmented Lagrangian method for solving convex optimization problems and beyond, see \cite{hestenes1969multiplier,powell1969method,rockafellar1976augmented}. Let $\widetilde{W}\in \mathbb{S}^{n+1}$ be fixed. The inner subproblem \eqref{alg-alm-prim} can be expressed as:
\begin{equation}\label{prob-alm-prim-v1}
    \min\left\{L_\sigma(Y,Z;\widetilde W)  :\  Y\in \M,\ Z\in \mathcal{P}\right\}.
\end{equation}
In \eqref{prob-alm-prim-v1}, we can first minimize with respective to $Z\in \mathcal{P}$ to get the following convex optimization problem related only to $Y$:
\begin{equation}\label{prob-alm-prim-v2}
    \min\left\{\phi(Y) := \<C,Y\>+\frac{\sigma}{2}\| \Pi_{\mathcal{P}^*}(\sigma^{-1} \wt W-Y) \|^2 :\  Y\in \mathcal{M} \right\},\tag{ALM-sub}
\end{equation}
where we use the Moreau decomposition theorem in \cite{moreau1962decomposition}, which states that $X=\Pi_{\mathcal{C}}(X)-\Pi_{\mathcal{C}^*}(-X)$ for any $X \in \S^{n+1}$ and any closed convex cone $\mathcal{C} \subseteq \S^{n+1}$. Once we obtain the optimal solution $\wt Y$ of \eqref{prob-alm-prim-v2}, we can recover the optimal solution $\wt Z=\Pi_P(\wt Y-\sigma^{-1}\wt W)$. The ALM framework for solving \eqref{prob-ALM-reform} is summarized in Algorithm \ref{alg1}, where $Y^{k+1}$ is obtained by the factorization described in the next subsection. 
\begin{algorithm}
\linespread{1.1}\selectfont
\caption{The RNNAL method}
\label{alg1}
\begin{algorithmic}[1]
\STATE {\bf Parameters:} Given $\sigma_0>0$, initial point $R^0\in\M_r$.
\STATE $k\gets 0$, $W^0=0$.
\WHILE{not converged}
    \STATE Obtain $R^{k+1}$ by solving \eqref{prob-Rie-sub} inexactly.
    \STATE $Y^{k+1}=\widehat R^{k+1}(\widehat R^{k+1})^T$.
    \STATE $Z^{k+1}=\Pi_{\mathcal{P}}( Y^{k+1}-\sigma_k^{-1} W^{k})$.
    \STATE $W^{k+1}=W^k-\sigma_k(Y^{k+1}-Z^{k+1})$.
    \STATE {Update $\sigma_k$}.
    \STATE $k\gets k+1$.
\ENDWHILE
\end{algorithmic}
\end{algorithm}

\subsection{The Burer-Monteiro factorization approach for solving 
\eqref{prob-alm-prim-v2}}
In this subsection, we discuss how to solve the ALM subproblem \eqref{prob-alm-prim-v2} by utilizing the BM factorization and Riemannian optimization algorithm. Assume that \eqref{prob-alm-prim-v2} possesses an optimal solution with rank no greater than $r$. Note that any optimal solution $Y\in \mathcal{M}$ with rank $r$ can be factorized as
\begin{equation*}
    Y=\begin{pmatrix}
        1 & x^{\top} \\
x & X
    \end{pmatrix}=\begin{pmatrix}
        e_1^{\top} \\
    R
    \end{pmatrix}
    \begin{pmatrix}
        e_1 & R^{\top}
    \end{pmatrix}=\RR,
\end{equation*}
where $R\in \R^{n\times r}$ and $e_1\in \R^{r}$. Thus, \eqref{prob-alm-prim-v2} corresponding to \eqref{ALM-sub} is equivalent to the following problem:
\begin{equation}\label{prob-sub-v1}
    \min\left\{
    f_r(R) := \phi(
    \RR
    )
      :\ 
    R\in \N_r
     \right\},
\end{equation}
where $\N_r $ is defined as in (\ref{algebraic-variety-N}). One common attempt to solve \eqref{prob-sub-v1} is by Riemannian optimization algorithms. However, the following lemma shows that any point $R\in\N_r$ does not satisfy the LICQ condition, which serves as the key assumption for many Riemannian optimization algorithms. 
\begin{lemma}\label{lemma-non-licq}
    For any $r>0$ and $R\in \N_r$, the LICQ condition of $\N_r$ at $R$ does not hold.
\end{lemma}

\proof{}
    For any $r>0$ and $R\in \N_r$, note that the first two types of constraints in $\mathcal{N}_r$ suggest that $(AR)(AR)^\top=bb^\top$, which indicates that $\|a_i^\top R\|=b_i,\ \forall i\in [n]$. Also note that $ARe_1=b$ indicates that $a_i^\top Re_1=b_i,\ \forall i\in [n]$. Combining the two equations, we must have $AR=be_1^\top$. The Frechét differential mapping of the constraints defining $\mathcal{N}_r$ at $R$ is given by
    \begin{align*}
    g_R(H):&=(AHe_1; ARH^\top+AHR^\top-be_1^\top H^\top; 2\diag_\B (HR^\top)-H_\B e_1)
     \\
    &=(AHe_1;AHR^\top;2\diag_\B (HR^\top)-H_\B e_1)
     \quad\quad\forall H \in \mathbb{R}^{n \times r},
    \end{align*}
    where we used the fact that $AR=be_1^\top$. The adjoint mapping of $g_R(\cdot)$ is
    \begin{equation*}
    g_R^*(\blambda_1,\blambda_2,\bmu) := A^\top\blambda_1 e_1^\top + A^\top\blambda_2 R +\diag_\B^*(\bmu)(2R-ee_1^\top)\quad \forall (\blambda_1,\blambda_2,\bmu)\in \R^{m}\times \R^{m\times n}\times \R^{|\B|}.
    \end{equation*}
    Since $(\blambda_1,\blambda_2,\bmu)=(b,-A,0)$ is a nonzero solution to the equation $g^*_R(\blambda_1,\blambda_2,\bmu)=0$, the operator $g_R^*(\blambda_1,\blambda_2,\bmu)$ is not injective, and hence the LICQ condition does not hold. 
\endproof

To ensure the LICQ condition for the problem \eqref{prob-sub-v1}, we convert the set $\mathcal{N}_r$ to the set $\mathcal{M}_r$ as defined in \eqref{algebraic-variety-M}, and consider the new subproblem:
\begin{equation}\label{prob-Rie-sub}
    \min \left\{f_r(R) :\  R \in \mathcal{M}_r\right\}.\tag{Rie-sub}
    \end{equation}
From the proof of Lemma \ref{lemma-non-licq}, we have $\M_r=\N_r$. Thus, \eqref{prob-Rie-sub} is equivalent to \eqref{prob-sub-v1}. We should note that while the LICQ condition does not hold at any $R\in\N_r$, the
situation for $\M_r$ is much better as we can see later in Section~\ref{sec-alg-var}.
In particular, the LICQ condition holds for all $R\in\M_r$ when the index set $\B=\emptyset$. For later usage, we note that the Fr\'{e}chet differential mapping of the constraints defining $\M_r$ at $R$ is given by
$h_R:\mathbb{R}^{n \times r}\rightarrow  \mathbb{R}^{m\times r}\times\mathbb{R}^{|\B|}$ 
such that 
\begin{eqnarray}
h_R(H) :=  (AH ; 2\operatorname{diag}_\B (H R^{\top})-H_\B  e_1).
\label{eq-hR}
\end{eqnarray}

We assume that the following condition holds.
\begin{assumption}\label{assumption-LICQ}
    There exists a positive rank bound $\bar r\in \mathbb{N}^+$ such that for any $r\geq \bar r$, the set $\M_r$ is non-empty and satisfies LICQ for every $R\in \M_r$.
\end{assumption}
While the LICQ property may not hold {for some points in} the set $\M_r$ in general, in Section \ref{sec-alg-var}, we will analyze the smoothness of the new set $\M_r$, and provide strategies to ensure Assumption \ref{assumption-LICQ}. Under the assumption, we can use the Riemannian gradient descent method with Barzilai-Borwein stepsize {\cite{iannazzo2018riemannian}} to solve \eqref{prob-Rie-sub}. The algorithm framework is presented in Algorithm \ref{alg2}.

\begin{algorithm}[h]
\linespread{1.1}\selectfont
\caption{Riemmannian gradient descent for \eqref{prob-Rie-sub}}
\label{alg2}
\begin{algorithmic}[1]
\STATE {\bf Parameters:} $\epsilon_g>0$, $\epsilon_H>0$, $\tau>0$, $\bar r\in \mathbb{N}^+$, $\{\epsilon_i\}_{i \geq 0} \subset \mathbb{R}_{+}$, initial point $R_0\in\M_r$.
\STATE $i\gets 0$,\ $r_0=\operatorname{rank}(R_0)$.
\WHILE{not converged}
    \IF{$\|\operatorname{grad} f_{r_i}(R_i)\|> \epsilon_g$}
        \STATE Obtain $R_{i+1}$ by the Riemannian gradient descent method.
        \STATE $r^+_{i}=r_i$.
    \ELSE
        \STATE Recover the dual variable $S_{i+1}$ by Theorem \ref{thm-recover}.
        \IF{$\lambda_{\min}(S_{i+1})<-\epsilon_H$}
            \STATE Obtain $R_{i+1}$ by escaping from a saddle point by Theorem \ref{thm-second-order-descent}.
            \STATE $r^+_{i}=r_i+\tau$.
        \ENDIF
    \ENDIF
    \STATE Find $R'\in\M_{r'}$ such that $\bar r\leq r'\leq r^+_{i}$ and 
    $f_{r'}(R')\leq\, f_{r^+_{i}}(R_{i+1})+\epsilon_i$. \COMMENT{reduce rank}
    \STATE $R_{i+1}\gets R'$,\ $r_{i+1}=r'$.
    \STATE $i\gets i+1$.
\ENDWHILE
\end{algorithmic}
\end{algorithm}

\section{Theoretical Analysis}\label{sec-theorey}
In this section, we provide theoretical guarantees to ensure that the non-smooth non-convex subproblem \eqref{prob-Rie-sub} can be solved to global optimality. We also establish the global convergence of the ALM framework.

\subsection{Recovering dual variables}
When applying our algorithm to solve the DNN problem \eqref{prob-dnn-new}, we must check the global optimality of both the DNN problem \eqref{prob-P-dnn-YZ} and the ALM subproblem \eqref{alg-alm-prim} by their respective KKT conditions. However, Algorithm \ref{alg1} only provides the primal variables $Y^k,Z^k$ and dual variable $W^k$, so we have to recover the remaining dual variables. First, we recover the dual variables
corresponding to $\M_r$ in \eqref{prob-Rie-sub}. Let $\sigma>0$ and $ \widetilde W\in \S^{n+1}$ be fixed. Define $W := \sigma\Pi_{\mathcal{P}^*}(\sigma^{-1}\widetilde W- Y)$ and two auxiliary variables:
    \begin{align}
    L &:= Q-\operatorname{diag}_\B^*({\bmu})-W_{22},\qquad
    \p := 2 c+\tmu_\B -2 W_{21}, \label{def-L}
    \end{align}
where $\tmu_\B := \diag_\B^*(\bmu)e$. Then the KKT conditions for the subproblem \eqref{prob-Rie-sub} are
\begin{subequations}
\label{KKT-SDPR}
\begin{align}
& 2LR+\p e_1^\top-A^{\top} {\blambda}=0, \label{KKT-SDPR-line1}\\
& A R=b e_1^{\top},\ \operatorname{diag}_\B (R R^{\top})-R_\B  e_1=0,\label{KKT-SDPR-line2}
\end{align}
\end{subequations}
where $({\blambda},{\bmu}) \in \mathbb{R}^{m \times r}\times \mathbb{R}^{|\B|}$ are the only unknown dual variables in the KKT system. We may solve \eqref{KKT-SDPR-line1} to get a least square solution of $({\blambda},{\bmu})$. When $R\in\M_r$ is regular, the solution $({\blambda},{\bmu})$ of \eqref{KKT-SDPR-line1} is unique. Moreover, the cost for recovering $({\blambda},{\bmu})$ can be ignored because $({\blambda},{\bmu})$ has been computed in the Riemmanian gradient in \eqref{eq-grad}. 

Next, we recover the dual variables corresponding to $\N_r$ in \eqref{prob-sub-v1}, \eqref{prob-alm-prim-v2} and \eqref{prob-P-dnn-YZ} from the KKT solutions of \eqref{prob-Rie-sub}. The results are summarized in Theorem \ref{thm-recover}. Before that, we write down the KKT conditions of these problems. The KKT conditions for \eqref{prob-sub-v1} are
\begin{equation}
\label{KKT-sub-R}
\begin{aligned}
&2LR+\p e_1^\top-(A^{\top} {\blambda_1}e_1^\top+A^\top\blambda_2 R+\blambda_2^\top A R-\blambda_2^\top be_1^\top)=0, \\
&ARe_1=b,\ ARR^\top=b(Re_1)^\top,\
     \operatorname{diag}_\B (RR^{\top})=R_\B  e_1, 
\end{aligned}
\end{equation}
where $(\blambda_1,\blambda_2,\bmu)\in \R^{m}\times\R^{m\times n}\times \R^{|\B|}$ are dual variables. The KKT conditions for \eqref{prob-alm-prim-v2} are
\begin{subequations}
\label{KKT-sub-RR}
    \begin{align}
        & Ax=b,\ \operatorname{vec}(AX-bx^\top)=0,\ \operatorname{diag}_\B (X)=x_\B ,\ Y=\begin{pmatrix}1& x^\top \\ x & X\end{pmatrix}\succeq 0,\label{KKT-sub-RR-line1}\\
        & C-\mathcal{A}^*(y)-S-W=0,\ \langle{S}, Y\rangle=0,\ S\succeq 0,\label{KKT-sub-RR-line2}
    \end{align}
\end{subequations}
where $y := (\blambda_1;\operatorname{vec}(\blambda_2);\bmu;\alpha)\in\R^{m+mn+|\B|+1}$ and $S$ are the dual variables. Now we show how to recover the dual variables of \eqref{prob-sub-v1}, \eqref{prob-alm-prim-v2}, and \eqref{prob-P-dnn-YZ} from the KKT solutions of \eqref{prob-Rie-sub} in the following theorem.
\begin{theorem}\label{thm-recover}
    Assume that $(R; \blambda, \bmu)$ satisfies the KKT conditions \eqref{KKT-SDPR} of \eqref{prob-Rie-sub}, define 
    \begin{align*}
    Y& := \RR,\quad
    Z := \Pi_\mathcal{P}(Y-\sigma^{-1}\widetilde W),\quad
    W := \sigma\Pi_{\mathcal{P}^*}(\sigma^{-1}\widetilde W- Y),\\
    \blambda_1& := (A^\dagger)^\top(\p +{L A^\dagger b}),\quad
    \blambda_2 :=  (A^\dagger)^\top L(2I-A^\dagger A),
    \\
    \alpha &:= -(J_ARe_1)^\top L(J_ARe_1)-W_{11},\quad
    y := (\blambda_1;\operatorname{vec}(\blambda_2);\bmu;\alpha),\quad
    S := C-\mathcal{A}^*(y)-W,
    \end{align*}
    where $J_A := I-A^\dagger A$ and $A^\dagger=A^\top (AA^\top)^{-1}$. Then the dual variable $S$ can be written as
    \begin{equation}\label{dual-var-S-short}
        S=\begin{pmatrix}
            -(Re_1)^\top\\I
        \end{pmatrix}
        J_ALJ_A\begin{pmatrix}
            -Re_1&I
        \end{pmatrix}.
    \end{equation}
    Moreover, the following statements hold:
    \begin{enumerate}
        \item $(R; \blambda_1,\blambda_2, \bmu)$ satisfies the KKT conditions \eqref{KKT-sub-R} of problem \eqref{prob-sub-v1};
        \item $(Y;y,S)$ satisfies the KKT conditions \eqref{KKT-sub-RR} of \eqref{prob-alm-prim-v2} except $S\succeq 0$;
        \item $(Y,Z;y,S,W)$ satisfies the KKT conditions \eqref{KKT-DNN} of \eqref{prob-P-dnn-YZ} except $S\succeq 0$ and $Y-W=0$.
    \end{enumerate}
\end{theorem}
\proof{}
    By the definition of $\p$ and $L$, the dual variable $S$ can be equivalently written as
    \begin{equation*}
        S=\begin{pmatrix}
            -\alpha-W_{11}
            &\cfrac{\p^\top-\blambda_1^\top A+b^\top \blambda_2}{2}\\
            \cfrac{\p-A^\top \blambda_1 + \blambda_2^\top b}{2}
            &L-\cfrac{A^\top \blambda_2+\blambda_2^\top A}{2}
        \end{pmatrix}.
    \end{equation*}
    Now we prove that it is equal to \eqref{dual-var-S-short}. The diagonal blocks can be easily verified. We only need to prove the left bottom block. By the properties of pseudoinverse, we have $A^\dagger A=(A^\dagger A)^\top$, $AA^\dagger A=A$, $J_A=J_A^\top$ and $J_AA^\top=0$. By multiplying $e_1$ throughout  \eqref{KKT-SDPR-line1}, we get $\p=A^\top{\blambda}e_1-2LRe_1$.
    Now we have
    \begin{align*}
        &\p-A^\top\blambda_1+\blambda_2^\top b
        =\p-A^\top(A^\dagger)^\top(\p +{L A^\dagger b})+(2I-A^\dagger A)LA^\dagger b\\
        &=J_A(\p +2LA^\dagger b)
        =J_A(A^\top{\blambda}e_1-2LRe_1+2LA^\dagger b)
        =2J_AL(A^\dagger b-Re_1)\\
        &=2J_AL(A^\dagger be_1^\top e_1-Re_1)
        =2J_AL(A^\dagger AR e_1-Re_1)
        =-2J_ALJ_ARe_1,
    \end{align*}
    where we have used the fact that $AR=be_1^\top$. Then \eqref{dual-var-S-short} is proven. Next, we prove the rest of the KKT results in the following three parts.
    
    (\textit{1}) For the KKT conditions \eqref{KKT-sub-R} of problem \eqref{prob-sub-v1}, by comparing the KKT conditions \eqref{KKT-sub-R} and \eqref{KKT-DNN}, and noting that $AR=be_1^\top$, it is sufficient to prove that
    $A^\top\blambda=A^{\top} {\blambda_1}e_1^\top+A^\top\blambda_2 R.$
    Indeed
    \begin{align*}
        &A^{\top} {\blambda_1}e_1^\top+A^\top\blambda_2 R 
        = A^{\top} {(A^\dagger)^\top(\p +{L A^\dagger b)}}e_1^\top+A^\top(A^\dagger)^\top L(2I-A^\dagger A) R\\
        &=(A^\dagger A)^\top (\p e_1^\top+2LR)+(A^\dagger A)^\top LA^\dagger ( be_1^\top-AR)
        =(A^\dagger A)^\top A^\top \blambda=A^\top  \blambda,
    \end{align*}
    where the third equality follows from \eqref{KKT-SDPR-line1}.

    (\textit{2}) For the KKT conditions \eqref{KKT-sub-RR} of \eqref{prob-alm-prim-v2},
    the primal feasibility conditions \eqref{KKT-sub-RR-line1} directly follow from \eqref{KKT-SDPR}. We only need to prove that $\langle S, Y\rangle=0$. Since 
    \begin{align*}
    \langle S, Y\rangle
        &=\left\langle
        \begin{pmatrix}
            -Re_1&I
        \end{pmatrix}^\top J_ALJ_A\begin{pmatrix}
            -Re_1&I
        \end{pmatrix},\RR\right\rangle\\
    &=\left\langle
        J_ALJ_A\begin{pmatrix}
            -Re_1&I
        \end{pmatrix}\widehat R,\begin{pmatrix}
            -Re_1&I
        \end{pmatrix}\widehat R
    \right\rangle
    =\left\langle
        J_ALJ_AR(I-e_1e_1^\top),R(I-e_1e_1^\top)
    \right\rangle,
    \end{align*}
    it is sufficient to prove that $J_ALJ_AR(I-e_1e_1^\top)=0$, which holds because
    \begin{align*}
        &2J_ALJ_AR(I-e_1e_1^\top)
        =2(I-A^\dagger A)L(I-A^\dagger A)R(I-e_1e_1^\top)
        =2(I-A^\dagger A)LR(I-e_1e_1^\top)\\
        &=(I-A^\dagger A)(A^\top\blambda-\p e_1^\top)(I-e_1e_1^\top)
        =(A^\top-A^\dagger AA^\top)\blambda(I-e_1e_1^\top)-(I-A^\dagger A)\p(e_1^\top-e_1^\top e_1e_1^\top)=0,
    \end{align*}
    where the second equality come from $AR(I-e_1e_1^\top)=be_1^\top(I-e_1e_1^\top)=0$, and the third equality comes from \eqref{KKT-SDPR-line1}.

    (\textit{3}) Since the KKT conditions \eqref{KKT-sub-RR} are satisfied, we only need to prove $\langle Z,W\rangle=0$, $Z\in\mathcal{P}$ and $W\in\mathcal{P}^*$, which hold because $Z=\Pi_\mathcal{P}(Y-\sigma^{-1}\widetilde W)$ and $W=\sigma\Pi_{\mathcal{P}^*}(\sigma^{-1}\widetilde W- Y)$.
    
\endproof

\begin{remark}
    It is impossible to prove that all KKT conditions of \eqref{prob-P-dnn-YZ} and \eqref{prob-alm-prim-v2} hold because we only update one iteration of the ALM and use the first order KKT condition of \eqref{prob-Rie-sub}, which cannot guarantee global optimality. 
    To prove ${S} \succeq 0$, we need to combine it with a saddle-point escaping strategy. To prove $Y-Z=0$, we need the convergence of ALM.
\end{remark}

\begin{remark}
    The choice of recovered dual variables $y$ and $S$ is not unique. For example, define $\Delta y=(-b;\operatorname{vec}( A);0_{|\B|};\|b\|^2)$, and assume that $y$ and $S$ are the dual variables recovered by Theorem \ref{thm-recover}, then for any $t\geq 0$, the dual variables $y'=y+t\Delta y $ and $S'=S-\mathcal{A}^*(t\Delta y)$ still satisfy all the KKT conditions mentioned in Theorem \ref{thm-recover}. 
\end{remark}

\subsection{Escaping from saddle points}\label{subsec-escape}
In this subsection, we show how to escape from a saddle point of \eqref{prob-Rie-sub} by increasing the rank. It should be mentioned that we cannot apply many existing results like those in \cite{wang2023decomposition,wang2023solving,tang2024feasible} because they need to assume that either the objective function $f_r(\cdot)$ is twice differentiable or that the feasible set $\N_r$ after direct BM factorization is smooth, which does not hold in our case. 
Recall that the function $f_r:\R^{n\times {r}}\rightarrow\R$ satisfies that for any $R\in\R^{n\times {r}}$,
\begin{equation*}
    f_r(R) := \phi(\widehat R\widehat R^\top)=\langle C,\widehat R\widehat R^\top\rangle +\frac{\sigma}{2}\| \Pi_{\mathcal{P}^*}(\sigma^{-1} \wt W-\widehat R\widehat R^\top) \|^2.
\end{equation*}
The Euclidean gradient of $\phi(Y)$ at the point $Y=\widehat R\widehat R^\top$ is
\begin{equation*}
    \nabla \phi (Y)=C-\sigma \Pi_{\mathcal{P}^*}(\sigma^{-1}  \wt W- Y)=C-W,
\end{equation*}
where $W$ coincides with the definition in the last subsection. The following lemma presents an approximation to the objective function.

\begin{lemma}\label{lemma-taylor}
    Suppose that Assumption \ref{assumption-LICQ} holds and $r\geq\bar r$. For any $R\in\Mab$, $\tau\in\mathbb{N}^+$ and $H\in \R^{n\times \tau}$ such that $AH=0$, define $P := [R,0_{n\times \tau}],\ U := [0_{n\times r},H]$, then the following approximation holds:
    \begin{align*}
        f_{r+\tau}\left(\Re_{P}(tU)\right)=f_r\left(R\right)+\< L, HH^\top \> \ t^2 + o(t^2),
    \end{align*}
    where $L$ is defined in \eqref{def-L}, wherein $\bmu$ is the unique solution of 
    $h_R h_R^* (\blambda,\bmu) = h_R (\nabla f_r(R))$. 
\end{lemma}

\proof{}
    Under Assumption \ref{assumption-LICQ}, $\M_{r+\tau}$ is a Riemannian submanifold of $\R^{n\times {(r+\tau)}}$, where $\M_{r+\tau}=\{\bar R\in\R^{n\times (r+\tau)}: A\bar R=b\bar e_1^\top,\diag_\B (\bar R\bar R^\top)=\bar R_\B  \bar e_1\}$ with $\bar{e}_1$ being the first standard unit vector in $\mathbb{R}^{r+\tau}$. We have that $P\in \M_{r+\tau},\ U\in \operatorname{T}_P\M_{r+\tau}$ and the projection retraction $\operatorname{Rtr}_{P}(\cdot)$ is a second order retraction such that for $t\in\R$,
\begin{equation}\label{eq-retrac-taylor}
    \operatorname{Rtr}_{P}(tU)= P+tU+\cfrac{t^2}{2}V+o(t^2)
\end{equation}
for some $V\in(\operatorname{T}_P\M_{r+\tau})^\perp$. Since $\operatorname{T}_P\M_{r+\tau}$ is the null space of the linear map $h_P$, there exists $(\hblambda,\hmu) \in \R^{m\times (r+\tau)}\times \R^{|\B|}$ such that
\begin{equation}\label{eq-retrac-second-order}
    V=h_P^*(\hblambda,\hmu)=A^\top\hblambda+\operatorname{diag}_\B^*(\hmu) (2P - ee_1^\top) .
\end{equation}
Since $\operatorname{Rtr}_{P}(tU)\in\M_{r+\tau}$, substituing \eqref{eq-retrac-taylor} and \eqref{eq-retrac-second-order} into the constraint equations in $\M_{r+\tau}$ and using the fact that the coefficient of $t^2$ is zero, we obtain that $(\hblambda,\hmu)$ is the unique solution of the following linear system:
\begin{equation}\label{eq-retrac-mu-lambda}
     h_{P}(h_{P}^*(\hblambda,\hmu))=(0_{m\times(r+\tau)},-2\diag_\B (HH^\top)).
\end{equation}
From \eqref{eq-retrac-second-order} and \eqref{eq-retrac-mu-lambda}, we know that $AV=0_{m\times (r+\tau)}$. Since the last $\tau$ columns of $AV$ and $2P-ee_1^\top$ are zeros, the last $\tau$ columns of $AA^\top\hblambda$ are also zeros, which further implies that the last $\tau$ columns of $\hblambda$ are zeros since $AA^\top$ is nonsingular. Thus we can assume that $V=[V_1,0_{n\times\tau}]$ for some $V_1\in \R^{n\times r}$. Because $\phi$ is continuously diffierentiable, $f_{r+\tau}$ is also continuously diffierentiable. From \eqref{eq-retrac-taylor}, we have the following result:
    \begin{align}
        f_{r+\tau}\left(\operatorname{Rtr}_{P}(tU)\right)&=f_{r+\tau}\left(P+tU+\frac{t^2}{2}V\right)+o(t^2)\notag\\
        &=f_{r+\tau}\left(\left(R+\frac{t^2}{2}V_1,tH\right)\right)+o(t^2)\notag\\
        &=\phi\left(\RR+\cfrac{t^2}{2}\L_R(V_1)+t^2\begin{pmatrix}
        0&0\\
        0&HH^\top 
    \end{pmatrix}\right)+o(t^2)\notag\\
    &=f_r(R)+t^2\left\langle \nabla \phi(\RR),\cfrac{1}{2}\L_R(V_1)+ \begin{pmatrix}
        0&0\\
        0&HH^\top 
    \end{pmatrix}\right\rangle +o(t^2),\label{eq-retrac-f-res}
    \end{align}
where $\L_R$ is defined in \eqref{def-Lmap}. By the definition of $\L^*_R$ in \eqref{def-Lmap*}, we have that
\begin{align}
    \< \nabla \phi(\RR),\L_R(V_1)\>&=\< \L_R^*(\nabla \phi(\RR)),V_1\>=\< 2\begin{pmatrix}
        0_n& I_n
    \end{pmatrix}\nabla \phi(\RR)\widehat R,V_1\>
    =\< \nabla f_r(R),V_1\>\notag\\
    &=\< \operatorname{grad} f_r(R)+h_R^*(\blambda,\bmu),V_1\>
    =\< \operatorname{grad} f_r(R),V_1\>+\< (\blambda,\bmu),h_R(V_1)\>,\label{eq-grad-L}
\end{align}
where $(\blambda,\bmu)$ is the unique solution to 
$h_R h_R^* (\blambda,\bmu) = h_R (\nabla f_r(R))$.
By \eqref{eq-retrac-second-order}, we know that $V_1=h_R^*(\hblambda,\hmu)\in (\T_R\M_r)^\perp$. Since $\operatorname{grad} f_r(R)\in \T_R\M_r$, we have 
\begin{equation}\label{eq-grad-V1-0}
    \< \operatorname{grad} f_r(R),V_1\>=0.
\end{equation}
Also, from \eqref{eq-retrac-second-order} and \eqref{eq-retrac-mu-lambda}, we have that
\begin{equation}\label{eq-hRV1}
h_R(V_1)=h_Rh^*_R(\hblambda,\hmu)=(0_{m\times r},-2\diag_\B (HH^\top)).
\end{equation}
Substituting \eqref{eq-grad-V1-0} and \eqref{eq-hRV1} into \eqref{eq-grad-L}, we get
\begin{equation}\label{eq-grad-L-final}
    \< \nabla \phi(\RR),\L_R(V_1)\>=-2\< \bmu,\diag_\B (HH^\top)\>.
\end{equation}
Substituting \eqref{eq-grad-L-final} into \eqref{eq-retrac-f-res}, we have
\begin{equation*}
    f_{r+\tau}\left(\operatorname{Rtr}_{P}(tU)\right)=f_r(R)+t^2 \< C_{22}-W_{22}- \diag_\B^*(\bmu),HH^\top\>+o(t^2).
\end{equation*}
By the definition of $L$ in \eqref{def-L}, we get the desired result. 
\endproof
With Lemma \ref{lemma-taylor}, we can escape from a saddle point by increasing the rank.
\begin{theorem}\label{thm-second-order-descent}
    Suppose that Assumption \ref{assumption-LICQ} holds and $r\geq\bar r$. Let $\tau\in\mathbb{N}^+$ be a positive number and $V\in\R^{n\times \tau}$ be a matrix whose columns consist of eigenvectors corresponding to negative eigenvalues of $S$ recovered from Theorem \ref{thm-recover}. Then $U := [0_{n\times r},J_A\begin{pmatrix}
        -Re_1&I
    \end{pmatrix}V]$ is a descent direction of \eqref{prob-Rie-sub} at the point $P := [R,0_{n\times \tau}]$, namely, for some $\beta<0$,
    \begin{equation}\label{eq-taylor-second-order}
        f_{r+\tau} (\operatorname{Rtr}_P(tU))=f_r(R)+\beta t^2+o(t^2).
    \end{equation}
\end{theorem}
\proof{}
    The Taylor expansion \eqref{eq-taylor-second-order} follows from Lemma \eqref{lemma-taylor}, where $\beta$ is computed by
    \begin{align*}
    \beta&=\< L, (J_A\begin{pmatrix}
        -Re_1&I
        \end{pmatrix}V)(J_A\begin{pmatrix}
        -Re_1&I
        \end{pmatrix}V)^\top \>\\
    &=\< (J_A\begin{pmatrix}
        -Re_1&I
        \end{pmatrix})^\top L(J_A\begin{pmatrix}
        -Re_1&I
        \end{pmatrix}) , VV^\top \>
    =\<S,VV^
    \top\><0.
    \end{align*}
    The third equality follows from \eqref{dual-var-S-short} in Theorem \ref{thm-recover}. 
\endproof

\begin{remark}
    By Theorem \eqref{thm-second-order-descent}, we can always find a descent direction if $S\not\succeq 0$. When $S\succeq 0$, by Theorem \eqref{thm-recover}, all KKT conditions of the subproblem \eqref{prob-alm-prim-v1} hold, and $Y=\RR$ is a global minimizer of \eqref{prob-alm-prim-v1}, hence $R$ is also a global minimizer of \eqref{prob-Rie-sub}. 
\end{remark}

\subsection{Convergence analysis of ALM} 
In this subsection, we establish the global convergence of the ALM outlined in 
Algorithm \ref{alg1} for solving \eqref{prob-P-dnn-YZ}. Denote the indicator function $p(Y,Z) := \delta_{\M}(Y)+\delta_{\mathcal{P}}(Z)$ and its conjugate function as $p^*$, which is the support function of $\M\times \mathcal{P}$. Consider the following equivalent form of \eqref{prob-P-dnn-YZ}:
\begin{equation}\label{P1}
    \min_{Y,Z}\left\{\<C,Y\>+p(Y,Z)  :\  Y-Z=0 \right\}. 
\end{equation}
It should be mentioned that the splitting technique in \eqref{P1} is a standard natural approach for applying a Riemannian-based ALM to solve an SDP problem. While such a technique has been used in \cite{wang2023decomposition,wang2023solving}, the key difference between our work and theirs is that they only keep simple linear constraints in the ALM subproblem for which the resulting feasible set after the BM factorization is a simple smooth manifold such as the oblique manifold or Stiefel manifold. In our case, we keep all the linear constraints in the ALM subproblem and the resulting feasible set $\N_r / \M_r$ after BM factorization is no longer a simple manifold. Thus, we cannot directly apply the existing convergence results.
Let $L : \S^{n+1} \times \S^{n+1}\times \S^{n+1} \rightarrow (-\infty,+\infty]$ be the Lagrangian function of \eqref{P1} in the extended form:
\begin{equation*}
L(Y, Z;W) :=  \begin{cases}\<C,Y\>-\langle W, Y-Z\rangle & (Y,Z)\in\operatorname{dom} p, \\ 
+\infty & \text{otherwise} .\end{cases}
\end{equation*}
The Lagrangian dual of \eqref{P1} takes the form of
\begin{equation}\label{D1}
    \max_{W}\left\{\inf_{Y,Z} \,\, L(Y,Z,W) \right\},
\end{equation}
which is equivalent to
\begin{equation}\label{D2}
    \max_{W,S_1,S_2}\left\{-p^*(S_1,S_2)  :\  C+S_1-W=0,\ S_2+W=0 \right\}. 
\end{equation}
We use $\Omega_{D}$ to denote the set of all optimal solutions to the dual problem \eqref{D2}.
For any given $k\geq 0$ and $W^k$, let 
\begin{equation}
\left\{\begin{aligned}
\wt{W}^k(Y,Z) :&=W^k-\sigma_k(Y-Z),  \\
\wt{S}^k(Y,Z) :&=\operatorname{Prox}_{p^*}(Y+\wt{W}^k-C,Z-\wt{W}^k), \\
E^k(Y,Z) :&=(Y,Z)-\operatorname{Prox}_{p}(Y+\wt{W}^k-C,Z-\wt{W}^k) \\
&=(C-\wt{W}^k,\wt{W}^k)+\wt{S}^k(Y,Z),
\end{aligned} \quad\quad (Y,Z) \in \operatorname{dom} p. \right.
\end{equation}
We use the following stopping criteria for the subproblem in Step 4 of Algorithm \ref{alg1}:
\begin{align}
    &f_{r_k}(R^{k+1})-\inf_{R\in \mathcal{M}_{r_k}}f_{r_k}(R)\leq \epsilon_k^2/2\sigma_k,
    \quad 
    &\epsilon_k\geq 0,\quad \sum _{k =0}^{\infty} \epsilon_k<\infty,
    \label{stop-cre-1}\tag{C1}\\
    &
    \begin{aligned}
        \|E^{k+1}\|\leq\ &
    \cfrac{\hat{\epsilon}_k^2 / \sigma_k}{1+\|(Y^{k+1},Z^{k+1})\|+\|(W^{k+1},S^{k+1})\|}\\
    &\min\left\{\cfrac{1}{\|W^{k+1}-W^{k} \|/\sigma_k + 1/\sigma_k}\ ,\ 1\right\},
    \end{aligned}\quad 
    &\hat\epsilon_k\geq 0,\quad \sum _{k =0}^{\infty} \hat\epsilon_k<\infty,
    \label{stop-cre-2}\tag{C2}
\end{align}
where 
\begin{align*}
    W^{k+1} := \wt W^k(Y^{k+1},Z^{k+1}),\ S^{k+1} := \wt S^k(Y^{k+1},Z^{k+1}),\ E^{k+1} := \wt E^k(Y^{k+1},Z^{k+1}),
\end{align*}
and $f_{r_k}:\R^{n\times r_k}\rightarrow\R$ denotes the objective function of the subproblem \eqref{prob-sub-v1} at the $k$-th iteration. The two criteria can be reached by applying the Riemannian gradient descent method together with the saddle-point escaping strategy described in Subsection \ref{subsec-escape}. Now the global convergence of ALM for \eqref{P1} with criteria \eqref{stop-cre-1} follows from \cite[Theorem 4]{rockafellar1976augmented} and \cite[Proposition 2]{cui2019r}.
\begin{theorem}
Suppose that Assumption \ref{assumption-mild} holds and $\Omega_{D}$ is nonempty. Let $\{(Y^{k},Z^k,W^k)\}$ be an infinite sequence
generated by the ALM for \eqref{P1} under criterion \eqref{stop-cre-1}. Then, the dual sequence $\{W^k\}$
converges to some optimal solution $W^*$, and the primal sequence $\{(Y^k,Z^k)\}$ satisfies that for all $k \geq 0$,
\begin{align*}
    & \|Y^k-Z^k\|=(1/\sigma_k)\|W^{k+1}-W^k\|\rightarrow 0,\\
    & \langle C,Y^{k+1}\rangle-\inf \eqref{P1} \leq f_{r_k}(R^{k+1})-\inf_{R\in \Mab}f_{r_k}(R) + (1/2\sigma_k)(\|W^k\|^2-\|W^{k+1}\|^2).
\end{align*}
Moreover, if \eqref{P1} admits a nonempty and bounded solution set, then the sequence 
$\{(Y^k, Z^k)\}$ is also bounded, and all of its accumulation points are optimal solutions to \eqref{P1}.
\end{theorem}

The stopping criterion \eqref{stop-cre-1} is usually difficult to verify. Next we establish the global convergence under a more practical criterion \eqref{stop-cre-2}. We first make some assumptions.
\begin{assumption}\label{assum-lip}
The function $p^*$ is globally Lipschitz continuous on $\operatorname{dom} p^*$.
\end{assumption}

\begin{assumption}\label{assum-RCQ}
    The set $\Omega_{D}$ is nonempty and the Robinson constraint qualification (RCQ) of \eqref{D2} holds at some optimal point $(\bar W,\bar S_1,\bar {S_2})$ (c.f. \cite[Section 3.4.1]{bonnans2013perturbation}).
\end{assumption}

\begin{assumption}\label{assum-dual-distance}
There exists a constant $\gamma$
such that for any $W\in \S^{n+1}$ and $(S1,S2)\in \operatorname{dom} p^*$,
\begin{equation*}
\operatorname{dist}\big(( W,S_1,S_2), \Omega_{D}\big) \leqslant \gamma(1+\|( W,S_1,S_2)\|)\left\|\begin{pmatrix}
    C+S_1-W\\
    S_2+W
\end{pmatrix}\right\|.
\end{equation*}
\end{assumption}
Assumption \ref{assum-lip} holds when $\M$ is bounded. Assumptions \ref{assum-lip} and \ref{assum-RCQ} are the same as the assumptions in \cite{cui2019r}. We make one additional Assumption \ref{assum-dual-distance}, which corresponds to Lemma 5 in \cite{cui2019r} due to the potential lack of interior point of $\M$. Now we state the global convergence of the ALM under criterion \eqref{stop-cre-2} in the following theorem, which is essentially adopted from \cite[Theorem 2]{cui2019r}. 
\begin{theorem}
Suppose that Assumptions \ref{assumption-mild} and \ref{assum-lip}-\ref{assum-dual-distance} hold. Let $\{(Y^{k},Z^k,W^k)\}$ be an infinite sequence generated by the ALM for \eqref{P1} under criterion \eqref{stop-cre-2}. Then, the dual sequence $\{W^k\}$
converges to some optimal solution $W^*$, the primal sequence $\{(Y^k,Z^k)\}$ is bounded, and all of its accumulation points are optimal solutions to \eqref{prob-P-dnn-YZ}.
\end{theorem}

\section{Geometric properties of $\M_r$}\label{sec-alg-var}

Our approach is different from existing works on using the BM factorization approach to solve the ALM subproblem \eqref{ALM-sub-LR}. Existing works only focus on feasible sets that are obvious Riemannian manifolds such as the oblique manifold or the Stiefel manifold, and do not consider more sophisticated feasible sets, which are no longer manifolds as in the case of $\N_r$ in our problem, under the direct BM factorization. Moreover, existing works either do not handle additional nonnegative cone constraints or have not demonstrated the effectiveness of their approaches in handling DNN problems. Here we conduct a more refined analysis of the smoothness property of the algebraic variety and propose an equivalent algebraic reformulation of the feasible set $\N_r$ to obtain better smoothness properties of the equivalently reformulated set $\M_r$.

In this section, we first present the smoothness analysis of the algebraic variety $\M_r$, and then propose an approach to avoid non-smoothness. Next, we show how to compute the projection onto the tangent space $\operatorname{T}_R{\M_r}$ efficiently. After that, we demonstrate how to compute the retraction onto $\M_r$ by transforming a non-convex projection problem onto $\M_r$ into a convex generalized geometric median problem. Finally, we design a generalized Weiszfeld algorithm with convergence guarantees to solve the latter convex problem. All the analysis in this section can be directly extended to a more general algebraic variety $\M_r^g$ defined in \eqref{manifold-general}.

\subsection{Smoothness analysis of $\M_r$}
Recall that the algebraic variety $\mathcal{M}_r$ in \eqref{prob-Rie-sub} is defined as
\begin{align}
\mathcal{M}_r := &\left\{R \in \mathbb{R}^{n \times r} :\  AR=be_1^\top,\ \operatorname{diag}_\B (RR^{\top})-R_\B  e_1=0 \right\}\label{mani-Mr}\\
=&\left\{R \in \mathbb{R}^{n \times r} :\  AR=be_1^\top,\ 
\operatorname{diag}_\B \big((2R-ee_1^\top)(2R-ee_1^\top)^{\top}\big)=e \right\}\label{mani-spherical}
\end{align}
with $\operatorname{dim}(\mathcal{M}_r)=nr-mr-|\B|$. The tangent space at $R \in \mathcal{M}_r$ is
\begin{equation*}
\mathrm{T}_R \mathcal{M}_r := \left\{H \in \mathbb{R}^{n\times r} :\  AH=0,\ 2\operatorname{diag}_\B (H R^{\top})-H_\B  e_1=0 \right\} .
\end{equation*}
Recall the linear operator 
$h_R:\mathbb{R}^{n \times r}\rightarrow  \mathbb{R}^{m\times r}\times\mathbb{R}^{|\B|} $
defined in \eqref{eq-hR} with 
\[
h_R(H) :=  (AH ; 2\operatorname{diag}_\B (H R^{\top})-H_\B  e_1), \quad H \in \mathbb{R}^{n \times r}.
\]
The adjoint mapping $h_R^*: \mathbb{R}^{m\times r}\times\mathbb{R}^{|\B|} \rightarrow \mathbb{R}^{n \times r}$ satisfies that for any $(\blambda,\bmu) \in \mathbb{R}^{m\times r}\times \mathbb{R}^{|\B|} $,
\begin{equation}\label{eq-lambda-mu-grad}
    h_R^*(\blambda, \bmu) := A^\top\blambda+\operatorname{diag}_\B^*(\bmu) (2R - ee_1^\top) .
\end{equation}
Then we have $H\in \mathrm{T}_R \mathcal{M}_r$ if and only if $h_R(H)=0$. The following proposition characterizes the smoothness of $\M_r$.
\begin{proposition}[Smoothness of $\M_r$]\label{prop-regular}
    For any $R \in \mathcal{M}_r$, define the linear operator $P:\R^{m\times r}\rightarrow\R^{m\times r}$ such that for any $\blambda\in \R^{m\times r}$,
    \begin{equation*}
        P(\blambda) := {A} {A}^\top\blambda-{A}\operatorname{diag}_\B^*\big(\operatorname{diag}_\B ({A}^\top\blambda (2R-ee_1^\top)^{\top})\big) (2R-ee_1^\top),
    \end{equation*}
    then the following statements are equivalent:
    \begin{enumerate}
        \item The point $R$ is regular;
        \item The operator $h_R^*(\blambda,\bmu)$ is injective; 
        \item The operator $P$ is positive definite, i.e., $\langle \blambda,P(\blambda)\rangle> 0$ for all non-zero $\blambda\in \R^{m\times r}$.
    \end{enumerate}
\end{proposition}
\proof{}
    For the case where $\B=\emptyset$, the proposition holds trivially since $AA^\top$ is positive definite. Thus we only need to consider the case where $B\not=\emptyset$.  The equivalence between \textit{1} and \textit{2} directly follows from the definition of the LICQ condition. We now prove the equivalence between \textit{1} and \textit{3}. For $\textit{1}\Rightarrow \textit{3}$, 
    suppose $R \in \mathcal{M}_r$ is regular, then the spherical constraints in \eqref{mani-spherical} indicates that $\|(2R-ee_1^\top)_i\|=1$ for any $i \in B$. Thus for any $\blambda \in \mathbb{R}^{m\times r}$, we have
    \begin{equation}\label{eq-thm1-1}
    \begin{aligned}
     \langle P(\blambda), \blambda\rangle&=\|A^\top \blambda\|^2-\|\operatorname{diag}_\B (A^\top \blambda (2R-ee_1^\top)^{\top})\|^2 \\
    & \geq\| A^\top\blambda\|^2-\sum_{i\in B}\|( A^\top \blambda)_i\|^2\|(2R-ee_1^\top)_i\|^2 \geq 0,
    \end{aligned}
    \end{equation}
    where the first inequality follows from the Cauchy-Schwartz inequality and the second one follows from $\|(2R-ee_1^\top)_i\|=1$ for any $i \in B$. Assume by contradiction that the operator $P(\cdot)$ is not positive definite, then there exists non-zero $\blambda \in \mathbb{R}^{m\times r}$ such that $ P(\blambda)=0$, which implies that equalities hold throughout \eqref{eq-thm1-1}. Thus, we must have $(A^\top\blambda)_{[n] \backslash B}=0$ and there exists $\bmu \in \mathbb{R}^{|\B|}$ such that ${A}^\top\blambda+\operatorname{diag}_\B^*(\bmu) (2R-ee_1^\top)=0$, which contradicts that $R$ is regular. 
    
    For $\textit{3}\Rightarrow\textit{1}$,
    assume by contradiction that $R$ is not regular, then there must exist non-zero $(\blambda,\bmu) \in$ $\mathbb{R}^{m\times r}\times\mathbb{R}^{|\B|} $ such that ${A}^\top\blambda+\operatorname{diag}_\B^*(\bmu) (2R-ee_1^\top)=0$, which indicates that $(A^\top \blambda)_{[n] \backslash B}=0$. If $\blambda=0$, then we have that $\bmu=0$ due to the fact that every row of $(2R-ee_1^\top)_\B $ has length 1, which contradicts that $(\blambda,\bmu) \neq 0$. Thus, $\blambda \neq 0$. Plug $A^\top \blambda=-\operatorname{diag}_\B^*(\bmu) (2R-ee_1^\top)$ into \eqref{eq-thm1-1}, we have that equalities hold throughout, hence $\langle P(\blambda),\blambda\rangle=0$. This contradicts that $P\succ 0$. 
\endproof
\begin{remark}
    The regularity of $R\in\Mab$ cannot be guaranteed in general. However, we can always avoid the non-regularity by equivalently reformulating the original problem \eqref{prob-MBQP} and the DNN problem \eqref{prob-dnn-new}, which is described in the next subsection.
\end{remark}
\begin{remark}
    All the analysis of $\M_r$ can be extended to a more general algebraic variety:
    \begin{equation}\label{manifold-general}
        \mathcal{M}^g_r := \left\{R \in \mathbb{S}^{n \times r} :\  \mathcal{A}(R)=b,\ \operatorname{diag}_\B (RR^{\top})=p \right\},
    \end{equation}
    where $\mathcal{A}:\R^{n\times r}\rightarrow \R^{m}$ is a surjective linear operator, $\B\subseteq[n]$ is an index set, and $p\in \R^{|\B|}_+$ is a nonnegative vector. The algebraic variety $\M_r$ is a special case of $\M^g_r$ by a linear transformation $R'=2R-ee_1^\top$. Another application of the set $\M_r^g$ can be seen in \cite{tang2024solving}.
\end{remark}

\subsection{Avoiding non-smoothness}\label{subsec-avoid-non-reuglar}
In this part, we provide strategies to avoid non-smoothness in $\M_r$. The algorithmic design and theoretical analysis are based on the smoothness of $\M_r$. Without the smoothness property, the singularity of the linear map $h_Rh_R^*$ will not only slow down the computation of the projection and retraction but also affect the feasibility of the dual variables. For some special structured $\M_r$, one can easily verify the smoothness according to Proposition \ref{prop-regular}. 
\begin{example}
    {$\M_r$ is smooth for any $r\geq 1$ if either $\B=\emptyset$ or $m=0$.}
\end{example}
However, $\M_r$ can be non-smooth in some cases, as shown in the following example.
\begin{example}\label{example-rank-1-nonsmooth}
{$\mathcal{M}_r$ is non-smooth if $r=1$, $\B=[n]$ and $m\geq 1$.}
\end{example}
\proof{}
    When $r=1$, $\B=[n]$ and $m\geq 1$, the set $\M_r$ can be simplified as
    \begin{equation}\label{set-M1}
        \M_1:=\{R\in\R^n :\ AR=b,\ (2R-e)^2=e \},
    \end{equation}
    where $(2R-e)^2=((2R_1-1)^2,\cdots,(2R_n-1)^2)^\top$. For any $R\in\M_1$, let $\blambda=e_1\in\R^{m}$ and $\bmu=-a_1\oslash (2R-e)\in\R^{n}$, where ``$\oslash$" is the element-wise division operator. The variable $\bmu$ is well-defined because \eqref{set-M1} indicates that every entry of $2R-e$ is non-zero. Then we have
    \begin{equation*}
        h_R^*(\blambda,\bmu)=A^\top\blambda+\diag^*(\bmu)(2R-e)=a_1-a_1 \oslash (2R-e) \circ (2R-e)=0.
    \end{equation*}
    Thus $h_R^*(\blambda,\bmu)$ is not injective. By Proposition \ref{prop-regular}, $\M_1$ is non-smooth. 
\endproof

\begin{figure}[h!]
    \centering
    \includegraphics[width=.5\textwidth]{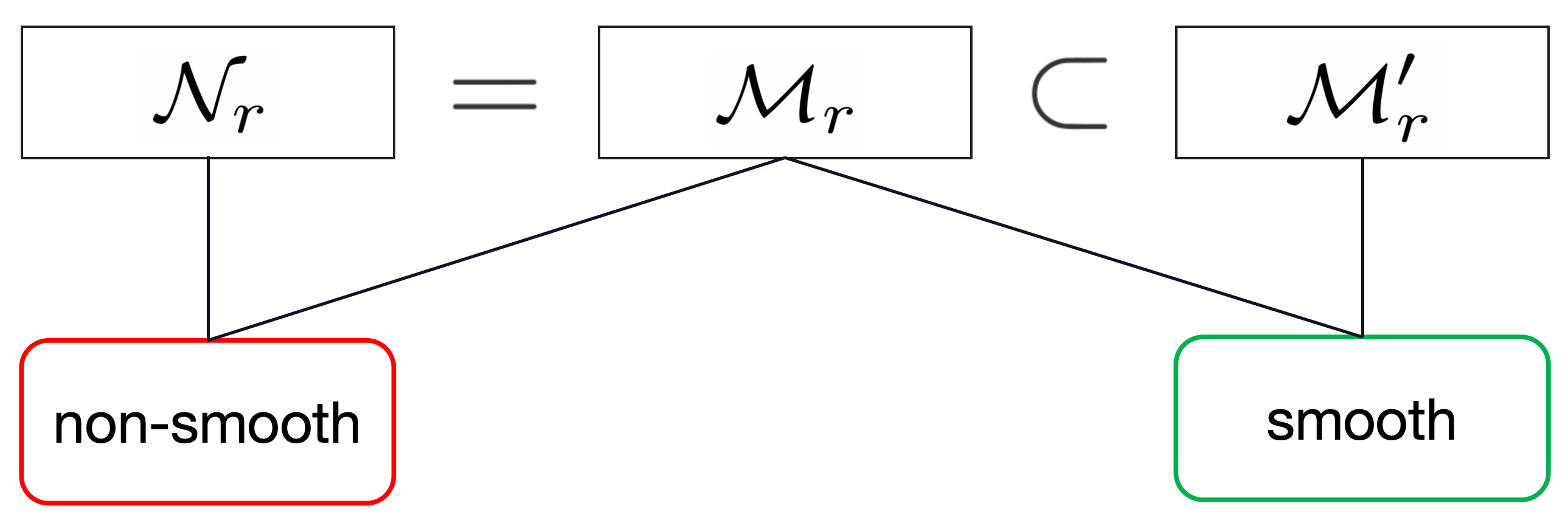}
    \caption{Relationship between different feasible sets and their smoothness properties. For $\N_r$, every point is nonsmooth, whereas for $\M_r'$, every point is smooth. For $\M_r$, it has smooth points but may also have nonsmooth points.}
    \label{fig:relation-manifold}
\end{figure}

Thus, we need to design some strategies to avoid non-smoothness in $\M_r$. One approach is to analyze the local geometric properties at non-regular points and design specialized strategies. For some special $\M_r$ like the one in the quadratic knapsack problem, it is proven in \cite{tang2024feasible} that the non-regular point of the algebraic variety is either the zero point or the integer feasible solutions of the original problem \eqref{prob-MBQP}. However, designing specialized strategies for every type of \eqref{prob-MBQP} problem is difficult and complicated. Here we provide another approach from the modeling perspective. By replacing equality constraints with inequality constraints and adding slack variables, \eqref{prob-MBQP} is equivalent to
%
%
\begin{equation}\label{prob-MBQP-reform}\tag{MBQP$'$}
\min\left\{x'^{\top} Q' x'+2 c'^{\top}x' : 
\begin{array}{l}
     A'x'=b',\ x_i \in\{0,1\},\ \forall i \in B, \\
     x_ix_j=0,\ \forall (i,j)\in E,\ x':=\begin{pmatrix}
    x\\s
\end{pmatrix}\in \R^{n+2m}_+
\end{array}
\right\},
\end{equation}
where $\B\subseteq[n]$, $E\subseteq \{(i,j)\mid 1\leq i<j\leq n \}$ and 
\[
Q':= \begin{pmatrix}
    Q&0_{n\times 2m}\\
    0_{2m\times n}& 0_{2m\times 2m}
\end{pmatrix},\ c':= \begin{pmatrix}
    c\\0_{2m}
\end{pmatrix},\ A' := \begin{pmatrix}
    A&I_{m}&0_{m\times m}\\
    A&0_{m\times m}&-I_{m}
\end{pmatrix},\ b' := \begin{pmatrix}
    b\\b
\end{pmatrix}.
\]
The DNN relaxation of the new equivalent problem \eqref{prob-MBQP-reform} is given by
\begin{equation}\label{prob-dnn-reform}
    \min\left\{
\<C',Y'\>  :\  
Y'\in \F'\cap\Z'\cap \mathbb{S}_{+}^{n+2m+1} \cap\mathbb{N}^{n+2m+1}
\right\},
\end{equation}
where $C':= [1,(c')^\top; c',Q']$ and 
\begin{align*}
    \Z'&:=\left\{\begin{pmatrix}z'& x'^\top \\ x' & X'\end{pmatrix}\in  \mathbb{S}^{n+2m+1} :\ X'_{ij}=0, \ \forall (i,j)\in E \right\},\\
    \F' &:= \left\{\begin{pmatrix}1& x'^\top \\ x' & X'\end{pmatrix}\in  \mathbb{S}^{n+2m+1} :\  A'x'=b',\ A'X'=b'(x')^\top,\ x'_i =X'_{ii}, \  \forall i \in B\right\}.
\end{align*}
The new algebraic variety for the subproblem is:
\[
\mathcal{M}_r' := \left\{R \in \mathbb{R}^{(n+2m) \times r} :\  A'R=b'e_1^\top,\ \operatorname{diag}_\B (RR^{\top})=R_\B  e_1 \right\}.
\]
We first show that the DNN relaxations of two different reformulations are equivalent.
\begin{lemma}\label{lemma-DNN-equivalent}
The new DNN relaxation \eqref{prob-dnn-reform} is equivalent to the original DNN relaxation \eqref{prob-dnn-new}.
\end{lemma}
\proof{}
We prove by showing that \eqref{prob-dnn-reform} can be equivalently reduced to \eqref{prob-dnn-new}. For any feasible point $Y'$ of \eqref{prob-dnn-reform}, we have $A'x'=b'$ and $x'\geq 0$. By eliminating the variable $x$, we get $(I_m,I_m)s=0$ and $s\geq 0$, hence $s=0$. Similarly, since $A'X'=b'(x')^\top$ and $X'\geq 0$, 
we can show that the last $2m$ rows of $X'$ are equal to $0$. Thus, we have
\begin{equation*}
    x'=\begin{pmatrix}
        x\\0_{2m}
    \end{pmatrix},\ X'=\begin{pmatrix}
        X&0_{n\times 2m}\\
        0_{2m\times n}&0_{2m\times 2m}
    \end{pmatrix}.
\end{equation*}
By plugging the formulas above into \eqref{prob-dnn-reform} and eliminating redundant constraints, we can easily verify that the reduced problem is equivalent to \eqref{prob-dnn-new}. 
\endproof
Next, we prove that the new algebraic variety $\M_r'$ is a smooth manifold.
\begin{lemma} \label{lemma-avoid-nonsmooth}
For any positive {integer} $r$, every point $R\in\mathcal{M}_r'$ satisfies the LICQ condition.
\end{lemma}
\proof{}
    Recall that $R \in \mathcal{M}_{r}'$ is regular if and only if the corresponding linear map $h_R^*$, defined in \eqref{eq-lambda-mu-grad} with $A$ replaced by $A'$, is injective. For any $R\in \mathcal{M}_{r}'$, $h_R^*(\blambda,\bmu)=0$ means
    \begin{equation}\label{eq-mani}
     \begin{pmatrix}
            A^\top&A^\top\\
            I_{m}&0_{m\times m}\\
            0_{m\times m}&-I_{m}
        \end{pmatrix}\blambda
        +
        \begin{pmatrix}
            \operatorname{diag}_\B^*(\bmu)&0_{n\times 2m}\\0_{2m\times n}&0_{2m\times 2m}
        \end{pmatrix}
        (2R - ee_1^\top)
        =0,
    \end{equation}
    where $R\in\mathbb{R}^{(n+2m)\times r}$, $\bmu\in\mathbb{R}^{|\B|}$, and $\blambda\in\mathbb{R}^{2m\times r}$. The last $2m$ rows of equation \eqref{eq-mani} indicate that $\blambda=0_{2m\times r}$, so \eqref{eq-mani} can be simplified as $(\operatorname{diag}_\B^*(\bmu), 0_{n\times 2m}) (2R - ee_1^\top)=0$. Since $R\in \mathcal{M}_{r}'$ indicates that $\operatorname{diag}_\B ((2R - ee_1^\top)(2R - ee_1^\top)^{\top})=e $, the $i$-th row of $2R-ee_1^\top$ cannot be $0_{1\times r}$ for any $i\in [B]$. Thus, we must have $\bmu=0_{|\B|}$. The above deviation shows that $h_R^*(\blambda,\bmu)=0$ if and only if $(\blambda,\bmu)=0$. Thus, $h(\blambda,\bmu)$ is injective, and every $R \in \mathcal{M}_{r}'$ is regular. 
\endproof
\begin{remark}
    Although problem \eqref{prob-MBQP} and \eqref{prob-MBQP-reform} are equivalent and their DNN relaxations \eqref{prob-dnn-new} and \eqref{prob-dnn-reform} are also equivalent, we cannot use the same argument to show that the algebraic varieties $\M_r$ and $\M_r'$ are equivalent because we do not have the condition that $R\geq 0$. The relationship between different feasible sets is shown in Figure \ref{fig:relation-manifold}, where we abuse the notation $\subset$ to indicate that for any $R\in \M_r$, we have $[R;0_{2m\times r}]\in \M_r'$. 
\end{remark}

In practice, we first apply RNNAL to solve the DNN problem \eqref{prob-dnn-new}. If the singularity issues occur during the projection and retraction step, we address them by considering the equivalent DNN problem \eqref{prob-dnn-reform} instead and applying RNNAL again. In Subsection \ref{subsec-qap}, we use the quadratic assignment problem (QAP) as an example to demonstrate the strategy's effectiveness.

\subsection{Projection}\label{subsec-proj}
In this subsection, we show how to compute the projection onto the tangent space $\operatorname{T}_R{\M_r}$ efficiently. For any $V \in \mathbb{R}^{n \times r}$, the projection mapping at the point $R\in \mathcal{M}_r$ is the orthogonal projection operator onto its tangent space, which can be computed by
\begin{equation*}
\operatorname{Proj}_{\T_R \mathcal{M}_r}(V) := V-h_R^*(h_Rh_R^*)^{\dagger}h_R(V),
\end{equation*}
where $(\cdot)^\dagger$ denotes the pseudo inverse. The Riemannian gradient is the projection of the Euclidean gradient onto the tangent space, which is given by
\begin{equation*}
    \operatorname{grad}f_r(R)=\operatorname{Proj}_{\T_R \mathcal{M}_r}(\nabla f_r(R))=\nabla f_r(R)-h_R^* (\blambda,\bmu),
\end{equation*}
where $(\blambda,\bmu)\in \R^{m\times r}\times\R^{|\B|}$ is the solution of the following linear system:
\begin{equation}\label{eq-grad}
    h_R( h_R^*(\blambda,\bmu))= h_R(\nabla f_r(R)).
\end{equation}
When $R$ is regular, $h_R^*$ is injective due to Proposition \ref{prop-regular}, and \eqref{eq-grad} has a unique solution. The most expensive step to compute the projection is to tackle the structured system of linear equations \eqref{eq-grad}. We show that equation \eqref{eq-grad} of size $mr+|\B|$ can be transformed into a smaller normal equation of size $|\B|$. First, denote $R':= 2R-ee_1^\top:= (R'_1,\cdots,R'_r)\in \R^{n\times r}$ and $D_i:= (\Diag(R'_i))_\B \in\R^{|\B|\times n}$ for any $i=1,\cdots,r$, then the matrix representation $H$ of the operator $h_R$ such that $\operatorname{vec}(h_R(X))=H\operatorname{vec}(X)$ for any $X\in\R^{n\times r}$ is given by
\begin{equation*}
H:= \begin{pmatrix}
    A&&\\
    &\ddots&\\
    &&A\\
    D_1 & \cdots &D_r
    \end{pmatrix}\quad \in \R^{(m r+|\B|)\times nr} .
\end{equation*}
Thus, \eqref{eq-grad} can be equivalently expressed in the following form:
\begin{equation}\label{eq-proj-linear-equation}
    HH^\top x=d,
\end{equation}
where $x=(\operatorname{vec}(\blambda);\bmu)\in \R^{m r+|\B|}$, $d=H\operatorname{vec}(\nabla f_r(R))\in\R^{m r+|\B|} $.
We aim to accelerate the computation by exploiting the special structure of the matrix $H$. The system of linear equations \eqref{eq-proj-linear-equation} can be equivalently written as
\begin{equation}\label{eq-proj-HHt-block}
    HH^\top x= \begin{pmatrix}
        E_1&E_2\\
        E_2^\top&E_3
    \end{pmatrix}\begin{pmatrix}
        x_1\\
        x_2
    \end{pmatrix} =\begin{pmatrix}
        d_1\\
        d_2
    \end{pmatrix},
\end{equation}
where $x:= (x_1;x_2)\in \R^{m r+|\B|}$, $d:= (d_1;d_2)\in\R^{m r+|\B|}$ and 
\begin{alignat*}{3}
    E_1:&= \Diag(AA^\top,\cdots,AA^\top) &\quad\in\R^{mr\times mr},\\
    E_2:&= (AD_1^\top;\cdots;AD_r^\top) &\quad\in\R^{mr\times |\B|},\\
    E_3:&= D_1D_1^\top +\cdots + D_r D_r^\top&\quad\in\R^{|\B|\times |\B|}.
\end{alignat*}
To reduce the size of the linear system, we next eliminate $x_1$ and compute $x_2$ by
\begin{equation}\label{eq-x2}
    Mx_2=d_3,
\end{equation}
where $d_3:= d_2-E_2^\top E_1^{-1}d_1\in \R^{|\B|}$ and $M:= E_3-E_2^\top E_1^{-1}E_2\in\R^{|\B|\times |\B|}$ is the Schur complement matrix such that
\begin{equation}\label{formula-M}
    M:=  \sum_{i=1}^{r} D_i(I-A^\dagger A)D^\top_i=(I-A^\dagger A)_{BB}\circ (R'R'^\top)_{BB} \quad\in \R^{|\B|\times |\B|},
\end{equation}
where $A^\dagger = A^\top (AA^\top)^{-1}$ and $(X)_{BB}$ is the submatrix of $X$ obtained by selecting the rows and columns from the index set $B$. When $R$ is regular, the coefficient matrix $M$ is positive definite. 

There are two ways of solving \eqref{eq-x2}. When $mr$ is large, we can compute the matrix $M$ using formula \eqref{formula-M} and solve (\ref{eq-x2}) by Cholesky decomposition. Note that when forming $M$, the matrix $(AA^\top)^{-1}$ and hence $I-A^\dagger A$ 
only need to be computed once before executing Algorithm \ref{alg1}, and the matrix-matrix product $RR^\top$ in $R'(R')^\top=4RR^\top-2Re_1e^\top-2ee_1^\top R+ee^\top$ is already calculated when computing the gradient $\nabla f_r(R)$. Therefore, the cost of computing the coefficient matrix only requires $\O(|\B|^2)$ extra arithmetic operations in every step. The total computational cost of solving the linear system is $\O(|\B|^3).$ When $mr$ is relatively smaller than $|\B|,$ we can use the Sherman–Morrison–Woodbury formula to solve (\ref{eq-x2}) with complexity $\O((mr)^2|\B|+(mr)^3)$ because $M=\DD (R'R'^\top)_{BB} -\sum_{i=1}^{r} D_i A^\dagger A D^\top_i$ is a diagonal matrix minus a rank-$mr$ matrix. Combining these two results, we know that the computational complexity of solving (\ref{eq-x2}) is $\O\big(\min\left\{ |\B|^3, (mr)^2|\B|+(mr)^3\right\}\big).$

Finally, it follows from \eqref{eq-proj-HHt-block} that 
$
    x_1=E_1^{-1}(d_1-E_2 x_2).
$
We can also compute $x_1$ efficiently due to the diagonal block structure of $E_1$. To see this, we denote $x_1:= (x_1^1;\cdots;x_1^r)\in \R^{mr}$ and $d_1:= (d_1^1;\cdots;d_1^r)\in \R^{mr}$, then 
\begin{equation}\label{eq-xt1}
    x_1^t=(AA^\top)^{-1}(d_1^t-AD^\top_t x_2),\quad t=1,\cdots,r.
\end{equation}
Therefore, the computational cost of computing $x_1$ from $x_2$ is $\O(m^2r+|\B|mr)$, given that the Cholesky factorization of $AA^\top$ has already been computed. From the above analysis, the total computational cost of a projection is 
\begin{equation}\label{projcost}
\O \big(\min\left\{ |\B|^3+m^2r+mr|\B|, (mr)^2|\B|+(mr)^3\right\}\big),
\end{equation}
which is much better than directly solving \eqref{eq-proj-HHt-block} at the cost of $\O \big((|\B|+mr)^3\big)$ when either $|\B|$ or $mr$ is small. Note that (\ref{projcost}) is just the worst case complexity bound. In practice, the projection could be further accelerated by making use of the sparsity of $A$ or solving \eqref{eq-proj-linear-equation} iteratively by using the PCG method with the diagonal preconditioner, when it is well-conditioned. Finally, we note that if $|\B|=0$, we can solve \eqref{eq-proj-HHt-block} via \eqref{eq-xt1} without the form involving $x_2$.

\subsection{Retraction}\label{subsec-retrac}
In this subsection, we show how to compute the {metric} projection onto $\M_r$ by transforming the non-convex projection problem into a convex problem. {This step is one of the key factors to make sure RNNAL is practically efficient.} For $\bar R\in\M_r$ and $H\in \T_{\bar R}\M_r$, let $V := \bar R+H$. Define the following metric projection:
\begin{equation}\label{metric-proj}
\operatorname{Rtr}_{\bar R}(H) := \operatorname{Proj}_{\mathcal{M}_r}(V)=\arg \min \left\{\|R-V\|_F^2  :\  R\in \mathcal{M}_r\right\}.
\end{equation}
Note that the definition of the retraction mapping is a set-valued mapping because problem \eqref{metric-proj} may have multiple optimal solutions. We always consider the solution set of a projection operator as a single point when the solution is unique. Problem \eqref{metric-proj} is non-convex because of the spherical constraints in $\M_r$. 
However, we can show in Theorem \ref{thm-retrac-comp} that under certain conditions, \eqref{metric-proj} can be transformed into a convex programming problem. Before that, we first define the following spherical manifold:
\begin{equation}\label{Br}
    \Br := \left\{R \in \mathbb{R}^{n \times r}: \operatorname{diag}_\B (RR^{\top})-R_\B  e_1=0\right\},
\end{equation}
which contains the algebraic variety $\mathcal{M}_r$ as a subset. For any $R$ such that every row of $(2R-ee_1^\top)_\B $ is non-zero, the unique projection onto $\Br $ is given by
\begin{equation*}
    \operatorname{Proj}_{\Br }(R)=\cfrac{1}{2}\left(\diag_\B {(\v)}({2R-ee_1^\top})+{ee_1^\top}\right),
\end{equation*}
where $\v\in\R^{n}$ is a vector such that $\v_i=1/\|(2R-ee_1^\top)_i\|$ for any $i\in [B]$, and $\v_i=1$ otherwise. The following lemma explains the relationship between $\operatorname{Proj}_{\Br }$ and $\operatorname{Proj}_{\mathcal{M}_r}$.
\begin{lemma}\label{lemma-retrac}
     For any $V \in \mathbb{R}^{n \times r}$, if there exists $\rv \in \mathbb{R}^{m\times r}$ such that $A\operatorname{Proj}_{\Br }(V+A^\top\rv)=b e_1^\top$, then $\operatorname{Proj}_{\mathcal{M}_r}(V)=\operatorname{Proj}_{\Br }(V+A^\top\rv)$.
\end{lemma}
\proof{}
    First, for any $\rv \in \mathbb{R}^{m\times r},\eqref{metric-proj}$ is equivalent to the following problem:
$$
\min \left\{\|R-(V+ A ^\top\rv)\|_F^2  :\  R\in \mathcal{M}_r\right\}
$$
because with the affine constraint $AR=be_1^\top,\|R-V\|_F^2 $ and $\|R-(V+ A ^\top\rv)\|_F^2 $ only differ by a constant. 
The equivalence indicates that 
\begin{equation}
\label{lemma1-proj-eq1}
\operatorname{Proj}_{\mathcal{M}_r}(V)=\operatorname{Proj}_{\mathcal{M}_r}(V+A^\top\rv).
\end{equation}
Next, because $A\operatorname{Proj}_{\Br }(V+A^\top\rv)=b e_1^\top$ by assumption, we have that $\operatorname{Proj}_{\Br }(V+A^\top\rv) \subseteq \mathcal{M}_r$.
Together with the fact that 
$\mathcal{M}_r \subseteq \Br $, we get
$
\operatorname{Proj}_{\mathcal{M}_r}(V+A^\top\rv)=\operatorname{Proj}_{\Br }(V+A^\top\rv),
$
which proves the desired result with \eqref{lemma1-proj-eq1}. 
\endproof
With Lemma \eqref{lemma-retrac}, we can now compute the retraction by solving a convex problem.

\begin{theorem}[Retraction Computation]\label{thm-retrac-comp}
\label{thm-retrac}
\begin{itemize}
    \item [(i)] For any $V \in \mathbb{R}^{n \times r}$, define $V'=V-\frac{ee_1^\top}{2}$. If there exists an optimal solution $\rv\in \R^{m\times r}$ of the following problem:
\begin{equation}\label{prob-lam}
    \min_{\rv\in\R^{m\times r}}\,\, 
    G(\rv) := \sum_{i \in B}\|(V'+A^\top \rv)_i\|+\sum_{i \in[n] \backslash B} {\|(V'+A^\top \rv)_i\|^2}+\langle (Ae-2b)e_1^\top,\rv\rangle 
\end{equation}
such that every row of $\left(V'+A^\top \rv\right)_\B $ is non-zero, then $\operatorname{Proj}_{\mathcal{M}_r}(V)=\operatorname{Proj}_{\Br }\left(V+A^\top \rv\right)$.
    \item [(ii)] For any regular point $R \in \mathcal{M}_r$, there exists $\epsilon>0$ such that for any $V \in \mathrm{B}_\epsilon(R)$, $\operatorname{Proj}_{\mathcal{M}_r}(V)=\operatorname{Proj}_{\Br }\left(V+A^\top \rv\right)$ for some optimal solution $\rv$ of \eqref{prob-lam} satisfying that every row of $\left(V'+A^\top \rv\right)_\B $ is non-zero.
\end{itemize}
\end{theorem}
\proof{}
    {(i)} 
    Since every row of $(V'+A^\top \rv)_\B $ is non-zero, $G(\rv)$ is differentiable at the point $\rv$. Since $\rv$ is an optimal solution, we have that $\nabla_{\rv} G(\rv)=0$, which implies  $ A \operatorname{Proj}_{\Br }(V+A^\top \rv)-be_1^\top=0$. Then the result follows from Lemma \ref{lemma-retrac}.
    
    {(ii)} For any $(V',\rv)$ such that every row of $\left(V+A^\top \rv\right)_\B $ is non-zero, define the nonlinear mapping $F(V, \rv) :=  A \operatorname{Proj}_{\Br }\left(V+A^\top \rv\right)-be_1^
    \top$. We have that $F(R, 0)=0$ and $F(R,\rv)$ is well-defined and smooth in a neighbourhood of $(R, 0)$. Also, the Jacobian ${\mathcal{J}_{F, \rv}(R, 0)=P(\cdot)\succ 0}$, where $P(\cdot)$ is defined in Proposition \ref{prop-regular}. Thus, by the implicit function theorem, there exists $\epsilon>0$ and a unique function $\rv_g: \mathrm{B}_\epsilon(R) \rightarrow \mathbb{R}^{m\times r}$ such that $\rv_g(R)=0$, $F\left(V, \rv_g(V)\right)=0$ and every row of $\left(V'+A^\top \rv_g(V)\right)_\B $ is non-zero for $V\in \mathrm{B}_\epsilon(R)$. From Lemma \ref{lemma-retrac} we have
    $
    \operatorname{Proj}_{\mathcal{M}_r}(V)=\operatorname{Proj}_{\Br }(V+ A ^\top\rv_g(V)).
    $
    Also, we have $\nabla G_{\rv}(\rv_g(V))=0$, which implies that $\rv_g(V)$ is an optimal solution of \eqref{prob-lam}. 
\endproof

\subsection{Retraction subproblem}
Next, we show how to solve the unconstrained non-smooth convex problem \eqref{prob-lam} by a generalized Weiszfeld algorithm and prove its convergence. Since $\Mab$ is nonempty, $b\in\operatorname{Range}(A)$, we can assume that there exists $b'\in\R^{m\times r}$ such that $(Ae-2b)e_1^\top=Ab'$. Then problem \eqref{prob-lam} can be equivalently written as 
\begin{equation}\label{prob-lam-new}
    \min_{\rv\in\R^{m\times r}}\,\, G(\rv) := \sum_{i \in B}\|(V'+A^\top \rv)_i\|+\sum_{i \in[n] \backslash B} {\|(V'+A^\top \rv)_i\|^2}+\langle b',A^\top\rv\rangle .\tag{Ret-sub}
\end{equation}
The problem \eqref{prob-lam-new} can be regarded as a generalization of the geometric median problem, which can be solved by the Weiszfeld algorithm \cite{weiszfeld1937point,weiszfeld2009point}. Specifically, note that if every row of $(V'+A^\top \rv)_\B $ is non-zero, the gradient of \eqref{prob-lam-new} at $\rv$ is
\begin{equation}\label{grad-geometric-medium}
    \nabla G(\rv)=A(b'+ \operatorname{diag}(\v)(V'+ A ^\top\rv)),
\end{equation}
where $\v\in\R^{n}$ is the vector such that $\v_i=1/\|(V'+ A ^\top\rv)_i\|$ for any $i\in [B]$, and $\v_i=2$ otherwise. By fixing $\v$ and letting $\nabla G(\rv)=0$, the generalized Weiszfeld algorithm (GWA) updates as follows:
\begin{equation}\label{alg-weiszfeld}
  \rv^{k+1}=T(\rv^k) := -( A \operatorname{diag}(\v^k)  A ^\top)^{-1}A(b'+ \operatorname{diag}(\v^k){V'}),\tag{GWA}
\end{equation}
where $\v^k\in\R^n$ denotes the vector just defined above corresponding to $\rv^k$. We assume that every row of $(V'+A^\top \rv^k)_\B $ is non-zero for all $k\geq 0$ so that $\v^k$ is well-defined. Note that 
\begin{equation*}
 \rv^{k+1}-\rv^k=-( A \operatorname{diag}(\v^k)  A ^\top)^{-1}A(b'+ \operatorname{diag}(\v^k)(V'+ A ^\top\rv^k))
=-( A \operatorname{diag}(\v^k)  A ^\top)^{-1}\nabla G(\rv^k),
\end{equation*}
which means we may use $\left\|\rv^{k+1}-\rv^k\right\|_F$ as the residue to measure the optimality of $\rv^{k+1}$. Now we establish the convergence results of \eqref{alg-weiszfeld} in Theorem \ref{thm-retrac}. Denote the sequence generated by \eqref{alg-weiszfeld} as $\{\rv^k\}_{k\geq0}$.

\begin{lemma}\label{lemma-retrac-sub}
    For any given $\rv^0$, if every row of $(V'+A^\top \rv^0)_\B $ is non-zero, then 
    \begin{itemize}
        \item [(i)] $G(T(\rv^0))\leq G(\rv^0)$. The equality holds if and only if $T(\rv^0)=\rv^0$;
        \item [(ii)] $T(\rv)$ is continuous in a neighborhood of $\rv^0$.
    \end{itemize}
    
\end{lemma}

\proof{}
Define an auxiliary function $H(\U,\rv)$ for $\U\in\R^{m\times r},\rv\in\R^{m\times r}$ as follows:
    \begin{multline*}
        H(\U,\rv)=\sum_{i \in[n] \backslash B} {\|(V'+A^\top \U)_i\|^2}+\langle b',A^\top \U\rangle\\
        +\sum_{i\in B}\left(\|(V'+A^\top \rv)_i\|+\cfrac{1}{2\|(V'+A^\top \rv)_i\|}\left(\|(V'+A^\top \U)_i\|^2-\|(V'+A^\top \rv)_i\|^2\right)\right).
    \end{multline*}
    The following properties hold:
    \begin{enumerate}
        \item $H(\rv^0,\rv^0)=G(\rv^0)$.
        \item $H(\U,\rv^0)$ is strongly convex and quadratic w.r.t. $\U$.\quad The strong convexity follows from the fact that the Hessian operator $\nabla_{\U\U}H(\U,\rv^0)=A \operatorname{diag}(\v)  A ^\top\succ 0$, where $\v\in\R^{n}$ is a vector such that $\v_i=1/\|(V'+ A ^\top\rv^0)_i\|$ for any $i\in [B]$, and $\v_i=2$ otherwise.
        \item $G(\U)\leq H(\U,\rv^0)$.\quad The inequality follows from ${y}\leq {x}+(y^2-x^2)/(2{x})$, which holds for any $x>0$ and $y\geq 0$.
        \item $T(\rv^0)=\arg\min_{\U} H(\U,\rv^0)$.\quad Assume that $\U^*=\arg\min_{\U} H(\U,\rv^0)$, then the first order optimality condition $\nabla_\U H(\U^*,\rv^0)=0$ holds, which is equivalent to
    \begin{equation}\label{eq-weisz-grad}
    A(b'+\diag(\v)(V'+A^\top \U^*))=0,    
    \end{equation}
    where $\v\in\R^{n}$ is the vector such that $\v_i=1/\|(V'+ A ^\top\rv^0)_i\|$ for any $i\in [B]$, and $\v_i=2$ otherwise. The equation \eqref{eq-weisz-grad} exactly implies that $\U^*=T(\rv^0)$.
    \end{enumerate}
    Now we can prove the lemma. For 
    (i),
    \begin{equation}\label{weis-ieq}
        G(T(\rv^0))\leq H(T(\rv^0),\rv^0)=\min_{\U}\{H(\U,\rv^0)\}\leq H(\rv^0,\rv^0)=G(\rv^0).
    \end{equation}
    The equality holds only if $T(\rv^0)=\rv^0$ because the optimal solution to $\min_{\U}\{H(\U,\rv^0)\}$ is unique. 
    
    For (ii), consider the nonlinear mapping $F(\U,\rv) := \nabla_\U H(\U,\rv)$. Since every row of $(V'+A^\top \rv^0)_\B $ is non-zero, $F$ is well-defined and smooth in a neighborhood of $(T(\rv^0),\rv^0)$. Since $F(T(\rv^0),\rv^0)=\nabla_\U H(T(\rv^0),\rv^0)=0$ and \[
    \nabla_\U F(\U,\rv^0)=\nabla_{\U\U}H(\U,\rv^0) = A \operatorname{diag}(\v)  A ^\top \succ 0,\] by the implicit function theorem, there exists a unique continuously differentiable function $T':B_\epsilon(\rv^0)\rightarrow \R^{m\times r}$ such that $F(T'(\rv),\rv)=0$ for any $\rv\in B_\epsilon(\rv^0)$. Note that $T(\rv)$ is the unique solution to $\min_\U H(\U,\rv)$, thus $T=T'$ is continuously differentiable in $B_\epsilon(\rv^0)$. 
    
\endproof

\begin{theorem}
    For any regular point $R \in \mathcal{M}_r$, there exists $\epsilon>0$ such that for any $V \in \mathrm{B}_\epsilon(R)$, if the sequence $\{\rv^k\}_{k\geq 0}$ generated by \eqref{alg-weiszfeld} has an accumulation point $\brv$, and every row of $(V'+A^\top \brv)_\B $ and $(V'+A^\top \rv^k)_\B $ is non-zero for all $k\geq 0$, then $\brv$ is a global minimum of $G(\rv)$.
\end{theorem}
\proof{}
    By (ii) of Theorem \ref{thm-retrac}, there exists $\epsilon>0$ such that for any $V \in \mathrm{B}_\epsilon(R)$, the optimal solution to \eqref{prob-lam-new} exists, so $G(\rv^k)$ is bounded below. Also, by (i) of Lemma \ref{lemma-retrac-sub}, $G(\rv^k)$ is non-increasing. Hence, $\lim_{k\rightarrow \infty}\{G(\rv^k)\}$ exists, and 
    \begin{equation}\label{eq-conv-weisz-1}
        \lim_{k\rightarrow \infty}\{G(T(\rv^k))-G(\rv^k)\}=0.
    \end{equation}
    Since every row of $(V'+A^\top \brv)_\B $ is non-zero, by (ii) of Lemma \ref{lemma-retrac-sub}, $T$ is continuous in a neighborhood of  $\brv$. Since $\brv$ is an accumulation point, we can assume that $\lim_{n\rightarrow \infty}\rv^{k_n}=\brv$. Then by the continuity of $G$ and \eqref{eq-conv-weisz-1} we have
    \begin{align*}
        G(T(\brv))=G(T(\lim_{n\rightarrow \infty}
        \rv^{k_n}))=\lim_{n\rightarrow \infty} G(T(\rv^{k_n}))=\lim_{n\rightarrow \infty} G(\rv^{k_n})= G(\lim_{n\rightarrow \infty}
        \rv^{k_n})=G(\brv).
    \end{align*}
    By (i) of Lemma \ref{lemma-retrac-sub}, we further have 
    \begin{equation}\label{eq-weisz-bar-lambda}
        \brv=T(\brv)=-( A \operatorname{diag}(\bv)  A ^\top)^{-1}A(b'+ \operatorname{diag}(\bv){V'}),
    \end{equation}
    where $\bv\in\R^{n}$ is the vector such that $\bv_i=1/\|(V'+ A ^\top\brv)_i\|$ for any $i\in [B]$, and $\bv_i=2$ otherwise. By multiplying $A \operatorname{diag}(\bv) A^\top$ on both sides of \eqref{eq-weisz-bar-lambda} and simplification, we get that $\nabla G(\brv)=0$. Thus, $\brv$ is a global minimum $\brv$ of $G(\rv)$. 
\endproof

\begin{remark}
    Theorem \ref{thm-retrac} assumes that \eqref{alg-weiszfeld} will not pass through any anchor point, i.e. every row of $(V'+A^\top \rv)_\B $ is non-zero for every $\rv\in\{\rv^k\}_{k\geq 0}$ and $\brv$. In practice, \eqref{alg-weiszfeld} never reaches an anchor point when $\M_r$ is locally smooth. Theoretically, we can remove the assumption by using a modified version of \eqref{alg-weiszfeld}, see for example, \cite{ostresh1978convergence,vardi2001modified,beck2015weiszfeld}. The idea of the modified algorithm is that we can always escape from the encountered anchor point along the unit direction with the smallest directional derivative to decrease the function value.
\end{remark}

\begin{remark}
    We can also use the Newton method to solve \eqref{prob-lam-new} as follows:
\begin{equation}
    \rv^{k+1}=\rv^k-( A \mathcal{P} A ^\top)^{-1}A(b'+ \operatorname{diag}(\v^k)(V'+ A ^\top\rv^k)),\tag{Newton}
\end{equation}
where $\mathcal{P}:\R^{n\times r}\rightarrow\R^{n\times r}$ is the linear operator defined as
\begin{equation*}
    (\mathcal{P}(X))_i=\begin{cases}
         X_i&i\notin B \\
         X_i\left(\cfrac{I}{\|(V'+A^\top\rv)_i\|}-\cfrac{(V'+A^\top\rv)_i^\top(V'+A^\top\rv)_i}{\|(V'+A^\top\rv)_i\|^3} \right)&i\in B.
    \end{cases}
\end{equation*}
Due to the fast local convergence, in practice, we can switch to the Newton method when the residue $\|\rv^{k+1}-\rv^k\|$ in \eqref{alg-weiszfeld} is small to accelerate the convergence of GWA.
\end{remark}

\section{Numerical experiments}\label{sec-numerical}
In this section, we conduct numerical experiments to demonstrate the effectiveness of RNNAL in solving the DNN problem \eqref{prob-P-dnn-YZ}. All the experiments are run using {\sc Matlab} R2023b on a workstation with Intel
Xeon E5-2680 (v3) cores with 96GB RAM.
\vspace{6pt}

\noindent\textbf{Baseline Solvers}. We compare the performance of RNNAL with another ALM-based algorithm RNNAL-Diag and the solver SDPNAL+ \cite{SDPNAL,SDPNALp1,SDPNALp2}. RNNAL-Diag uses the same framework as RNNAL but also penalizes all the equality constraints in $\F$, except $\diag_\B (X)=x_\B $. Thus, the feasible set of the RNNAL-Diag subproblem after BM factorization is the simple spherical manifold $\Br$ as defined in \eqref{Br}. RNNAL-Diag is a representative of many low-rank SDP algorithms that require the smoothness of the feasible set after direct BM factorization. A comparison of different ALM-based algorithms is provided in Table \ref{tab:alm-summary}. The reason for not using an algorithm in the last row of Table \ref{tab:alm-summary} is that the rank of $R$ must be fixed due to the fixed dimension of the multiplier of $AR=be_1^\top$. There are also other low-rank SDP solvers. The reason for not using ManiSDP \cite{wang2023solving}, HALLaR \cite{monteiro2024low}, or LoRADS \cite{han2024low} is that they cannot handle the {nonnegativity} constraint $Y\geq 0$. The reason for not using SDPDAL \cite{wang2023decomposition} is that the code is unavailable and the framework is similar to RNNAL-Diag. The reason for not using SDPLR \cite{BM1} or SketchyCGAL \cite{yurtsever2021scalable} is that their performance can not measure up to RNNAL-Diag during our preliminary tests.
\begin{center}
\begin{table}[h!]
\begin{center}
    \def\arraystretch{1.1}
    \begin{tabular}{|l|c|c|c|c|} 
     \hline
     \multicolumn{1}{|c|}{penalty term} & manifold  & algorithm & issue  \\ 
     \hline
     $Y\geq 0$, $X_E=0$ & $\N_r$ &  - & nonsmooth  \\ 
     \hline
     $Y\geq 0$, $X_E=0$ & $\M_r$ &  RNNAL & -  \\ 
     \hline
     $Y\geq 0$, $X_E=0$, $Ax=b$, $AX=bx^\top$ & $\Br$  & RNNAL-Diag & -  \\
     \hline
     $Y\geq 0$, $X_E=0$, $Ax=b$, $AX=bx^\top$, $\diag_\B (X)=x_\B $ & $\R^{n\times r}$ & SDPLR & slow \\
     \hline
     $Y\geq 0$, $X_E=0$, $AR=be_1^\top$, $\diag_\B (X)=x_\B $ & $\R^{n\times r}$ & - & fixed rank\\
     \hline
    \end{tabular}
    \caption{Comparison of different ALM-based methods.}
    \label{tab:alm-summary}
\end{center}
\end{table}
\end{center}

\noindent\textbf{Stopping Conditions}. Based on the KKT conditions \eqref{KKT-DNN} for \eqref{prob-P-dnn-YZ}, we define the following relative KKT residue to measure the accuracy of the solution obtained by RNNAL and RNNAL-Diag:
\begin{align*}
\operatorname{R_p} := \max\left\{\cfrac{\|\mathcal{A}(Y)- {d} \|}{1+\| {d} \|},\cfrac{\|Y-Z\|}{1+\|Y\|+ \|Z\|}\right\},\ 
\operatorname{R_d} :=  \frac{\|\Pi_{\mathcal{K}}(-S)\|}{1+\|S\|},\
\operatorname{R_{c}} := \frac{|\langle Y, S\rangle|}{1+\|Y\|+\|S\|} .
\end{align*}
For a given tolerance $\tol> 0$, we terminate RNNAL and RNNAL-Diag when the maximum residue $\operatorname{R_{max}}:=\max\{\operatorname{R_p},\operatorname{R_d},\operatorname{R_c}\}<\tol$ or the maximum time $\timelimit$ is reached. We choose $\tol=10^{-6}$ and $\timelimit=3600\,\tt{(secs)}$ for all solvers in our experiments.
\vspace{6pt}

\noindent\textbf{Implementation}. 
For RNNAL and RNNAL-Diag, we use a Riemannian gradient descent method with BB step and non-monotone line search to solve the augmented Lagrangian subproblems (see \cite{iannazzo2018riemannian,gao2021riemannian,tang2023feasible}). The penalty parameter $\sigma_k$ is initialized as $\sigma_0=1$ and increased by a factor of $1.25$ if the 
reduction in the primal infeasibility $\operatorname{R_p}$ is not significant enough. {The initial rank $r_0$ is set to $\min\{200, \lceil n/5 \rceil\}$ by default. For rank adjustments, the rank is increased by the default value of $\tau = 1$ to escape saddle points using Armijo line search and is decreased following the procedure outlined in \cite{tang2023feasible}. However, choosing problem-specific values for $r_0$ and $\tau$ may lead to improved performance. We should mention
that the rank adjustment is performed only once after solving 
the ALM subproblem in each ALM iteration.
} The initial point $R_0$ is randomly selected from the feasible region $\M_{r_0}$ for RNNAL and $\B_{r_0}$ for RNNAL-Diag.
\vspace{6pt}

\noindent\textbf{Table Notations.} We use `-' to indicate that an algorithm does not reach the required tolerance $\tol$ within the given maximum time of $\timelimit$. For SDPNAL+, we use ``it'', ``itsub'', ``itA'', which are the same notations as \cite{SDPNALp1}, to denote the number of outer iterations, the total number of semismooth Newton inner iterations, the total number of iterations for ADMM+, respectively. For RNNAL and RNNAL-Diag, we use ``it'', ``itsub'', ``$r$'' to denote the number of ALM iterations, the total number of Riemannian gradient descent iterations, the final rank of the output matrix $R$, respectively. The objective function is denoted by ``obj''. The computation time (in seconds) is reported in the last column of the table.

\subsection{Quadratic assignment problems}\label{subsec-qap}
\begin{figure}[h!]
  \begin{minipage}[h]{.5\linewidth}
    \centering
    \includegraphics[width=0.8\textwidth]{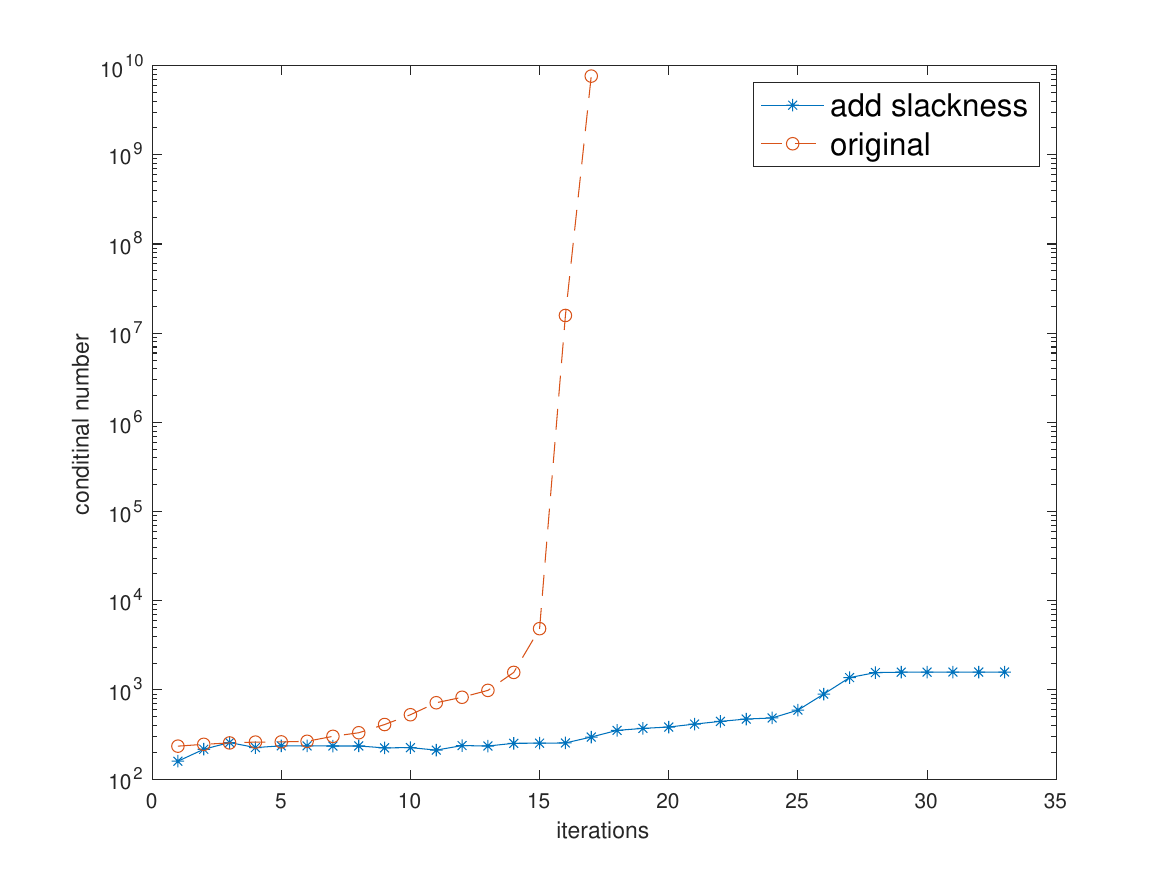}%
    \caption
      {%
        Comparison of the condition number of $h_Rh_R^*$ for the chr12a dataset.%
        \label{fig:cond-number}%
      }%
  \end{minipage}\hfill
  \begin{minipage}[h]{.4\linewidth}
    \centering
    \small
    \begin{tabular}{clcc}
    \toprule
    problem & algorithm & time & $\operatorname{R_{max}}$ \\
    \midrule
    \multirow{2}{*}{\eqref{prob-qap-dnn}}&
     RNNAL    & -     & -   \\
    &SDPNAL+ & 20.1 & 2.1e-7\vspace{4pt}\\
    \multirow{2}{*}{\eqref{prob-qap-dnn-v2}}&
     RNNAL    & 16.7 & 1.7e-7\\
    &SDPNAL+  & 66.9 & 3.2e-7\\
      \bottomrule
    \end{tabular}
    \captionof{table}
      {%
        Time comparison of RNNAL and SDPNAL+ for different formulations.%
        \label{tab:qap}%
      }
  \end{minipage}
\end{figure}

In this subsection, we use the quadratic assignment problem (QAP) as an example to show that RNNAL can avoid non-smoothness by applying the constraint-relaxation strategy proposed in Subsection \ref{subsec-avoid-non-reuglar}. Let $\Pi$ be the set of $p\times p$ permutation matrices. Given matrices $W,D\in\S^p$, the QAP is given by
\begin{equation}\label{prob-QAP-v1}
\min\left\{\langle Y,WYD \rangle :\ Y\in \Pi\right\}.
\end{equation}
Denote $n\coloneqq p^2,\ x:=\operatorname{vec}(\Pi)\in\R^{n}$, $Q:=D \otimes W\in \R^{n\times n}$ and 
$A=(e_p\otimes I_p, I_p\otimes e_p)^\top \in\R^{2p\times n}$.
Since $\Pi=\{Y\in \{0,1\}^{p\times p}:  Ye=e, Y^\top e=e\}$, \eqref{prob-QAP-v1} can be equivalently expressed in the form of \eqref{prob-MBQP} as follows:
\begin{equation}\label{prob-qap-v2}
    \min\left\{
    x^\top Q x: \
    Ax=e,\ x\in \{0,1\}^{n}
    \right\}.
\end{equation}
The corresponding DNN relaxation in the form of \eqref{prob-dnn-new} is
\begin{equation}\label{prob-qap-dnn}
    \min\left\{\<Q,X\>:\ Ax=e, AX=ex^\top , \diag(X)=x, \begin{pmatrix}1& x^\top \\ x & X\end{pmatrix}\in \mathbb{S}_{+}^{n+1} \cap \mathbb{N}^{n+1}\right\}.
\end{equation}
The algebraic variety $\M_r$ of \eqref{prob-qap-dnn} is given by
\eqref{algebraic-variety-M} with $b=e$ and $\B=[n]$.
When the DNN relaxation \eqref{prob-qap-dnn} is tight, a rank-one solution exists and is the exact solution to \eqref{prob-MBQP}. Applying RNNAL directly to solve \eqref{prob-qap-dnn} may lead to numerical issues since $\M_1$ is non-smooth as shown in Example \ref{example-rank-1-nonsmooth}. However, the non-smoothness can be avoided by applying the constraint-relaxation strategy. We first define  
\begin{equation*}
    \bar Q:=\begin{pmatrix}
        Q&0_{n\times 4p}\\
        0_{4p\times n}&0_{4p\times 4p}
    \end{pmatrix},\  \bar A:=\begin{pmatrix}
        A&I_{2p}&0_{2p\times 2p}\\
        A&0_{2p\times 2p}&-I_{2p}
    \end{pmatrix},\ B:=[n].
\end{equation*}
Then problem \eqref{prob-qap-v2} can be equivalently written as 
\begin{equation*}
    \min\left\{
    x^\top Q x: \
    \bar A\begin{pmatrix}
        x\\s
    \end{pmatrix}=e, x\in \{0,1\}^{n},\ \begin{pmatrix}
        x\\s
    \end{pmatrix}\in\R^{n+4p}_+
    \right\}.
\end{equation*}
The corresponding DNN relaxation is
\begin{equation}\label{prob-qap-dnn-v2}
    \min\left\{\< \bar Q,X\>:\ \bar A\bar x=e, \bar AX=e\bar x^\top , \diag_\B (X)=\bar x_\B , \begin{pmatrix}1& \bar x^\top \\ \bar x & X\end{pmatrix}\in \mathbb{S}_{+}^{n+4p+1} \cap \mathbb{N}^{n+4p+1}\right\}.
\end{equation}
According to Lemma \ref{lemma-DNN-equivalent}, problem \eqref{prob-qap-dnn} and \eqref{prob-qap-dnn-v2} are equivalent. Moreover, Lemma \ref{lemma-avoid-nonsmooth} ensures that the new manifold $\M'_r$ of \eqref{prob-qap-dnn-v2} is smooth for any positive integer $r$. Thus, we can apply RNNAL to solve the new DNN problem \eqref{prob-qap-dnn-v2}. To demonstrate the effectiveness of our avoiding-non-smoothness strategy, we use RNNAL to solve problems \eqref{prob-qap-dnn} and \eqref{prob-qap-dnn-v2} and compare the condition number of $h_Rh_R^*$ to assess the level of smoothness. The {\tt chr12a} dataset with $p=12$ from the QAP Library \cite{qaplib} is selected as an example due to the tightness of its DNN relaxation \eqref{prob-qap-dnn}. Similar behaviors are observed for other datasets with exact DNN relaxations. The results are shown in Figure \ref{fig:cond-number} and Table \ref{tab:qap}. 

In Figure \ref{fig:cond-number}, the $x$-axis represents the ALM outer iteration number. For each ALM subproblem, $h_Rh_R^*$ is computed several times, but we only report the one with the maximum condition number. It can be observed that the condition number for the original DNN problem \eqref{prob-qap-dnn} (orange curve) increases rapidly near the optimal solution. Conversely, the condition number for the equivalent DNN problem \eqref{prob-qap-dnn-v2} remains relatively small. Table \ref{tab:qap} shows that RNNAL failed to solve \eqref{prob-qap-dnn} due to the singularity issues in the projection and retraction near the rank-one solution. However, RNNAL successfully achieved the required accuracy with the new reformulation \eqref{prob-qap-dnn-v2}. This reformulation increases the variable dimension and the number of constraints, causing SDPNAL+ to take longer to converge. 
For the computational results on more instances from the QAP Library, we refer the readers to Appendix \ref{appendix-QAP}. It is important to note that the speed comparison between RNNAL and SDPNAL+ is not the primary focus here, as the solutions to the QAP instances are typically of high rank. Consequently, SDPNAL+ is anticipated to exhibit superior computational speed compared to RNNAL for these instances.

\subsection{Binary integer nonconvex quadratic programming}

The binary integer nonconvex quadratic (BIQ) programming is a special case of \eqref{prob-MBQP} without linear constraints:
\begin{equation}
\min\left\{x^{\top} Q x+2 c^{\top}x:\ x \in\{0,1\}^n\right\}.
\end{equation}
The corresponding DNN relaxation in the form of \eqref{prob-dnn-new} is
\begin{equation} \label{DNN-BIQ}
\min\left\{\<Q,X\>+2 c^{\top}x:\ \diag(X)=x,\ \begin{pmatrix}1& x^\top \\ x & X\end{pmatrix}\in \mathbb{S}_{+}^{n+1} \cap \mathbb{N}^{n+1}\right\}.
\end{equation}
The algebraic variety $\mathcal{M}_r$ 
corresponding to \eqref{DNN-BIQ} is the oblique manifold:
\begin{equation}\label{OB-manifold}
\mathcal{O B}(n, r) := \left\{R \in \mathbb{R}^{n \times r}:\diag(RR^\top)=e \right\}.
\end{equation}
We choose the test data for $Q$ and $c$ from the library ORLIB\footnote{Dataset from \url{http://people.brunel.ac.uk/~mastjjb/jeb/info.html}.} maintained by J.E. Beasley with problem dimension $n\in\{1000,2500\}$. Each dimension has ten instances, but we only report one in this subsection because the performance on other instances is similar. {The complete results for all instances are available in Appendix \ref{appendix-BIQ}.} Due to the lack of available data on larger dimensions, we randomly generate two sets of data with problem dimensions $n\in\{5000,10000\}$ following the same generation procedure proposed by J.E. Beasley in \cite{beasley1998heuristic}. Since RNNAL-Diag is equivalent to RNNAL for solving BIQ problems, we present results for RNNAL and SDPNAL+ only. 

\begin{tiny}
\begin{longtable}[c]{clrrrccccc}
\caption{Computational results for BIQ problems. \label{tab-BIQ}} \\
\toprule
problem  & algorithm & it & itsub & $r$/itA & $\operatorname{R_p}$ & $\operatorname{R_d}$ & $\operatorname{R_c}$ & obj & time   \\
\midrule
\endfirsthead

\multicolumn{10}{c}%
{{ Table \thetable\ continued from previous page}} \\
\toprule
problem  & algorithm & it & itsub & $r$/itA & $\operatorname{R_p}$ & $\operatorname{R_d}$ & $\operatorname{R_c}$ & obj & time   \\
\midrule
\endhead
\midrule
\multicolumn{10}{r}{{Continued on next page}} \\
\midrule
\endfoot

\bottomrule
\endlastfoot
$n=1000$ & RNNAL & 11 & 911 & 55 & 9.04e-07 & 3.17e-07 & 3.80e-07 & -3.9849472e+05 & 9.40e+00 \\ 
 & SDPNAL+ & 118 & 172 & 2752 & 8.92e-07 & 9.69e-07 & 5.54e-07 & -3.9849494e+05 & 2.61e+02 \\ [3pt] 

$n=2500$ & RNNAL & 9 & 817 & 97 & 6.11e-07 & 4.53e-07 & 5.57e-08 & -1.6354913e+06 & 1.03e+02 \\ 
 & SDPNAL+ & - & - & - & - & - & - & - & - \\ [3pt] 

$n=5000$& RNNAL & 6 & 1091 & 147 & 5.52e-07 & 9.46e-08 & 1.51e-07 & -4.7435656e+06 & 5.61e+02 \\ 
& SDPNAL+ & - & - & - & - & - & - & - & - \\ [3pt] 

$n=10000$& RNNAL & 5 & 1078 & 216 & 4.80e-07 & 1.80e-07 & 9.93e-07 & -1.3832829e+07 & 2.32e+03 \\ 
& SDPNAL+ & - & - & - & - & - & - & - & - \\

\end{longtable}
\end{tiny}

From Table \ref{tab-BIQ}, we can see that RNNAL can solve all the problems to the required accuracy, while SDPNAL+ fails to solve problems with dimension $n\geq 2500$ within the 1-hour time limit. RNNAL is about 30 times faster than SDPNAL+ for medium size problems. Moreover, RNNAL can handle problems with dimensions up to $n=10000$ in 40 minutes. Such observations show the effectiveness of RNNAL and its potential to solve large-scale BIQ problems. {Since any max-cut problem can be formulated as a BIQ problem, we have also tested RNNAL on max-cut instances. Preliminary results indicate that RNNAL is as effective for max-cut problems as it is for BIQ problems.}

{We observe that the solution ranks are typically around $50-200$, which are not particularly small. Additionally, computing the projection and retraction on the oblique manifold is straightforward and computationally efficient. Consequently, the majority of the computational time is spent on evaluating the objective function $f_r(R)$ and its gradient $\nabla f_r(R)$. A small portion of the time (approximately $10\%-20\%$) is spent on computing the smallest eigenvalue and its corresponding eigenvector of the dual variable $S$, which is used to assess dual infeasibility and identify the direction to escape from saddle points. This portion increases as the problem dimension grows.
}

\subsection{Maximum stable set problems}\label{subsec-theta}
Given a graph $G$ with $n$ vertices and $m$ edges, denote its edge set as $E\subseteq \{(i,j)\mid 1\leq i<j\leq n \}$. The maximum stable set problem is given as follows:
\begin{equation}\label{prob-theta-v1}
\max\left\{x^{\top} x:\ x_ix_j=0, \ \forall (i,j)\in E,\  x \in\{0,1\}^n\right\}.
\end{equation}
The DNN relaxation $\theta_+(G)$ of \eqref{prob-theta-v1} in the form of \eqref{prob-dnn-new} is
\begin{equation}\label{prob-theta-dnn}
    \theta_+(G)=\min\left\{-\<I,X\>:\ \diag(X)=x,\ X_{ij}=0,\ \forall (i,j)\in E,\ \begin{pmatrix}1& x^\top \\ x & X\end{pmatrix}\in \mathbb{S}_{+}^{n+1} \cap \mathbb{N}^{n+1}\right\},
\end{equation}
which has $m+n+1$ equality constraints. The algebraic variety $\M_r$ corresponding to \eqref{prob-theta-dnn} is the oblique manifold $\mathcal{O B}(n, r)$ defined in \eqref{OB-manifold}. We select large sparse graphs from Gset\footnote{Dataset from \url{https://web.stanford.edu/~yyye/yyye/Gset/}.}. The complete results for all graphs are available in Appendix \ref{appendix-theta}. Since RNNAL-Diag is equivalent to RNNAL for solving $\theta_+$ problems, we present results for RNNAL and SDPNAL+ only. 

\begin{tiny}
\begin{longtable}[c]{clrrrccccc}
\caption{Computational results for $\theta_+$ problems. \label{tab-theta}} \\
\toprule
problem $(n,m)$  & algorithm & it & itsub & $r$/itA & $\operatorname{R_p}$ & $\operatorname{R_d}$ & $\operatorname{R_c}$ & obj & time   \\
\midrule
\endfirsthead

\multicolumn{10}{c}%
{{ Table \thetable\ continued from previous page}} \\
\toprule
problem  & algorithm & it & itsub & $r$/itA & $\operatorname{R_p}$ & $\operatorname{R_d}$ & $\operatorname{R_c}$ & obj & time   \\
\midrule
\endhead
\midrule
\multicolumn{10}{r}{{Continued on next page}} \\
\midrule
\endfoot

\bottomrule
\endlastfoot

G43& RNNAL & 10 & 795 & 80 & 7.48e-07 & 4.18e-08 & 2.09e-08 & -2.7973625e+02 & 9.61e+00 \\ 
(1000,9990)  & SDPNAL+ & 48 & 61 & 1250 & 4.58e-07 & 6.61e-07 & 9.64e-14 & -2.7973595e+02 & 1.43e+02 \\ [3pt] 

G34& RNNAL & 10 & 2188 & 11 & 3.79e-07 & 1.48e-08 & 5.58e-07 & -9.9999198e+02 & 9.91e+01 \\ 
(2000,4000)& SDPNAL+ & - & - & - & - & - & - & - & - \\ [3pt] 

G48& RNNAL & 8 & 1497 & 21 & 9.25e-07 & 4.67e-08 & 3.09e-07 & -1.4999238e+03 & 2.44e+02 \\ 
(3000,6000)& SDPNAL+ & - & - & - & - & - & - & - & - \\ [3pt] 

G55& RNNAL & 20 & 2130 & 353 & 9.96e-07 & 1.85e-07 & 4.84e-08 & -2.3230485e+03 & 1.45e+03 \\ 
(5000,12498)& SDPNAL+ & - & - & - & - & - & - & - & - \\ 

\end{longtable}
\end{tiny}

From Table \ref{tab-theta}, we can see that RNNAL can solve all the problems to the required accuracy, while SDPNAL+ fails to solve problems with dimension $n\geq 2000$ within the 1-hour time limit. RNNAL is over 10 times faster than SDPNAL+ for all problems. In particular, RNNAL is faster than SDPNAL+ by at least a factor of 30 times for the G34 problem. {The time-consuming steps are similar to those encountered in BIQ problems.}

\subsection{Quadratic knapsack problems}\label{subsec-QKP}
The binary quadratic knapsack problem (QKP), introduced in \cite{gallo1980quadratic},  is as follows:
\begin{equation}\label{prob-qkp-v1}
\max\left\{x^{\top} Q x:\ a^\top x\leq \tau,\ x \in\{0,1\}^n\right\},
\end{equation}
where $Q\in \S^n$ is a nonnegative profit matrix, $a\in \R^n_{++}$ is the weight vector and $\tau>0$ is the knapsack capacity. To derive the DNN relaxation of \eqref{prob-qkp-v1}, we convert the inequality constraint to an equality constraint, and then focus on the new problem:
\begin{equation}\label{prob-qkp-v2}
\max\left\{x^{\top} Q x:\ a^\top x= \tau,\ x \in\{0,1\}^n\right\}.
\end{equation}
Problem \eqref{prob-qkp-v1} and \eqref{prob-qkp-v2} are not equivalent in general. 
However, under mild conditions, {SDP} relaxations of \eqref{prob-qkp-v1} and \eqref{prob-qkp-v2} are equivalent, see \cite{tang2024feasible} for more details. When $a=e$ and $\tau=k$, \eqref{prob-qkp-v2} reduces to the $k$-subgraph problem. The DNN relaxation of \eqref{prob-qkp-v2} in the form of \eqref{prob-dnn-new} is
\begin{equation}\label{prob-qkp-dnn}
    \min\left\{-\<Q,X\>:\ a^\top x=\tau,\ a^\top X=\tau x^\top,\ \diag(X)=x,\ \begin{pmatrix}1& x^\top \\ x & X\end{pmatrix}\in \mathbb{S}_{+}^{n+1} \cap \mathbb{N}^{n+1}\right\},
\end{equation}
which has $2n+2$ equality constraints.
For \eqref{prob-qkp-dnn}, the feasible set $\M_r$ \eqref{mani-Mr} corresponding to the low-rank factorization form of the ALM subproblem in \eqref{prob-Rie-sub} is given by
\begin{equation}\label{manifold-QKP}
    \M_r=\left\{R\in \R^{n\times r}:\ a^\top R=\tau e_1^\top,\diag (RR^\top)-Re_1=0 \right\}.
\end{equation}
We randomly generate the profit matrix $Q$ and weight vector $a$ following the procedure proposed by Gallo et al. in \cite{gallo1980quadratic}, which has been widely used in the literature (see, for example, \cite{caprara1999exact,billionnet2004exact,pisinger2007quadratic,tang2024feasible}). The entries of the profit matrix $Q_{ij}=Q_{ji}$ are integers randomly generated uniformly in the range $[1, 100]$ with probability $p$ and zero otherwise. The elements of the coefficient vector $a$ are randomly selected integers in the range $[1, 50]$. The knapsack capacity is chosen to be $0.9 \cdot e^{\top} a$. The probability $p$ is chosen in $\{0.1, 0.5, 0.9\}$. The dimension $n$ is chosen in $\{500,1000,5000,10000\}$. 

\begin{tiny}
\begin{longtable}[c]{clrrrccccc}
\caption{Computational results for QKP problems. \label{tab-QKP}} \\
\toprule
 problem  & algorithm & it & itsub & $r$/itA & $\operatorname{R_p}$ & $\operatorname{R_d}$ & $\operatorname{R_c}$ & obj & time   \\
\midrule
\endfirsthead

\multicolumn{10}{c}%
{{ Table \thetable\ continued from previous page}} \\
\toprule
problem & algorithm & it & itsub & $r$/itA & $\operatorname{R_p}$ & $\operatorname{R_d}$ & $\operatorname{R_c}$ & obj & time   \\
\midrule
\endhead
\midrule
\multicolumn{10}{r}{{Continued on next page}} \\
\midrule
\endfoot
\bottomrule
\endlastfoot

$n=500$ & RNNAL & 7 & 256 & 12 & 9.10e-07 & 5.21e-07 & 1.82e-08 & -1.1421503e+06 & 1.85e+00 \\ 
$p=0.1$ & SDPNAL+ & 56 & 77 & 2940 & 6.91e-07 & 4.79e-07 & 4.92e-14 & -1.1421502e+06 & 9.99e+01 \\ 
& RNNAL-Diag & 60 & 171299 & 19 & 8.93e-07 & 5.52e-07 & 2.25e-16 & -1.1421487e+06 & 4.78e+02 \\ [3pt] 

$n=500$ & RNNAL & 7 & 197 & 17 & 8.71e-07 & 3.31e-07 & 1.22e-09 & -5.6677125e+06 & 1.77e+00 \\ 
$p=0.5$ & SDPNAL+ & 67 & 139 & 3864 & 1.02e-06 & 7.41e-07 & 6.65e-07 & -5.6677127e+06 & 1.48e+02 \\ 
& RNNAL-Diag & 50 & 189687 & 24 & 9.48e-07 & 9.28e-07 & 2.32e-16 & -5.6677064e+06 & 5.41e+02 \\ [3pt] 

$n=500$ & RNNAL & 4 & 184 & 21 & 6.15e-07 & 3.34e-07 & 2.57e-08 & -1.0261057e+07 & 1.47e+00 \\ 
$p=0.9$ & SDPNAL+ & 56 & 102 & 3375 & 6.79e-07 & 9.82e-07 & 5.01e-07 & -1.0261059e+07 & 1.19e+02 \\ 
& RNNAL-Diag & 50 & 163983 & 19 & 4.83e-07 & 4.16e-07 & 2.50e-15 & -1.0261052e+07 & 4.49e+02 \\ [3pt] 

$n=1000$ & RNNAL & 7 & 341 & 24 & 9.54e-07 & 9.64e-07 & 5.09e-09 & -4.6071170e+06 & 1.16e+01 \\ 
$p=0.1$ & SDPNAL+ & 79 & 138 & 4862 & 2.75e-07 & 9.41e-07 & 1.53e-06 & -4.5995222e+06 & 7.44e+02 \\ 
& RNNAL-Diag & - & - & - & - & - & - & - & - \\ [3pt] 

$n=1000$ & RNNAL & 4 & 154 & 36 & 9.64e-07 & 2.72e-07 & 1.81e-09 & -2.2759420e+07 & 5.96e+00 \\ 
$p=0.5$ & SDPNAL+ & 97 & 295 & 6220 & 8.72e-07 & 9.32e-07 & 8.14e-09 & -2.2747378e+07 & 1.05e+03 \\ 
& RNNAL-Diag & - & - & - & - & - & - & - & - \\ [3pt]

$n=1000$ & RNNAL & 4 & 198 & 24 & 5.23e-07 & 2.17e-09 & 2.02e-09 & -4.0961669e+07 & 6.79e+00 \\ 
$p=0.9$ & SDPNAL+ & 113 & 310 & 7179 & 8.21e-07 & 9.83e-07 & 2.91e-07 & -4.0934695e+07 & 1.13e+03 \\ 
& RNNAL-Diag & - & - & - & - & - & - & - & - \\ [3pt]

$n=5000$ & RNNAL & 3 & 379 & 130 & 6.78e-07 & 4.63e-07 & 5.47e-09 & -1.1383194e+08 & 2.55e+02 \\ 
$p=0.1$ & SDPNAL+ & - & - & - & - & - & - & - & - \\ 
& RNNAL-Diag & - & - & - & - & - & - & - & - \\ [3pt]

$n=5000$ & RNNAL & 2 & 242 & 97 & 5.27e-07 & 1.19e-07 & 1.51e-10 & -5.6708048e+08 & 1.64e+02 \\ 
$p=0.5$ & SDPNAL+ & - & - & - & - & - & - & - & - \\ 
& RNNAL-Diag & - & - & - & - & - & - & - & - \\ [3pt]

$n=5000$ & RNNAL & 1 & 177 & 143 & 7.97e-07 & 4.54e-08 & 1.96e-09 & -1.0209621e+09 & 1.21e+02 \\ 
$p=0.9$ & SDPNAL+ & - & - & - & - & - & - & - & - \\ 
& RNNAL-Diag & - & - & - & - & - & - & - & - \\ [3pt]

$n=10000$ & RNNAL & 2 & 414 & 202 & 5.02e-07 & 1.46e-07 & 7.28e-10 & -4.5412188e+08 & 1.16e+03 \\ 
$p=0.1$ & SDPNAL+ & - & - & - & - & - & - & - & - \\ 
& RNNAL-Diag & - & - & - & - & - & - & - & - \\ [3pt]

$n=10000$ & RNNAL & 1 & 197 & 202 & 6.78e-07 & 1.09e-07 & 1.29e-10 & -2.2698395e+09 & 5.62e+02 \\ 
$p=0.5$ & SDPNAL+ & - & - & - & - & - & - & - & - \\ 
& RNNAL-Diag & - & - & - & - & - & - & - & - \\ [3pt]

$n=10000$ & RNNAL & 1 & 226 & 202 & 3.36e-07 & 3.41e-08 & 1.73e-10 & -4.0846437e+09 & 6.33e+02 \\ 
$p=0.9$ & SDPNAL+ & - & - & - & - & - & - & - & - \\ 
& RNNAL-Diag & - & - & - & - & - & - & - & - \\ 

\end{longtable}
\end{tiny}

As shown in Table \ref{tab-QKP}, RNNAL is significantly faster than RNNAL-Diag and SDPNAL+ for most cases. For some instances, RNNAL can be more than 200 times faster than SDPNAL+ and 1000 times faster than RNNAL-Diag. Moreover, RNNAL can solve large QKP problems with $n=10000$ in 10 to 20 minutes. {Additionally, the ranks of the solutions of the large DNN 
instances with $n\geq 5000$ are roughly in the range of 100--200, which is not considered as small. Note that computing the projection and retraction on the manifold $\M_r$ is generally not straightforward and can be computationally expensive. However, by employing the strategies outlined in Subsections \ref{subsec-proj} and \ref{subsec-retrac}, the time spent on these operations is reduced to approximately $10\%$ of the total time. This small proportion 
of time spent highlights the efficiency of our proposed method for computing the projection and retraction.}

\subsection{Disjunctive quadratic knapsack problems}\label{subsec-DQKP}

We consider the disjunctive quadratic knapsack problem (DQKP) introduced in \cite{yamada2002heuristic,saracc2014generalized} as follows:
\begin{equation}\label{prob-DQKP-v1}
\max\left\{x^{\top} Q x:\ a^\top x\leq \tau,\ x_ix_j=0,\ (i,j)\in E,\ x \in\{0,1\}^n\right\},
\end{equation}
where $E\subseteq \{(i,j)\mid 1\leq i<j\leq n \}$ denotes the set of incompatible pairs, and other notations are the same as those in subsection \ref{subsec-QKP}. Similar to QKP, {by replacing the inequality constraint with an equality constraint,} the DNN relaxation of \eqref{prob-DQKP-v1} in the form of \eqref{prob-P-dnn-YZ} is
\begin{equation}\label{prob-DQKP-dnn}
    \min\left\{-\<Q,X\>:\ Y-Z=0,\ Y\in \F\cap \K,\ Z\in \mathcal{P} \right\},
\end{equation}
where $\K=\S^{n+1}_+$,
$\mathcal{P}=\mathcal{Z}\cap \mathbb{N}^{n+1}$ with $\mathcal{Z}$ defined as in Section~\ref{subsec-DNN}, and 
\begin{align*}
    &\F:=\left\{\begin{pmatrix}
        1&x^\top\\
        x&X
    \end{pmatrix}\in \S^{n+1}:\ a^\top x=\tau,\ a^\top X=\tau x^\top,\ \diag(X)=x  \right\}.
\end{align*}
Compared to the DNN relaxation \eqref{prob-qkp-dnn} for QKP, \eqref{prob-DQKP-dnn} includes $|E|$ additional equality constraints in ${\cal P}$. The feasible set $\M_r$ \eqref{mani-Mr} corresponding to the low-rank factorization form of the ALM subproblem in \eqref{prob-Rie-sub} is the same as \eqref{manifold-QKP}.

We randomly generate the profit matrix $Q$ and weight vector $a$ following the same procedure as QKP in Subsection \ref{subsec-QKP}. The probability $p$ is chosen to be 0.9. For a given conflict density ratio $d$, the disjunctive set $E$ corresponds to the edge index of a randomly generated graph $G$ with $n$ nodes and $dn$ edges. The knapsack capacity is set to be $(e^\top a)/\operatorname{\triangle(G)}$, where $\operatorname{\triangle(G)}$ denotes the maximal degree of the graph $G$. This generation procedure ensures that the capacity constraint is valid. We choose the dimension $n\in\{1000,2000,5000\}$ and $d\in\{2,5\}$. For more results on different choices of $d$, we refer the readers to Appendix \ref{appendix-DQKP}. We do not report the results of RNNAL-Diag because it cannot reach the required accuracy within the 1-hour time limit even for the smallest problem with $n=1000$.

\begin{tiny}
\begin{longtable}[c]{clrrrccccc}
\caption{Computational results for DQKP problems. \label{tab-DQKP}} \\
\toprule
 problem & algorithm & it & itsub & $r$/itA & $\operatorname{R_p}$ & $\operatorname{R_d}$ & $\operatorname{R_c}$ & obj & time   \\
\midrule
\endfirsthead

\multicolumn{10}{c}%
{{ Table \thetable\ continued from previous page}} \\
\toprule
problem & algorithm & it & itsub & $r$/itA & $\operatorname{R_p}$ & $\operatorname{R_d}$ & $\operatorname{R_c}$ & obj & time   \\
\midrule
\endhead
\midrule
\multicolumn{10}{r}{{Continued on next page}} \\
\midrule
\endfoot
\bottomrule
\endlastfoot

$n=1000$ & RNNAL & 20 & 524 & 50 & 6.36e-07 & 1.39e-08 & 1.34e-07 & -2.5519045e+06 & 1.54e+01 \\ 
$d=2$ & SDPNAL+ & 99 & 715 & 5949 & 2.98e-07 & 9.97e-07 & 1.27e-06 & -2.5519060e+06 & 2.69e+03 \\ [3pt] 

$n=1000$ & RNNAL & 35 & 8270 & 116 & 6.80e-07 & 9.44e-07 & 1.53e-10 & -1.2078151e+06 & 2.48e+02 \\ 
$d=5$ & SDPNAL+ & - & - & - & - & - & - & - & - \\ [3pt] 

$n=2000$ & RNNAL & 15 & 1164 & 77 & 7.06e-07 & 3.15e-07 & 4.33e-08 & -9.5394855e+06 & 1.00e+02 \\ 
$d=2$ & SDPNAL+ & - & - & - & - & - & - & - & - \\ [3pt] 

$n=2000$ & RNNAL & 27 & 5059 & 148 & 7.04e-07 & 9.38e-07 & 2.59e-07 & -4.8842259e+06 & 4.61e+02 \\ 
$d=5$ & SDPNAL+ & - & - & - & - & - & - & - & - \\ [3pt] 

$n=5000$ & RNNAL & 10 & 642 & 148 & 9.41e-07 & 4.52e-07 & 3.89e-09 & -5.6073030e+07 & 4.86e+02 \\ 
$d=2$ & SDPNAL+ & - & - & - & - & - & - & - & - \\ [3pt] 

$n=5000$ & RNNAL & 14 & 2727 & 178 & 7.52e-07 & 2.76e-08 & 3.19e-09 & -2.7536080e+07 & 1.81e+03 \\ 
$d=5$ & SDPNAL+ & - & - & - & - & - & - & - & - \\

\end{longtable}
\end{tiny}

As shown in Table \ref{tab-QKP}, RNNAL outperforms SDPNAL+ on all problems. In particular, RNNAL is faster than SDPNAL+ by a factor of 170 times and RNNAL-Diag by a factor of at least 230 times for the first instance. {The time-consuming steps are similar to those encountered in QKP problems.}

\subsection{Gromov-Wasserstein distance}
The Gromov-Wasserstein distance (GWD) aims to solve the following nonconvex QP problem:
\begin{equation}\label{eq-GWD-v1}
\min \left\{-\langle D_X \Pi D_Y,\Pi \rangle:\ \Pi e^k=a,\, \Pi^{\top} e^l=b,\, \Pi\in \R^{l\times k}_+\right\},
\end{equation}
where $e^k$ and $e^l$ are all ones vectors, $D_X\in\S^{l}$ and $D_Y\in\S^{k}$ are two symmetric distance matrices,  $a \in \mathbb{R}_{+}^l$ and $b \in \mathbb{R}_{+}^k$ are two discrete probability distributions satisfying 
$a^\top e^l= 1 = b^\top e^k=1$. Denote $n\coloneqq lk,\ x:=\operatorname{vec}(\Pi)\in\R^{n}$, $Q:=-D_Y \otimes D_X\in \R^{n\times n}$, $d:=(a ; b) \in \mathbb{R}^{l+k}$ and 
$A = (e^k\otimes I_l, I_k\otimes e^l)^\top.$
Then \eqref{eq-GWD-v1} can be equivalently written as \eqref{prob-MBQP} without binary constraints: 
\begin{equation}\label{eq-GWD-v2}
\min\left\{x^{\top} Q x:\ Ax=d,\ x\in \R^{n}_+\right\}.
\end{equation} 
The corresponding DNN relaxation in the form of \eqref{prob-dnn-new} is
\begin{equation}\label{prob-dnn-GWD}
\min\left\{\<Q,X\>:\ Ax=d,AX=dx^\top,\ \begin{pmatrix}1& x^\top \\ x & X\end{pmatrix}\in \mathbb{S}_{+}^{n+1} \cap \mathbb{N}^{n+1}\right\}.
\end{equation}
For \eqref{prob-dnn-GWD}, the feasible set $\M_r$ is an affine space given by
    $\M_r=\left\{ R\in\R^{n\times r}:\ AR=de_1^\top \right\},$
which is smooth if $A$ has full row rank. This can be ensured by removing the last row of $A$ in the prepossessing phase. 

\begin{figure}[H]
    \centering
    \begin{subfigure}[b]{0.13\textwidth}
        \centering
        \includegraphics[width=0.4\textwidth]{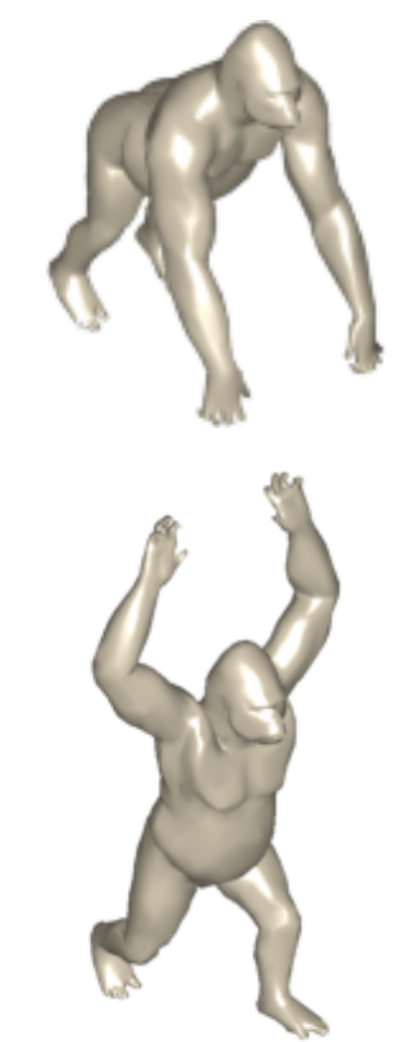}
        \caption{Gorilla}
        \label{fig:Gorilla}
    \end{subfigure}%
    \begin{subfigure}[b]{0.13\textwidth}
    	\centering
        \includegraphics[width=0.4\textwidth]{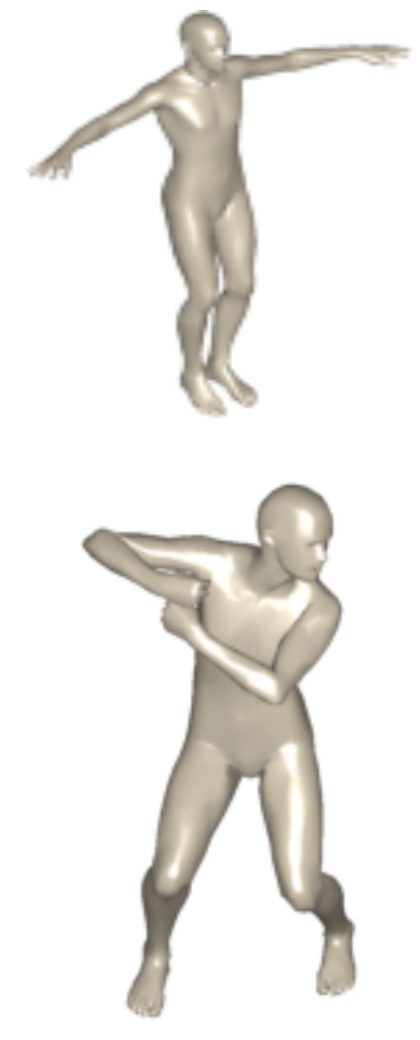}
        \caption{David}
        \label{fig:David}
    \end{subfigure}%
    \begin{subfigure}[b]{0.13\textwidth}
    	\centering
        \includegraphics[width=0.4\textwidth]{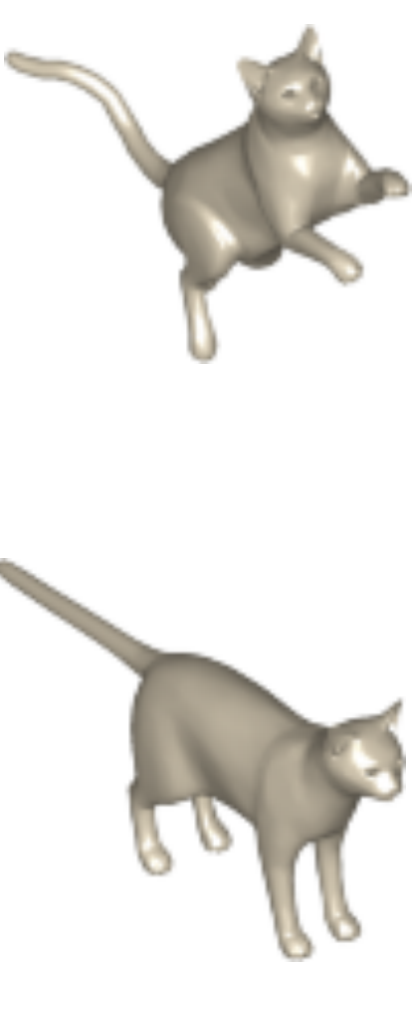}
        \caption{Cat}
        \label{fig:Cat}
    \end{subfigure}%
    \begin{subfigure}[b]{0.13\textwidth}
    	\centering
        \includegraphics[width=0.4\textwidth]{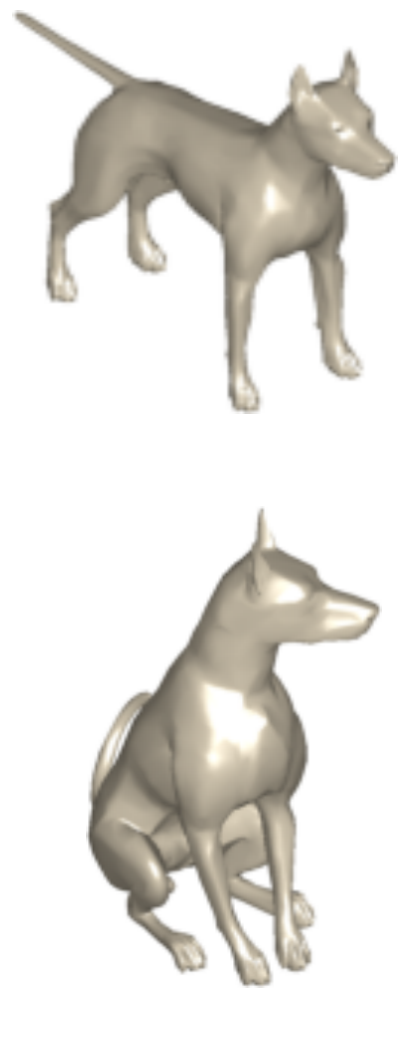}
        \caption{Dog}
        \label{fig:Dog}
    \end{subfigure}%
    \begin{subfigure}[b]{0.13\textwidth}
    	\centering
        \includegraphics[width=0.4\textwidth]{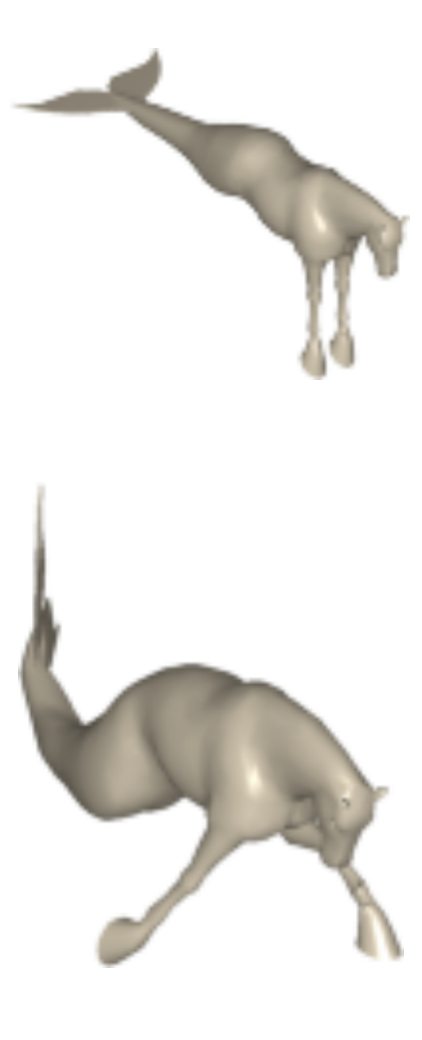}
        \caption{Seahorse}
        \label{fig:Seahorse}
    \end{subfigure}
    \caption{Selected shapes from the TOSCA dataset \cite{bronstein2008numerical}.}
    \label{fig-GWD-animal}
\end{figure}

We use GWD for the shape correspondence task. GWD can also be used for the graph partition task; we refer the readers to Appendix \ref{appendix-GWD-partition} for the computational results on the latter task. Shape correspondence aims to match two similar shapes with the same number of nodes. We use shapes from the TOSCA dataset \cite{bronstein2008numerical}, which includes eight classes of 3D shapes in various poses. First, we select five classes and randomly choose two shapes with different poses within each class. See Figure \ref{fig-GWD-animal} for the selected shapes. Then, we sample $l$ nodes uniformly for each shape and set the $(i,j)$-th element of the distance matrix as the shortest path between the $i$-th and $j$-th nodes. $a$ and $b$ follow discrete uniform distributions. The number of nodes $k=l$ is chosen from $\{35, 45\}$. For more results with different $l$ values, see Appendix \ref{appendix-GWD}. We do not include the results of RNNAL-Diag as it fails to achieve the required accuracy within the 1-hour time limit, even for the smallest problem with $l=30$.

\begin{tiny}
\begin{longtable}[c]{clrrrccccc}
\caption{Computational results for GWD shape correspondence problems. \label{tab-GWD-shape}} \\
\toprule
 problem & algorithm & it & itsub & $r$/itA & $\operatorname{R_p}$ & $\operatorname{R_d}$ & $\operatorname{R_c}$ & obj & time   \\
\midrule
\endfirsthead

\multicolumn{10}{c}%
{{ Table \thetable\ continued from previous page}} \\
\toprule
problem & algorithm & it & itsub & $r$/itA & $\operatorname{R_p}$ & $\operatorname{R_d}$ & $\operatorname{R_c}$ & obj & time   \\
\midrule
\endhead
\midrule
\multicolumn{10}{r}{{Continued on next page}} \\
\midrule
\endfoot
\bottomrule
\endlastfoot

Cat& RNNAL & 42 & 10862 & 42 & 4.66e-07 & 3.09e-07 & 6.89e-07 & 1.3986652e+05 & 2.26e+02 \\ 
$n=1225$ & SDPNAL+ & 572 & 1002 & 10988 & 5.30e-07 & 1.85e-07 & 8.72e-15 & 1.3990409e+05 & 3.54e+03 \\ [3pt] 

David& RNNAL & 76 & 12694 & 14 & 2.72e-07 & 9.36e-07 & 1.77e-07 & 2.8688769e+05 & 2.32e+02 \\ 
$n=1225$ & SDPNAL+ & - & - & - & - & - & - & - & - \\ [3pt] 

Dog& RNNAL & 41 & 5503 & 7 & 4.72e-07 & 3.05e-08 & 4.75e-07 & 9.4759701e+04 & 1.38e+02 \\ 
$n=1225$ & SDPNAL+ & 491 & 819 & 11122 & 2.32e-07 & 3.90e-07 & 3.35e-15 & 9.4851634e+04 & 3.18e+03 \\ [3pt] 

Gorilla& RNNAL & 71 & 13766 & 13 & 7.83e-07 & 5.85e-07 & 3.86e-07 & 2.3616473e+05 & 2.56e+02 \\ 
$n=1225$ & SDPNAL+ & 376 & 585 & 11275 & 5.80e-07 & 4.78e-07 & 3.21e-07 & 2.3615272e+05 & 2.96e+03 \\ [3pt] 

Seahorse& RNNAL & 74 & 13887 & 12 & 9.17e-07 & 2.89e-07 & 7.48e-08 & 8.8129167e+05 & 2.64e+02 \\ 
$n=1225$ & SDPNAL+ & - & - & - & - & - & - & - & - \\ [3pt] 

Cat& RNNAL & 108 & 21820 & 15 & 6.26e-07 & 1.39e-07 & 9.49e-07 & 3.4775634e+05 & 1.33e+03 \\ 
$n=2025$ & SDPNAL+ & - & - & - & - & - & - & - & - \\ [3pt] 

David& RNNAL & 64 & 14851 & 37 & 5.97e-07 & 9.46e-07 & 7.85e-07 & 2.4304548e+05 & 9.35e+02 \\ 
$n=2025$ & SDPNAL+ & - & - & - & - & - & - & - & - \\ [3pt] 

Dog& RNNAL & 73 & 16132 & 17 & 9.41e-07 & 8.11e-07 & 2.19e-07 & 1.1140722e+05 & 9.21e+02 \\ 
$n=2025$ & SDPNAL+ & - & - & - & - & - & - & - & - \\ [3pt] 

Gorilla& RNNAL & 105 & 18979 & 6 & 7.23e-07 & 2.54e-07 & 8.01e-07 & 3.6744730e+05 & 1.10e+03 \\ 
$n=2025$ & SDPNAL+ & - & - & - & - & - & - & - & - \\ [3pt] 

Seahorse& RNNAL & 84 & 17938 & 35 & 7.56e-07 & 7.92e-07 & 9.41e-07 & 6.8249936e+05 & 1.11e+03 \\ 
$n=2025$ & SDPNAL+ & - & - & - & - & - & - & - & - \\ 
\end{longtable}
\end{tiny}

\begin{figure}[ht]
    \centering
    \begin{subfigure}{0.32\textwidth}
        \centering
        \includegraphics[width=\linewidth]{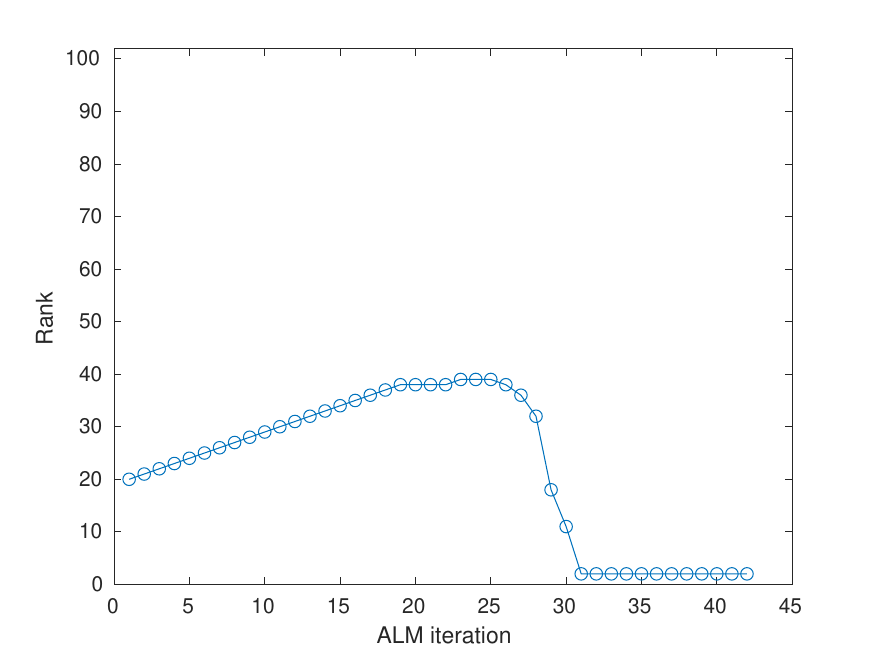}
        \caption{$r_0=n/5$}
        \label{fig:sub1}
    \end{subfigure}
    \begin{subfigure}{0.32\textwidth}
        \centering
        \includegraphics[width=\linewidth]{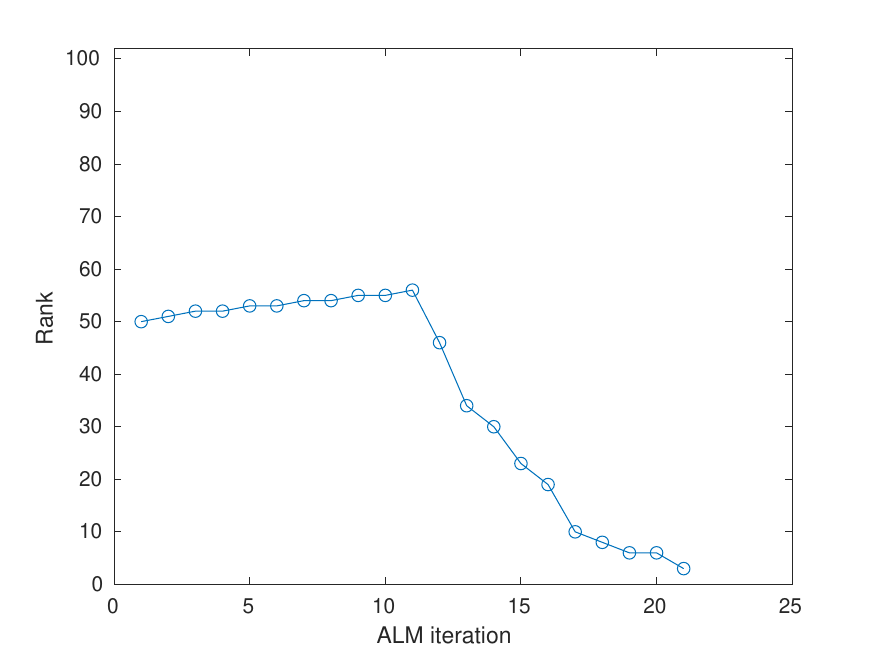}
        \caption{$r_0=n/2$}
        \label{fig:sub2}
    \end{subfigure}
    \begin{subfigure}{0.32\textwidth}
        \centering
        \includegraphics[width=\linewidth]{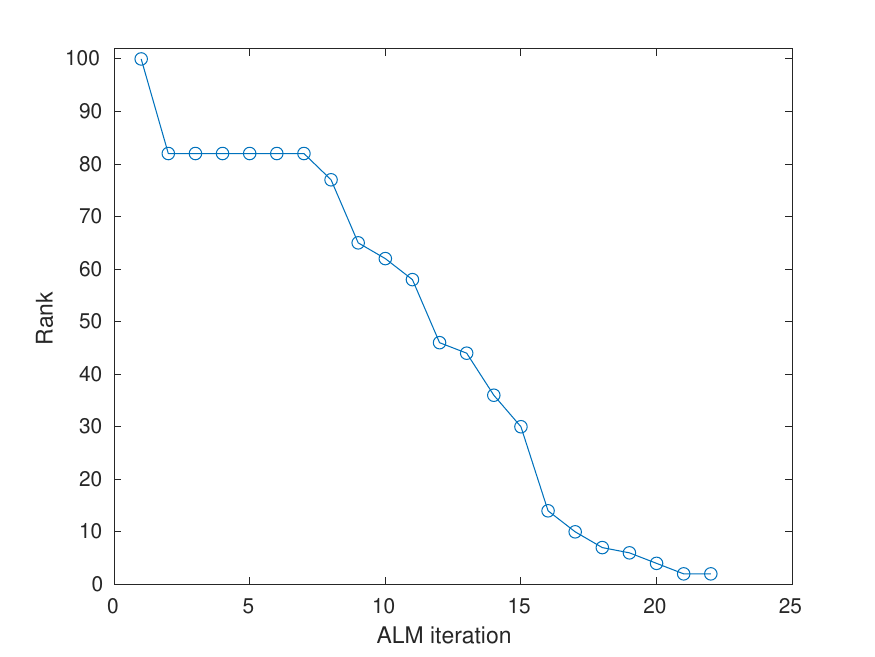}
        \caption{$r_0=n$}
        \label{fig:sub3}
    \end{subfigure}
    \caption{The rank evolvement under different initial rank $r_0$.}
    \label{fig-rank}
\end{figure}

As shown in Table \ref{tab-GWD-shape}, RNNAL achieves the required accuracy within the 1 hour limit, while SDPNAL+ fails for some problems with dimensions $n=1225$ and RNNAL-Diag fails for all problems. RNNAL outperforms SDPNAL+ and RNNAL-Diag on all problems. In particular, RNNAL is faster than SDPNAL+ and RNNAL-Diag by at least a factor of 25 times for the third instance. {Note that computing the projection and retraction onto the affine space $\M_r$ is both straightforward and computationally efficient, accounting for less than $10\%$ of the total computation time. To demonstrate the effectiveness of our rank-adaptive procedure, we apply RNNAL to solve a GWD problem with $n = 100$, which is generated by sampling 10 points from the ``Cat" dataset in TOSCA. We select the initial rank $r_0 \in \{n/5, n/2, n\}$. The evolution of the rank of the iterates during the procedure is depicted in Figure \ref{fig-rank}.
As shown in Figure~\ref{fig-rank}, regardless of the initial rank, our rank-adaptive procedure automatically adjusts the rank to achieve convergence.
}

\section{Conclusion}\label{sec-conclusion}
In this paper, we propose an augmented Lagrangian method utilizing low-rank factorization to solve DNN relaxations of large-scale mixed-binary quadratic programs with guaranteed convergence to global optimal solutions. We avoid the non-smoothness of the feasible set after BM factorization, which arises from the violation of Slater's condition, by making the key observation that we can 
reformulate most of the quadratic constraints into fewer, more manageable affine constraints, as well as through applying a new constraint-relaxation strategy. We offer theoretical analysis and practical methods to accelerate the computation of the projection and retraction for a class of algebraic varieties $\M_r^g$. Our numerical experiments on solving different classes of large-scale DNN instances have demonstrated the excellent practical performance of the proposed method when compared to other state-of-the-art solvers. Our work provides a prototype to solve general SDPs with additional polyhedral constraints \eqref{prob-P-dnn-YZ} including DNN problems.


\appendix

\section{Appendix}

\subsection{Experiments on QAP problems}
\label{appendix-QAP}

\begin{tiny}
\begin{longtable}[c]{clrrrccccc}
\caption{Computational results for QAP problems. \label{tab-QAP}} \\
\toprule
 problem & algorithm & it & itsub & $r$/itA & $\operatorname{R_p}$ & $\operatorname{R_d}$ & $\operatorname{R_c}$ & obj & time   \\
\midrule
\endfirsthead

\multicolumn{10}{c}%
{{ Table \thetable\ continued from previous page}} \\
\toprule
problem & algorithm & it & itsub & $r$/itA & $\operatorname{R_p}$ & $\operatorname{R_d}$ & $\operatorname{R_c}$ & obj & time   \\
\midrule
\endhead
\midrule
\multicolumn{10}{r}{{Continued on next page}} \\
\midrule
\endfoot
\bottomrule
\endlastfoot

chr22a& RNNAL & 20 & 3640 & 272 & 3.57e-08 & 3.67e-10 & 5.32e-09 & 6.1560007e+03 & 1.96e+02 \\ 
$n=484$ & SDPNAL+ & 82 & 151 & 3100 & 1.55e-07 & 1.49e-08 & 1.18e-14 & 6.1560002e+03 & 1.64e+02 \\ 
& RNNAL-Diag & 11 & 45775 & 251 & 9.61e-07 & 2.90e-10 & 9.63e-19 & 6.1560188e+03 & 3.86e+02 \\ [5pt] 

chr22b& RNNAL & 28 & 1656 & 151 & 9.88e-07 & 5.01e-09 & 1.33e-08 & 6.1940017e+03 & 1.03e+02 \\ 
$n=484$ & SDPNAL+ & 90 & 164 & 3795 & 8.29e-07 & 9.10e-07 & 6.49e-07 & 6.1940560e+03 & 2.13e+02 \\ 
& RNNAL-Diag & 14 & 77461 & 254 & 9.87e-07 & 2.09e-09 & 7.18e-18 & 6.1940192e+03 & 6.62e+02 \\ [5pt] 

chr25a& RNNAL & 27 & 2518 & 169 & 8.23e-07 & 8.73e-10 & 3.88e-11 & 3.7960037e+03 & 1.80e+02 \\ 
$n=625$ & SDPNAL+ & 62 & 113 & 2450 & 8.90e-08 & 2.69e-07 & 1.59e-15 & 3.7959993e+03 & 2.00e+02 \\ 
& RNNAL-Diag & 13 & 75720 & 325 & 1.00e-06 & 6.15e-10 & 2.54e-17 & 3.7960159e+03 & 1.03e+03 \\ [5pt] 

esc32a& RNNAL & 92 & 16028 & 832 & 8.85e-07 & 9.89e-07 & 5.28e-09 & 1.0331998e+02 & 1.63e+03 \\ 
$n=1024$ & SDPNAL+ & 70 & 70 & 1341 & 7.46e-07 & 9.98e-07 & 4.84e-08 & 1.0332040e+02 & 3.09e+02 \\ 
& RNNAL-Diag & - & - & - & - & - & - & - & - \\ [5pt] 

esc32b& RNNAL & 69 & 6218 & 834 & 3.84e-07 & 1.00e-06 & 1.21e-09 & 1.3188403e+02 & 6.55e+02 \\ 
$n=1024$ & SDPNAL+ & 22 & 22 & 547 & 9.65e-07 & 8.56e-07 & 8.98e-08 & 1.3188506e+02 & 1.18e+02 \\ 
& RNNAL-Diag & - & - & - & - & - & - & - & - \\ [5pt] 

esc32c& RNNAL & 30 & 1084 & 592 & 5.07e-07 & 1.19e-07 & 1.50e-09 & 6.1518097e+02 & 9.44e+01 \\ 
$n=1024$ & SDPNAL+ & 441 & 597 & 6769 & 9.33e-07 & 8.32e-07 & 4.97e-07 & 6.1517700e+02 & 1.95e+03 \\ 
& RNNAL-Diag & - & - & - & - & - & - & - & - \\ [5pt] 

esc32d& RNNAL & 49 & 1819 & 910 & 3.77e-07 & 9.60e-07 & 1.39e-09 & 1.9022708e+02 & 2.42e+02 \\ 
$n=1024$ & SDPNAL+ & 5 & 5 & 271 & 3.06e-07 & 9.29e-07 & 3.97e-15 & 1.9022755e+02 & 5.44e+01 \\ 
& RNNAL-Diag & - & - & - & - & - & - & - & - \\ [5pt] 

esc32e& RNNAL & 35 & 601 & 547 & 5.56e-08 & 5.09e-07 & 5.18e-11 & 1.9000341e+00 & 6.98e+01 \\ 
$n=1024$ & SDPNAL+ & 2 & 2 & 313 & 5.45e-08 & 1.38e-07 & 6.68e-16 & 1.8999252e+00 & 6.11e+01 \\ 
& RNNAL-Diag & - & - & - & - & - & - & - & - \\ [5pt] 

esc32g& RNNAL & 33 & 772 & 562 & 8.00e-07 & 5.04e-07 & 3.18e-11 & 5.8332343e+00 & 7.86e+01 \\ 
$n=1024$ & SDPNAL+ & 2 & 2 & 251 & 4.34e-08 & 3.36e-07 & 1.18e-14 & 5.8333664e+00 & 4.87e+01 \\ 
& RNNAL-Diag & - & - & - & - & - & - & - & - \\ [5pt] 

esc32h& RNNAL & 99 & 11042 & 916 & 1.53e-07 & 9.61e-07 & 8.90e-10 & 4.2440276e+02 & 1.28e+03 \\ 
$n=1024$ & SDPNAL+ & 34 & 34 & 747 & 5.58e-07 & 9.99e-07 & 8.38e-09 & 4.2440191e+02 & 1.63e+02 \\ 
& RNNAL-Diag & - & - & - & - & - & - & - & - \\ [5pt] 

kra30a& RNNAL & 88 & 14834 & 529 & 7.21e-07 & 9.36e-07 & 8.69e-10 & 8.6837278e+04 & 9.74e+02 \\ 
$n=900$ & SDPNAL+ & 33 & 61 & 1054 & 7.09e-07 & 9.93e-07 & 1.01e-08 & 8.6837314e+04 & 2.52e+02 \\ 
& RNNAL-Diag & - & - & - & - & - & - & - & - \\ [5pt] 

kra30b& RNNAL & 90 & 14296 & 543 & 4.55e-07 & 9.50e-07 & 2.15e-10 & 8.7857917e+04 & 9.44e+02 \\ 
$n=900$ & SDPNAL+ & 39 & 80 & 1024 & 7.71e-07 & 8.74e-07 & 6.65e-08 & 8.7857949e+04 & 3.47e+02 \\ 
& RNNAL-Diag & - & - & - & - & - & - & - & - \\ [5pt] 

kra32& RNNAL & 87 & 15597 & 621 & 8.12e-07 & 9.31e-07 & 1.14e-09 & 8.5774986e+04 & 1.28e+03 \\ 
$n=1024$ & SDPNAL+ & 31 & 39 & 781 & 9.68e-07 & 5.94e-07 & 2.20e-07 & 8.5775114e+04 & 2.55e+02 \\ 
& RNNAL-Diag & - & - & - & - & - & - & - & - \\ [5pt] 

lipa30a& RNNAL & 30 & 2304 & 282 & 2.53e-08 & 2.87e-09 & 7.10e-08 & 1.3178000e+04 & 1.91e+02 \\ 
$n=900$ & SDPNAL+ & 44 & 75 & 1276 & 5.77e-07 & 4.93e-08 & 2.82e-14 & 1.3178004e+04 & 2.56e+02 \\ 
& RNNAL-Diag & 8 & 35753 & 453 & 8.08e-07 & 1.70e-10 & 2.20e-16 & 1.3178004e+04 & 1.01e+03 \\ [5pt] 

lipa30b& RNNAL & 17 & 665 & 497 & 1.25e-08 & 2.64e-10 & 7.34e-09 & 1.5142600e+05 & 7.89e+01 \\ 
$n=900$ & SDPNAL+ & 26 & 135 & 800 & 1.00e-08 & 1.17e-07 & 2.28e-15 & 1.5142600e+05 & 3.15e+02 \\ 
& RNNAL-Diag & 3 & 4267 & 443 & 3.89e-07 & 1.34e-11 & 5.88e-17 & 1.5142599e+05 & 1.24e+02 \\ [5pt] 

nug21& RNNAL & 57 & 25733 & 283 & 8.70e-07 & 9.90e-07 & 9.98e-10 & 2.3819302e+03 & 4.57e+02 \\ 
$n=441$ & SDPNAL+ & 43 & 103 & 1608 & 1.02e-06 & 8.70e-07 & 1.22e-07 & 2.3819331e+03 & 1.04e+02 \\ 
& RNNAL-Diag & - & - & - & - & - & - & - & - \\ [5pt] 

nug22& RNNAL & 80 & 26509 & 283 & 7.83e-07 & 9.44e-07 & 1.42e-09 & 3.5286803e+03 & 5.19e+02 \\ 
$n=484$ & SDPNAL+ & 51 & 133 & 1753 & 9.89e-07 & 9.80e-07 & 3.49e-08 & 3.5286821e+03 & 1.67e+02 \\ 
& RNNAL-Diag & - & - & - & - & - & - & - & - \\ [5pt] 

nug24& RNNAL & 59 & 11317 & 315 & 6.53e-07 & 7.27e-07 & 3.39e-10 & 3.4010722e+03 & 3.01e+02 \\ 
$n=576$ & SDPNAL+ & 43 & 101 & 1227 & 1.82e-06 & 9.13e-07 & 7.91e-08 & 3.4010705e+03 & 1.80e+02 \\ 
& RNNAL-Diag & - & - & - & - & - & - & - & - \\ [5pt] 

nug25& RNNAL & 66 & 17045 & 363 & 7.44e-07 & 8.34e-07 & 1.83e-09 & 3.6258728e+03 & 5.31e+02 \\ 
$n=625$ & SDPNAL+ & 42 & 105 & 1591 & 1.49e-06 & 8.63e-07 & 3.41e-07 & 3.6258759e+03 & 2.31e+02 \\ 
& RNNAL-Diag & - & - & - & - & - & - & - & - \\ [5pt] 

nug27& RNNAL & 79 & 25466 & 426 & 8.11e-07 & 9.80e-07 & 1.11e-09 & 5.1296119e+03 & 1.03e+03 \\ 
$n=729$ & SDPNAL+ & 42 & 106 & 1804 & 1.02e-06 & 9.34e-07 & 9.24e-08 & 5.1296145e+03 & 2.90e+02 \\ 
& RNNAL-Diag & - & - & - & - & - & - & - & - \\ [5pt] 

nug28& RNNAL & 58 & 14747 & 471 & 8.94e-07 & 9.30e-07 & 3.09e-09 & 5.0256690e+03 & 6.86e+02 \\ 
$n=784$ & SDPNAL+ & 62 & 159 & 1469 & 1.12e-06 & 8.92e-07 & 1.49e-07 & 5.0256725e+03 & 4.36e+02 \\ 
& RNNAL-Diag & - & - & - & - & - & - & - & - \\ [5pt] 

nug30& RNNAL & 59 & 15565 & 532 & 9.94e-07 & 8.39e-07 & 3.86e-09 & 5.9494594e+03 & 9.66e+02 \\ 
$n=900$ & SDPNAL+ & 43 & 134 & 1144 & 9.86e-07 & 9.28e-07 & 2.54e-07 & 5.9494661e+03 & 6.07e+02 \\ 
& RNNAL-Diag & - & - & - & - & - & - & - & - \\ [5pt] 

tai25a& RNNAL & 956 & 11319 & 430 & 5.48e-07 & 9.34e-07 & 2.61e-10 & 1.1131360e+06 & 1.08e+03 \\ 
$n=625$ & SDPNAL+ & 131 & 131 & 2301 & 6.48e-07 & 1.00e-06 & 4.60e-08 & 1.1130116e+06 & 1.78e+02 \\ 
& RNNAL-Diag & - & - & - & - & - & - & - & - \\ [5pt] 

tai25b& RNNAL & 257 & 102893 & 464 & 9.14e-07 & 8.33e-07 & 2.53e-10 & 3.3802176e+08 & 3.50e+03 \\ 
$n=625$ & SDPNAL+ & 360 & 1026 & 19840 & 1.28e-06 & 1.19e-06 & 1.19e-07 & 3.3801861e+08 & 2.75e+03 \\ 
& RNNAL-Diag & - & - & - & - & - & - & - & - \\ [5pt] 

tai30a& RNNAL & 34 & 1402 & 524 & 7.29e-07 & 9.14e-07 & 4.34e-09 & 1.7068712e+06 & 1.15e+02 \\ 
$n=900$ & SDPNAL+ & 21 & 24 & 499 & 9.43e-07 & 5.93e-07 & 1.30e-15 & 1.7068712e+06 & 1.26e+02 \\ 
& RNNAL-Diag & - & - & - & - & - & - & - & - \\ [5pt] 

tai35a& RNNAL & 33 & 1235 & 703 & 6.44e-07 & 8.05e-07 & 3.84e-09 & 2.2166461e+06 & 1.98e+02 \\ 
$n=1225$ & SDPNAL+ & 10 & 10 & 430 & 9.36e-07 & 9.86e-07 & 4.80e-08 & 2.2166458e+06 & 1.51e+02 \\ 
& RNNAL-Diag & - & - & - & - & - & - & - & - \\ [5pt] 

tho30& RNNAL & 63 & 20510 & 534 & 6.88e-07 & 9.94e-07 & 6.99e-10 & 1.4357596e+05 & 1.22e+03 \\ 
$n=900$ & SDPNAL+ & 41 & 108 & 1454 & 1.12e-06 & 9.67e-07 & 9.58e-08 & 1.4357613e+05 & 5.34e+02 \\ 
& RNNAL-Diag & - & - & - & - & - & - & - & - \\

\end{longtable}
\end{tiny}

\subsection{Experiments on BIQ problems}
\label{appendix-BIQ}

\begin{tiny}
\begin{longtable}[c]{clrrrccccc}
\caption{Computational results for BIQ problems. \label{tab-appendix-BIQ}} \\
\toprule
problem  & algorithm & it & itsub & $r$/itA & $\operatorname{R_p}$ & $\operatorname{R_d}$ & $\operatorname{R_c}$ & obj & time   \\
\midrule
\endfirsthead

\multicolumn{10}{c}%
{{ Table \thetable\ continued from previous page}} \\
\toprule
problem  & algorithm & it & itsub & $r$/itA & $\operatorname{R_p}$ & $\operatorname{R_d}$ & $\operatorname{R_c}$ & obj & time   \\
\midrule
\endhead
\midrule
\multicolumn{10}{r}{{Continued on next page}} \\
\midrule
\endfoot

\bottomrule
\endlastfoot

bqp1000.1& RNNAL & 11 & 911 & 55 & 9.04e-07 & 3.17e-07 & 3.80e-07 & -3.9849472e+05 & 9.40e+00 \\ 
$n= 1000$& SDPNAL+ & 118 & 172 & 2752 & 8.92e-07 & 9.69e-07 & 5.54e-07 & -3.9849494e+05 & 2.61e+02 \\ [5pt] 

bqp1000.2& RNNAL & 11 & 1047 & 55 & 8.29e-07 & 5.32e-07 & 8.11e-08 & -3.8430752e+05 & 1.06e+01 \\ 
$n= 1000$& SDPNAL+ & 117 & 178 & 2750 & 4.34e-07 & 4.21e-07 & 5.55e-07 & -3.8430730e+05 & 2.71e+02 \\ [5pt] 

bqp1000.3& RNNAL & 12 & 942 & 59 & 9.04e-07 & 9.05e-07 & 6.61e-08 & -3.9883820e+05 & 9.92e+00 \\ 
$n= 1000$& SDPNAL+ & 127 & 193 & 2931 & 5.96e-08 & 9.66e-07 & 2.41e-07 & -3.9883798e+05 & 2.85e+02 \\ [5pt] 

bqp1000.4& RNNAL & 11 & 864 & 53 & 6.59e-07 & 2.00e-07 & 4.51e-08 & -3.9868711e+05 & 9.21e+00 \\ 
$n= 1000$& SDPNAL+ & 115 & 169 & 2600 & 3.63e-07 & 4.76e-07 & 4.34e-07 & -3.9868689e+05 & 2.52e+02 \\ [5pt] 

bqp1000.5& RNNAL & 16 & 1183 & 51 & 5.48e-07 & 7.18e-07 & 7.25e-09 & -3.8297578e+05 & 1.28e+01 \\ 
$n= 1000$& SDPNAL+ & 117 & 165 & 2636 & 5.41e-08 & 9.22e-07 & 3.37e-07 & -3.8297575e+05 & 2.52e+02 \\ [5pt] 

bqp1000.6& RNNAL & 12 & 717 & 58 & 7.66e-07 & 8.32e-07 & 8.22e-11 & -3.8617571e+05 & 8.12e+00 \\ 
$n= 1000$& SDPNAL+ & 117 & 165 & 2772 & 5.28e-09 & 9.99e-07 & 5.20e-07 & -3.8617568e+05 & 2.63e+02 \\ [5pt] 

bqp1000.7& RNNAL & 12 & 866 & 57 & 9.01e-07 & 5.16e-07 & 3.69e-08 & -3.9951364e+05 & 9.16e+00 \\ 
$n= 1000$& SDPNAL+ & 117 & 179 & 2885 & 1.24e-08 & 9.98e-07 & 9.23e-07 & -3.9951374e+05 & 2.72e+02 \\ [5pt] 

bqp1000.8& RNNAL & 12 & 820 & 56 & 8.06e-07 & 5.98e-07 & 4.91e-09 & -3.8355968e+05 & 8.62e+00 \\ 
$n= 1000$& SDPNAL+ & 126 & 195 & 2928 & 4.96e-08 & 9.90e-07 & 5.28e-07 & -3.8355966e+05 & 2.81e+02 \\ [5pt] 

bqp1000.9& RNNAL & 12 & 770 & 57 & 7.29e-07 & 4.38e-07 & 5.50e-07 & -3.7902773e+05 & 8.34e+00 \\ 
$n= 1000$& SDPNAL+ & 128 & 182 & 2932 & 3.03e-08 & 9.70e-07 & 3.18e-07 & -3.7902777e+05 & 2.78e+02 \\ [5pt] 

bqp1000.10& RNNAL & 11 & 725 & 53 & 7.77e-07 & 2.95e-07 & 3.82e-08 & -3.7949962e+05 & 7.77e+00 \\ 
$n= 1000$& SDPNAL+ & 117 & 159 & 2652 & 3.91e-08 & 9.22e-07 & 2.65e-07 & -3.7949932e+05 & 2.51e+02 \\ [5pt] 

bqp2500.1& RNNAL & 9 & 817 & 97 & 6.11e-07 & 4.53e-07 & 5.57e-08 & -1.6354913e+06 & 1.03e+02 \\ 
$n= 2500$& SDPNAL+ & - & - & - & - & - & - & - & - \\ [5pt] 

bqp2500.2& RNNAL & 7 & 703 & 109 & 9.71e-07 & 8.28e-07 & 4.84e-09 & -1.5975405e+06 & 8.66e+01 \\ 
$n= 2500$& SDPNAL+ & - & - & - & - & - & - & - & - \\ [5pt] 

bqp2500.3& RNNAL & 7 & 674 & 109 & 6.73e-07 & 5.73e-07 & 1.21e-07 & -1.5392082e+06 & 8.35e+01 \\ 
$n= 2500$& SDPNAL+ & - & - & - & - & - & - & - & - \\ [5pt] 

bqp2500.4& RNNAL & 7 & 757 & 104 & 9.85e-07 & 3.14e-07 & 1.02e-07 & -1.6247929e+06 & 9.19e+01 \\ 
$n= 2500$& SDPNAL+ & - & - & - & - & - & - & - & - \\ [5pt] 

bqp2500.5& RNNAL & 8 & 731 & 95 & 7.23e-07 & 2.60e-07 & 4.39e-08 & -1.6089563e+06 & 9.11e+01 \\ 
$n= 2500$& SDPNAL+ & - & - & - & - & - & - & - & - \\ [5pt] 

bqp2500.6& RNNAL & 8 & 937 & 99 & 7.79e-07 & 2.41e-07 & 7.12e-08 & -1.5912762e+06 & 1.15e+02 \\ 
$n= 2500$& SDPNAL+ & - & - & - & - & - & - & - & - \\ [5pt] 

bqp2500.7& RNNAL & 8 & 713 & 103 & 9.24e-07 & 8.08e-07 & 2.60e-07 & -1.6018316e+06 & 9.08e+01 \\ 
$n= 2500$& SDPNAL+ & - & - & - & - & - & - & - & - \\ [5pt] 

bqp2500.8& RNNAL & 7 & 1075 & 105 & 8.88e-07 & 7.39e-07 & 1.11e-07 & -1.5981038e+06 & 1.27e+02 \\ 
$n= 2500$& SDPNAL+ & - & - & - & - & - & - & - & - \\ [5pt] 

bqp2500.9& RNNAL & 8 & 745 & 111 & 9.26e-07 & 7.76e-07 & 1.63e-07 & -1.6041088e+06 & 9.43e+01 \\ 
$n= 2500$& SDPNAL+ & - & - & - & - & - & - & - & - \\ [5pt] 

bqp2500.10& RNNAL & 8 & 897 & 102 & 7.30e-07 & 7.27e-07 & 9.34e-08 & -1.6081578e+06 & 1.10e+02 \\ 
$n= 2500$& SDPNAL+ & - & - & - & - & - & - & - & - \\ 

\end{longtable}
\end{tiny}

\subsection{Experiments on $\theta_+$ problems}
\label{appendix-theta}

\begin{tiny}
\begin{longtable}[c]{clrrrccccc}
\caption{Computational results for $\theta_+$ problems. \label{theta}} \\
\toprule
problem  & algorithm & it & itsub & $r$/itA & $\operatorname{R_p}$ & $\operatorname{R_d}$ & $\operatorname{R_c}$ & obj & time   \\
\midrule
\endfirsthead

\multicolumn{10}{c}%
{{ Table \thetable\ continued from previous page}} \\
\toprule
problem  & algorithm & it & itsub & $r$/itA & $\operatorname{R_p}$ & $\operatorname{R_d}$ & $\operatorname{R_c}$ & obj & time   \\
\midrule
\endhead
\midrule
\multicolumn{10}{r}{{Continued on next page}} \\
\midrule
\endfoot

\bottomrule
\endlastfoot

G1& RNNAL & 17 & 1197 & 119 & 2.95e-07 & 5.66e-08 & 3.14e-07 & -1.4424460e+02 & 1.16e+01 \\ 
$n=  800$& SDPNAL+ & 36 & 46 & 1100 & 5.01e-07 & 1.65e-07 & 3.44e-14 & -1.4424460e+02 & 1.07e+02 \\ [3pt] 

G2& RNNAL & 15 & 938 & 120 & 6.49e-07 & 3.60e-08 & 3.71e-07 & -1.4456426e+02 & 9.71e+00 \\ 
$n=  800$& SDPNAL+ & 25 & 26 & 870 & 1.21e-07 & 9.62e-07 & 9.49e-15 & -1.4456409e+02 & 7.90e+01 \\ [3pt] 

G3& RNNAL & 15 & 844 & 119 & 3.57e-07 & 7.42e-08 & 4.46e-07 & -1.4447616e+02 & 9.46e+00 \\ 
$n=  800$& SDPNAL+ & 35 & 35 & 964 & 8.06e-07 & 7.30e-07 & 8.70e-07 & -1.4447651e+02 & 8.76e+01 \\ [3pt] 

G4& RNNAL & 14 & 1118 & 117 & 1.58e-07 & 7.37e-08 & 7.75e-09 & -1.4457533e+02 & 1.18e+01 \\ 
$n=  800$& SDPNAL+ & 38 & 40 & 1100 & 4.88e-07 & 2.41e-07 & 3.74e-14 & -1.4457530e+02 & 1.06e+02 \\ [3pt] 

G5& RNNAL & 18 & 961 & 119 & 3.09e-07 & 1.68e-07 & 3.10e-08 & -1.4449466e+02 & 1.05e+01 \\ 
$n=  800$& SDPNAL+ & 35 & 35 & 949 & 6.96e-07 & 8.89e-07 & 3.25e-07 & -1.4449488e+02 & 8.56e+01 \\ [3pt] 

G6& RNNAL & 17 & 1197 & 119 & 2.95e-07 & 5.66e-08 & 3.14e-07 & -1.4424460e+02 & 1.16e+01 \\ 
$n=  800$& SDPNAL+ & 36 & 46 & 1100 & 5.01e-07 & 1.65e-07 & 3.44e-14 & -1.4424460e+02 & 1.07e+02 \\ [3pt] 

G7& RNNAL & 15 & 938 & 120 & 6.49e-07 & 3.60e-08 & 3.71e-07 & -1.4456426e+02 & 9.82e+00 \\ 
$n=  800$& SDPNAL+ & 25 & 26 & 870 & 1.21e-07 & 9.62e-07 & 9.49e-15 & -1.4456409e+02 & 7.88e+01 \\ [3pt] 

G8& RNNAL & 15 & 844 & 119 & 3.57e-07 & 7.42e-08 & 4.46e-07 & -1.4447616e+02 & 9.30e+00 \\ 
$n=  800$& SDPNAL+ & 35 & 35 & 964 & 8.06e-07 & 7.30e-07 & 8.70e-07 & -1.4447651e+02 & 8.69e+01 \\ [3pt] 

G9& RNNAL & 14 & 1118 & 117 & 1.58e-07 & 7.37e-08 & 7.75e-09 & -1.4457533e+02 & 1.21e+01 \\ 
$n=  800$& SDPNAL+ & 38 & 40 & 1100 & 4.88e-07 & 2.41e-07 & 3.74e-14 & -1.4457530e+02 & 1.06e+02 \\ [3pt] 

G10& RNNAL & 18 & 961 & 119 & 3.09e-07 & 1.68e-07 & 3.10e-08 & -1.4449466e+02 & 1.05e+01 \\ 
$n=  800$& SDPNAL+ & 35 & 35 & 949 & 6.96e-07 & 8.89e-07 & 3.25e-07 & -1.4449488e+02 & 8.56e+01 \\ [3pt] 

G11& RNNAL & 14 & 5645 & 4 & 5.49e-07 & 4.07e-09 & 4.63e-07 & -3.9999913e+02 & 2.39e+01 \\ 
$n=  800$& SDPNAL+ & 801 & 880 & 19399 & 1.88e-05 & 5.45e-05 & 1.14e-12 & -3.9981544e+02 & 8.86e+02 \\ [3pt] 

G12& RNNAL & 9 & 3529 & 12 & 3.66e-07 & 7.80e-09 & 2.99e-07 & -3.9999982e+02 & 1.69e+01 \\ 
$n=  800$& SDPNAL+ & 70 & 203 & 2650 & 2.08e-07 & 4.91e-07 & 3.34e-12 & -3.9999995e+02 & 1.66e+02 \\ [3pt] 

G13& RNNAL & 17 & 1583 & 19 & 7.84e-07 & 8.42e-07 & 8.82e-10 & -3.9841674e+02 & 9.71e+00 \\ 
$n=  800$& SDPNAL+ & 63 & 196 & 2339 & 7.38e-13 & 6.49e-07 & 1.99e-06 & -3.9841542e+02 & 1.88e+02 \\ [3pt] 

G14& RNNAL & 15 & 2855 & 148 & 9.58e-07 & 2.38e-07 & 5.11e-10 & -2.7899999e+02 & 2.43e+01 \\ 
$n=  800$& SDPNAL+ & 131 & 372 & 6700 & 1.02e-06 & 6.92e-07 & 6.80e-08 & -2.7900027e+02 & 6.81e+02 \\ [3pt] 

G15& RNNAL & 27 & 31842 & 148 & 6.70e-07 & 6.39e-07 & 2.70e-11 & -2.8374869e+02 & 2.34e+02 \\ 
$n=  800$& SDPNAL+ & 260 & 793 & 13186 & 4.68e-08 & 9.99e-07 & 1.48e-07 & -2.8374853e+02 & 1.50e+03 \\ [3pt] 

G16& RNNAL & 42 & 57832 & 224 & 9.17e-07 & 9.77e-07 & 1.85e-11 & -2.8511897e+02 & 4.86e+02 \\ 
$n=  800$& SDPNAL+ & 237 & 813 & 14404 & 6.36e-07 & 8.93e-07 & 1.99e-06 & -2.8511755e+02 & 1.39e+03 \\ [3pt] 

G17& RNNAL & 30 & 43153 & 169 & 8.86e-07 & 5.69e-07 & 2.65e-11 & -2.8612382e+02 & 3.37e+02 \\ 
$n=  800$& SDPNAL+ & 300 & 732 & 19900 & 5.09e-07 & 1.85e-06 & 3.02e-06 & -2.8612513e+02 & 1.77e+03 \\ [3pt] 

G18& RNNAL & 15 & 2855 & 148 & 9.58e-07 & 2.38e-07 & 5.11e-10 & -2.7899999e+02 & 2.39e+01 \\ 
$n=  800$& SDPNAL+ & 131 & 372 & 6700 & 1.02e-06 & 6.92e-07 & 6.80e-08 & -2.7900027e+02 & 6.71e+02 \\ [3pt] 

G19& RNNAL & 27 & 31842 & 148 & 6.70e-07 & 6.39e-07 & 2.70e-11 & -2.8374869e+02 & 2.34e+02 \\ 
$n=  800$& SDPNAL+ & 260 & 793 & 13186 & 4.68e-08 & 9.99e-07 & 1.48e-07 & -2.8374853e+02 & 1.48e+03 \\ [3pt] 

G20& RNNAL & 42 & 57832 & 224 & 9.17e-07 & 9.77e-07 & 1.85e-11 & -2.8511897e+02 & 4.86e+02 \\ 
$n=  800$& SDPNAL+ & 237 & 813 & 14404 & 6.36e-07 & 8.93e-07 & 1.99e-06 & -2.8511755e+02 & 1.40e+03 \\ [3pt] 

G21& RNNAL & 30 & 43153 & 169 & 8.86e-07 & 5.69e-07 & 2.65e-11 & -2.8612382e+02 & 3.40e+02 \\ 
$n=  800$& SDPNAL+ & 300 & 732 & 19900 & 5.09e-07 & 1.85e-06 & 3.02e-06 & -2.8612513e+02 & 1.78e+03 \\ [3pt] 

G22& RNNAL & 10 & 1400 & 109 & 6.64e-07 & 2.29e-08 & 7.09e-08 & -5.7740156e+02 & 7.63e+01 \\ 
$n= 2000$& SDPNAL+ & 73 & 117 & 2080 & 3.30e-07 & 7.95e-07 & 3.17e-14 & -5.7740063e+02 & 1.13e+03 \\ [3pt] 

G23& RNNAL & 11 & 1388 & 110 & 2.63e-07 & 2.63e-08 & 1.65e-08 & -5.7655216e+02 & 8.29e+01 \\ 
$n= 2000$& SDPNAL+ & 54 & 88 & 1450 & 1.06e-06 & 7.61e-07 & 5.17e-08 & -5.7654759e+02 & 7.79e+02 \\ [3pt] 

G24& RNNAL & 10 & 1266 & 112 & 9.48e-07 & 2.41e-08 & 2.05e-07 & -5.7891540e+02 & 7.61e+01 \\ 
$n= 2000$& SDPNAL+ & 72 & 105 & 2070 & 3.86e-07 & 7.54e-07 & 1.53e-12 & -5.7891434e+02 & 1.11e+03 \\ [3pt] 

G25& RNNAL & 12 & 1313 & 104 & 5.82e-07 & 1.89e-08 & 5.15e-07 & -5.7704288e+02 & 8.04e+01 \\ 
$n= 2000$& SDPNAL+ & 66 & 92 & 1600 & 5.62e-07 & 5.08e-07 & 1.97e-12 & -5.7704237e+02 & 9.10e+02 \\ [3pt] 

G26& RNNAL & 17 & 1113 & 101 & 6.83e-07 & 8.67e-09 & 4.36e-07 & -5.7691745e+02 & 7.40e+01 \\ 
$n= 2000$& SDPNAL+ & 71 & 111 & 2060 & 2.18e-07 & 2.55e-07 & 2.10e-13 & -5.7691681e+02 & 1.17e+03 \\ [3pt] 

G27& RNNAL & 11 & 1319 & 108 & 3.57e-07 & 4.04e-08 & 1.13e-07 & -5.7740101e+02 & 7.97e+01 \\ 
$n= 2000$& SDPNAL+ & 73 & 117 & 2080 & 3.30e-07 & 7.97e-07 & 8.49e-13 & -5.7740063e+02 & 1.12e+03 \\ [3pt] 

G28& RNNAL & 11 & 1334 & 108 & 5.13e-07 & 2.60e-08 & 2.49e-07 & -5.7683206e+02 & 8.21e+01 \\ 
$n= 2000$& SDPNAL+ & 54 & 91 & 1450 & 1.83e-06 & 9.02e-07 & 1.92e-07 & -5.7682931e+02 & 7.96e+02 \\ [3pt] 

G29& RNNAL & 10 & 1266 & 112 & 9.48e-07 & 2.41e-08 & 2.05e-07 & -5.7891540e+02 & 7.65e+01 \\ 
$n= 2000$& SDPNAL+ & 72 & 114 & 2090 & 6.45e-07 & 9.80e-07 & 1.33e-12 & -5.7891477e+02 & 1.13e+03 \\ [3pt] 

G30& RNNAL & 12 & 1313 & 104 & 5.82e-07 & 1.89e-08 & 5.15e-07 & -5.7704288e+02 & 8.04e+01 \\ 
$n= 2000$& SDPNAL+ & 66 & 92 & 1600 & 5.62e-07 & 5.08e-07 & 1.97e-12 & -5.7704237e+02 & 9.28e+02 \\ [3pt] 

G31& RNNAL & 17 & 1113 & 101 & 6.83e-07 & 8.67e-09 & 4.36e-07 & -5.7691745e+02 & 7.38e+01 \\ 
$n= 2000$& SDPNAL+ & 71 & 111 & 2060 & 2.19e-07 & 2.55e-07 & 1.21e-12 & -5.7691681e+02 & 1.19e+03 \\ [3pt] 

G32& RNNAL & 13 & 7402 & 21 & 3.23e-07 & 2.53e-08 & 1.84e-07 & -9.9998584e+02 & 2.99e+02 \\ 
$n= 2000$& SDPNAL+ & - & - & - & - & - & - & - & - \\ [3pt] 

G33& RNNAL & 19 & 2381 & 36 & 4.76e-07 & 8.61e-07 & 5.30e-09 & -9.9604130e+02 & 1.17e+02 \\ 
$n= 2000$& SDPNAL+ & 60 & 127 & 1500 & 1.62e-07 & 8.70e-07 & 4.84e-13 & -9.9604179e+02 & 8.30e+02 \\ [3pt] 

G34& RNNAL & 10 & 2188 & 11 & 3.79e-07 & 1.48e-08 & 5.58e-07 & -9.9999198e+02 & 9.91e+01 \\ 
$n= 2000$& SDPNAL+ & - & - & - & - & - & - & - & - \\ [3pt] 

G35& RNNAL & 19 & 35486 & 367 & 7.25e-07 & 8.81e-07 & 5.58e-12 & -7.1823685e+02 & 2.19e+03 \\ 
$n= 2000$& SDPNAL+ & - & - & - & - & - & - & - & - \\ [3pt] 

G36& RNNAL & 19 & 33698 & 445 & 6.54e-07 & 8.15e-07 & 1.77e-11 & -6.9600062e+02 & 2.23e+03 \\ 
$n= 2000$& SDPNAL+ & - & - & - & - & - & - & - & - \\ [3pt] 

G37& RNNAL & 16 & 7089 & 416 & 7.11e-07 & 8.94e-07 & 1.77e-10 & -7.0800000e+02 & 4.71e+02 \\ 
$n= 2000$& SDPNAL+ & - & - & - & - & - & - & - & - \\ [3pt] 

G38& RNNAL & 17 & 12945 & 375 & 7.38e-07 & 8.23e-07 & 7.14e-12 & -7.1600032e+02 & 8.21e+02 \\ 
$n= 2000$& SDPNAL+ & - & - & - & - & - & - & - & - \\ [3pt] 

G39& RNNAL & 19 & 35486 & 367 & 7.25e-07 & 8.81e-07 & 5.58e-12 & -7.1823685e+02 & 2.20e+03 \\ 
$n= 2000$& SDPNAL+ & - & - & - & - & - & - & - & - \\ [3pt] 

G40& RNNAL & 19 & 33698 & 445 & 6.54e-07 & 8.15e-07 & 1.77e-11 & -6.9600062e+02 & 2.23e+03 \\ 
$n= 2000$& SDPNAL+ & - & - & - & - & - & - & - & - \\ [3pt] 

G41& RNNAL & 16 & 7089 & 416 & 7.11e-07 & 8.94e-07 & 1.77e-10 & -7.0800000e+02 & 4.72e+02 \\ 
$n= 2000$& SDPNAL+ & - & - & - & - & - & - & - & - \\ [3pt] 

G42& RNNAL & 17 & 12945 & 375 & 7.38e-07 & 8.23e-07 & 7.14e-12 & -7.1600032e+02 & 8.19e+02 \\ 
$n= 2000$& SDPNAL+ & - & - & - & - & - & - & - & - \\ [3pt] 

G43& RNNAL & 10 & 795 & 80 & 7.48e-07 & 4.18e-08 & 2.09e-08 & -2.7973625e+02 & 9.61e+00 \\ 
$n= 1000$& SDPNAL+ & 48 & 61 & 1250 & 4.58e-07 & 6.61e-07 & 1.28e-13 & -2.7973595e+02 & 1.43e+02 \\ [3pt] 

G44& RNNAL & 13 & 1045 & 83 & 9.54e-07 & 6.82e-08 & 2.02e-08 & -2.7974645e+02 & 1.35e+01 \\ 
$n= 1000$& SDPNAL+ & 47 & 62 & 1250 & 8.32e-07 & 9.03e-07 & 7.58e-14 & -2.7974580e+02 & 1.42e+02 \\ [3pt] 

G45& RNNAL & 13 & 1077 & 79 & 4.23e-07 & 2.82e-08 & 6.10e-08 & -2.7931767e+02 & 1.37e+01 \\ 
$n= 1000$& SDPNAL+ & 46 & 59 & 1250 & 2.33e-07 & 6.91e-07 & 8.09e-13 & -2.7931751e+02 & 1.42e+02 \\ [3pt] 

G46& RNNAL & 16 & 1095 & 76 & 4.32e-07 & 2.96e-08 & 1.36e-07 & -2.7903270e+02 & 1.43e+01 \\ 
$n= 1000$& SDPNAL+ & 54 & 74 & 1250 & 1.36e-06 & 8.99e-07 & 1.84e-07 & -2.7903228e+02 & 1.46e+02 \\ [3pt] 

G47& RNNAL & 13 & 1076 & 87 & 6.98e-07 & 4.94e-08 & 4.59e-07 & -2.8089197e+02 & 1.34e+01 \\ 
$n= 1000$& SDPNAL+ & 47 & 63 & 1250 & 2.96e-07 & 9.33e-07 & 1.54e-09 & -2.8089134e+02 & 1.42e+02 \\ [3pt] 

G48& RNNAL & 8 & 1497 & 21 & 9.25e-07 & 4.67e-08 & 3.09e-07 & -1.4999238e+03 & 2.44e+02 \\ 
$n= 3000$& SDPNAL+ & - & - & - & - & - & - & - & - \\ [3pt] 

G49& RNNAL & 12 & 11438 & 17 & 2.34e-07 & 1.36e-08 & 1.95e-07 & -1.4999873e+03 & 1.69e+03 \\ 
$n= 3000$& SDPNAL+ & - & - & - & - & - & - & - & - \\ [3pt] 

G50& RNNAL & 13 & 2007 & 71 & 5.16e-07 & 9.85e-07 & 7.40e-08 & -1.4940618e+03 & 3.65e+02 \\ 
$n= 3000$& SDPNAL+ & - & - & - & - & - & - & - & - \\ [3pt] 

G51& RNNAL & 21 & 11691 & 239 & 7.89e-07 & 6.89e-07 & 3.04e-11 & -3.4900026e+02 & 1.53e+02 \\ 
$n= 1000$& SDPNAL+ & 99 & 352 & 3500 & 8.01e-07 & 6.63e-07 & 8.62e-13 & -3.4899993e+02 & 7.07e+02 \\ [3pt] 

G52& RNNAL & 21 & 34697 & 223 & 7.24e-07 & 9.45e-07 & 4.54e-11 & -3.4838649e+02 & 4.38e+02 \\ 
$n= 1000$& SDPNAL+ & 140 & 506 & 6063 & 1.69e-07 & 9.99e-07 & 4.86e-07 & -3.4838599e+02 & 1.19e+03 \\ [3pt] 

G53& RNNAL & 22 & 75737 & 251 & 8.02e-07 & 8.86e-07 & 3.32e-10 & -3.4821204e+02 & 9.74e+02 \\ 
$n= 1000$& SDPNAL+ & 333 & 1287 & 12531 & 1.42e-06 & 9.99e-07 & 1.49e-07 & -3.4821136e+02 & 2.74e+03 \\ [3pt] 

G54& RNNAL & 19 & 4045 & 237 & 7.51e-07 & 8.13e-07 & 2.61e-11 & -3.4100004e+02 & 5.85e+01 \\ 
$n= 1000$& SDPNAL+ & 130 & 480 & 5245 & 1.59e-07 & 9.88e-07 & 1.31e-06 & -3.4100088e+02 & 1.08e+03 \\ [3pt] 

G55& RNNAL & 20 & 2130 & 353 & 9.96e-07 & 1.85e-07 & 4.84e-08 & -2.3230485e+03 & 1.45e+03 \\ 
$n= 5000$& SDPNAL+ & - & - & - & - & - & - & - & - \\ [3pt] 

G56& RNNAL & 20 & 2130 & 353 & 9.96e-07 & 1.85e-07 & 4.84e-08 & -2.3230485e+03 & 1.44e+03 \\ 
$n= 5000$& SDPNAL+ & - & - & - & - & - & - & - & - \\ [3pt] 

G57& RNNAL & 15 & 4640 & 13 & 3.87e-07 & 1.14e-08 & 9.20e-07 & -2.4999927e+03 & 2.15e+03 \\ 
$n= 5000$& SDPNAL+ & - & - & - & - & - & - & - & - \\ [3pt]

\end{longtable}
\end{tiny}

\subsection{Experiments on DQKP problems}\label{appendix-DQKP}
\begin{tiny}
\begin{longtable}[c]{clrrrccccc}
\caption{Computational results for DQKP problems. \label{tab-appendix-DQKP}} \\
\toprule
 problem & algorithm & it & itsub & $r$/itA & $\operatorname{R_p}$ & $\operatorname{R_d}$ & $\operatorname{R_c}$ & obj & time   \\
\midrule
\endfirsthead

\multicolumn{10}{c}%
{{ Table \thetable\ continued from previous page}} \\
\toprule
problem & algorithm & it & itsub & $r$/itA & $\operatorname{R_p}$ & $\operatorname{R_d}$ & $\operatorname{R_c}$ & obj & time   \\
\midrule
\endhead
\midrule
\multicolumn{10}{r}{{Continued on next page}} \\
\midrule
\endfoot
\bottomrule
\endlastfoot

$n=1000$ & RNNAL & 13 & 268 & 31 & 7.70e-07 & 4.50e-07 & 3.20e-09 & -2.1640571e+07 & 8.50e+00 \\ 
$d=0.1$ & SDPNAL+ & 114 & 320 & 6700 & 3.78e-07 & 9.99e-07 & 2.98e-13 & -2.1640569e+07 & 1.22e+03 \\ [5pt] 

$n=1000$ & RNNAL & 17 & 398 & 14 & 6.13e-07 & 7.05e-07 & 2.06e-09 & -6.5085331e+06 & 9.46e+00 \\ 
$d=0.5$ & SDPNAL+ & 74 & 532 & 4253 & 9.81e-07 & 1.10e-07 & 6.48e-14 & -6.5085176e+06 & 1.53e+03 \\ [5pt] 

$n=1000$ & RNNAL & 18 & 628 & 33 & 5.92e-07 & 1.00e-08 & 2.88e-08 & -4.4788736e+06 & 1.55e+01 \\ 
$d=1.0$ & SDPNAL+ & - & - & - & - & - & - & - & - \\ [5pt] 

$n=2000$ & RNNAL & 13 & 242 & 55 & 8.89e-07 & 1.45e-07 & 5.34e-09 & -8.6175006e+07 & 3.15e+01 \\ 
$d=0.1$ & SDPNAL+ & - & - & - & - & - & - & - & - \\ [5pt] 

$n=2000$ & RNNAL & 23 & 750 & 56 & 2.17e-07 & 8.39e-07 & 4.80e-08 & -2.6600038e+07 & 7.39e+01 \\ 
$d=0.5$ & SDPNAL+ & - & - & - & - & - & - & - & - \\ [5pt] 

$n=2000$ & RNNAL & 14 & 708 & 33 & 7.44e-07 & 7.21e-07 & 2.23e-09 & -1.8362081e+07 & 6.07e+01 \\ 
$d=1.0$ & SDPNAL+ & - & - & - & - & - & - & - & - \\ [5pt] 

$n=5000$ & RNNAL & 54 & 1065 & 127 & 3.18e-07 & 8.15e-07 & 3.35e-09 & -3.6381731e+08 & 1.58e+03 \\ 
$d=0.1$ & SDPNAL+ & - & - & - & - & - & - & - & - \\ [5pt] 

$n=5000$ & RNNAL & 63 & 2534 & 147 & 1.18e-08 & 3.97e-07 & 3.35e-08 & -1.6001202e+08 & 2.40e+03 \\ 
$d=0.5$ & SDPNAL+ & - & - & - & - & - & - & - & - \\ [5pt] 

$n=5000$ & RNNAL & 39 & 1350 & 144 & 3.66e-08 & 5.44e-07 & 7.25e-08 & -1.0264044e+08 & 1.40e+03 \\ 
$d=1.0$ & SDPNAL+ & - & - & - & - & - & - & - & - \\

\end{longtable}
\end{tiny}

\subsection{Experiments on GWD shape correspondence problems}\label{appendix-GWD}

\begin{tiny}
\begin{longtable}[c]{clrrrccccc}
\caption{Computational results for GWD shape correspondence problems. \label{tab-GWD-appendix-shape}} \\
\toprule
 problem & algorithm & it & itsub & $r$/itA & $\operatorname{R_p}$ & $\operatorname{R_d}$ & $\operatorname{R_c}$ & obj & time   \\
\midrule
\endfirsthead

\multicolumn{10}{c}%
{{ Table \thetable\ continued from previous page}} \\
\toprule
problem & algorithm & it & itsub & $r$/itA & $\operatorname{R_p}$ & $\operatorname{R_d}$ & $\operatorname{R_c}$ & obj & time   \\
\midrule
\endhead
\midrule
\multicolumn{10}{r}{{Continued on next page}} \\
\midrule
\endfoot
\bottomrule
\endlastfoot

Cat& RNNAL & 67 & 21884 & 135 & 9.85e-07 & 7.84e-07 & 1.65e-07 & 1.3170169e+05 & 3.17e+02 \\ 
$n=900$ & SDPNAL+ & 301 & 455 & 15122 & 3.82e-07 & 8.60e-07 & 1.37e-07 & 1.3170079e+05 & 1.95e+03 \\ [5pt] 

David& RNNAL & 46 & 8844 & 5 & 9.21e-07 & 4.10e-08 & 3.91e-07 & 1.2642888e+05 & 8.85e+01 \\ 
$n=900$ & SDPNAL+ & 255 & 429 & 5693 & 1.78e-08 & 9.31e-07 & 5.60e-14 & 1.2643058e+05 & 7.98e+02 \\ [5pt] 

Dog& RNNAL & 53 & 10808 & 8 & 3.57e-07 & 4.80e-08 & 3.99e-07 & 8.8541417e+04 & 1.03e+02 \\ 
$n=900$ & SDPNAL+ & 265 & 362 & 6851 & 3.44e-07 & 9.48e-07 & 4.97e-07 & 8.8543446e+04 & 9.13e+02 \\ [5pt] 

Gorilla& RNNAL & 35 & 7527 & 6 & 8.01e-07 & 8.73e-08 & 5.80e-08 & 1.4314155e+05 & 8.21e+01 \\ 
$n=900$ & SDPNAL+ & 198 & 324 & 5905 & 5.40e-09 & 9.59e-08 & 4.08e-15 & 1.4356859e+05 & 7.83e+02 \\ [5pt] 

Seahorse& RNNAL & 42 & 9266 & 6 & 4.05e-07 & 9.40e-08 & 7.36e-07 & 3.4701539e+05 & 8.65e+01 \\ 
$n=900$ & SDPNAL+ & 431 & 681 & 9004 & 7.05e-07 & 7.98e-07 & 1.83e-14 & 3.4703056e+05 & 1.28e+03 \\ [5pt] 

Cat& RNNAL & 35 & 6448 & 27 & 8.45e-07 & 2.24e-07 & 7.00e-07 & 1.2788956e+05 & 2.80e+02 \\ 
$n=1600$ & SDPNAL+ & - & - & - & - & - & - & - & - \\ [5pt] 

David& RNNAL & 131 & 27042 & 6 & 3.33e-07 & 5.50e-07 & 4.26e-07 & 1.9694194e+05 & 8.68e+02 \\ 
$n=1600$ & SDPNAL+ & - & - & - & - & - & - & - & - \\ [5pt] 

Dog& RNNAL & 80 & 15312 & 7 & 1.11e-07 & 1.19e-07 & 5.32e-08 & 1.2109773e+05 & 5.22e+02 \\ 
$n=1600$ & SDPNAL+ & - & - & - & - & - & - & - & - \\ [5pt] 

Gorilla& RNNAL & 66 & 12705 & 22 & 7.15e-07 & 2.06e-07 & 1.46e-07 & 4.2828956e+05 & 4.88e+02 \\ 
$n=1600$ & SDPNAL+ & - & - & - & - & - & - & - & - \\ [5pt] 

Seahorse& RNNAL & 150 & 39803 & 28 & 5.40e-07 & 9.00e-07 & 6.18e-07 & 5.0244394e+05 & 1.24e+03 \\ 
$n=1600$ & SDPNAL+ & - & - & - & - & - & - & - & - \\ [5pt] 

Cat& RNNAL & 53 & 15128 & 23 & 6.19e-07 & 8.29e-07 & 5.37e-07 & 2.3319179e+05 & 2.16e+03 \\ 
$n=2500$ & SDPNAL+ & - & - & - & - & - & - & - & - \\

\end{longtable}
\end{tiny}

\subsection{Experiments on graph partition problems}\label{appendix-GWD-partition}
\normalsize

Graph partition aims to match the source graph having $l$ nodes with a disconnected target graph having $k$ isolated and self-connected super nodes, where $k$ is the number of partitions. We choose synthetic datasets similar to the procedure in \cite{xu2019scalable}. Specifically, the source graph is a Gaussian random partition graph with $l$ nodes and $k$ partitions. The size of each cluster is drawn from a normal distribution $\mathcal{N}(l/k, l/100)$. The size of the last cluster is adjusted to make the total number of nodes equal to $l$. The nodes are connected within the partitions 
with the probability of 0.9 and between partitions with the probability of 0.1. $D_X$ and $D_Y$ are the adjacency matrices of the source graph and target graph, respectively. We set the distribution $a$ to be the normalized cluster size of the target graph and $b$ to be the empirical distribution of the source graph. We choose the partition number $k=3$ and the number of nodes $l\in\{300,600,900,1200,1500\}$.

\begin{tiny}
\begin{longtable}[c]{clrrrccccc}
\caption{Computational results for GWD graph partition problems. \label{tab-GWD-partition}} \\
\toprule
 problem & algorithm & it & itsub & $r$/itA & $\operatorname{R_p}$ & $\operatorname{R_d}$ & $\operatorname{R_c}$ & obj & time   \\
\midrule
\endfirsthead

\multicolumn{10}{c}%
{{ Table \thetable\ continued from previous page}} \\
\toprule
problem & algorithm & it & itsub & $r$/itA & $\operatorname{R_p}$ & $\operatorname{R_d}$ & $\operatorname{R_c}$ & obj & time   \\
\midrule
\endhead
\midrule
\multicolumn{10}{r}{{Continued on next page}} \\
\midrule
\endfoot
\bottomrule
\endlastfoot

& RNNAL & 53 & 1507 & 95 & 9.29e-07 & 1.75e-08 & 8.08e-10 & 2.0125132e-01 & 2.31e+01 \\ 
$n=900$ & SDPNAL+ & 77 & 293 & 4011 & 1.02e-06 & 6.13e-08 & 3.93e-12 & 2.0127614e-01 & 1.17e+03 \\ 
& RNNAL-Diag & 76 & 52654 & 5 & 8.52e-07 & 5.60e-08 & 2.85e-19 & 2.0126906e-01 & 5.49e+02 \\ [3pt] 

& RNNAL & 77 & 1806 & 116 & 5.86e-07 & 2.13e-09 & 3.17e-11 & 2.0324038e-01 & 1.51e+02 \\ 
$n=1800$ & SDPNAL+ & - & - & - & - & - & - & - & - \\ 
& RNNAL-Diag & - & - & - & - & - & - & - & - \\ [3pt] 

& RNNAL & 89 & 2075 & 177 & 6.57e-07 & 4.66e-07 & 9.65e-12 & 2.0248959e-01 & 5.21e+02 \\ 
$n=2700$ & SDPNAL+ & - & - & - & - & - & - & - & - \\ 
& RNNAL-Diag & - & - & - & - & - & - & - & - \\ [3pt] 

& RNNAL & 105 & 2412 & 233 & 3.58e-07 & 3.40e-07 & 7.69e-11 & 2.0344995e-01 & 1.65e+03 \\ 
$n=3600$ & SDPNAL+ & - & - & - & - & - & - & - & - \\ 
& RNNAL-Diag & - & - & - & - & - & - & - & - \\ [3pt] 

& RNNAL & 112 & 3042 & 304 & 2.49e-07 & 8.17e-07 & 1.35e-11 & 2.0350253e-01 & 3.51e+03 \\ 
$n=4500$ & SDPNAL+ & - & - & - & - & - & - & - & - \\ 
& RNNAL-Diag & - & - & - & - & - & - & - & - \\

\end{longtable}
\end{tiny}

Table \eqref{tab-GWD-partition} shows the numerical results on the GWD problems. Observe that RNNAL achieves the required accuracy within the 1 hour limit, while SDPNAL+ and RNNAL-Diag fail for problems with dimensions $n\geq 1800$. RNNAL is nearly 50 and 20 times faster than SDPNAL+ and RNNAL-Diag, respectively, as seen in the case of $n=900$. RNNAL-Diag is slow mainly because the number of ALM subproblems is usually much larger than that of RNNAL. Note that the final rank of the solutions obtained by RNNAL may be larger than the smallest rank, such as when $n=900$. This occurs because RNNAL adjusts the rank adaptively. This approach assists in escaping saddle points and balancing the tradeoff between subproblem iteration count and computational cost per iteration.

\bibliographystyle{abbrv}
\bibliography{RNNAL}

\end{document}